%

\def\hefresh{\setminus}


%
%
\catcode`\@=11
\message{*** Israel Journal of Mathematics package, Version 0.9, %
October 1995 ***}
%

\newbox\abstractbox
\newbox\authorbox
\newbox\secbox
\newbox\recdatebox
\setbox\strutbox\hbox{\vrule height9.5pt depth3.5pt width0pt}

\newcount\firstpage 
\newcount\lastpage 
\newskip\sectionskipamount

\newif\ifextra \extratrue

\newdimen\itemindent \itemindent=25pt
\newdimen\oldindent
\newdimen\newvsize

%
\font\twelverm=cmr10 at 12pt
\font\tenrm=cmr10
\font\ninerm=cmr9
\font\eightrm=cmr8
\font\sevenrm=cmr7
\font\sixrm=cmr6
\font\fiverm=cmr5
%
\font\twelvesl=cmsl10 at 12pt
\font\tensl=cmsl10
\font\ninesl=cmsl9
\font\eightsl=cmsl8
%
\font\twelveit=cmti10 at 12pt
\font\tenit=cmti10
\font\nineit=cmti9
\font\eightit=cmti8

%
\font\twelvebf=cmbx10 at 12pt
\font\tenbf=cmbx10
\font\ninebf=cmbx9
\font\eightbf=cmbx8
\font\sevenbf=cmbx7
\font\sixbf=cmbx6
\font\fivebf=cmbx5
%
\font\tensmc=cmcsc10
\font\ninesmc=cmcsc9
\font\eightsmc=cmcsc8
%
\font\twelvemi=cmmi10 at 12pt
\font\tenmi=cmmi10
\font\ninemi=cmmi9
\font\eightmi=cmmi8
\font\sevenmi=cmmi7
\font\sixmi=cmmi6
\font\fivemi=cmmi5
%
\font\twelvesy=cmsy10 at 12pt
\font\tensy=cmsy10
\font\ninesy=cmsy9
\font\eightsy=cmsy8
\font\sevensy=cmsy7
\font\sixsy=cmsy6
\font\fivesy=cmsy5
%
\font\tenex=cmex10
\font\nineex=cmex9
%
\font\twelvemsb=msbm10 at 12pt
\font\tenmsb=msbm10
\font\ninemsb=msbm9
\font\eightmsb=msbm8
\font\sevenmsb=msbm7
\font\sixmsb=msbm6
\font\fivemsb=msbm5
%

%
\font\twelveeufm=eufm10 at 12pt
\font\teneufm=eufm10
\font\nineeufm=eufm9
\font\eighteufm=eufm8
\font\seveneufm=eufm7
\font\sixeufm=eufm6
\font\fiveeufm=eufm5

\def\twelvepoint{%
\baselineskip 18pt%
\def\rm{\fam0\twelverm}%
\textfont0=\twelverm \scriptfont0=\ninerm \scriptscriptfont0=\sixrm%
\textfont1=\twelvemi \scriptfont1=\ninemi \scriptscriptfont1=\sixmi%
\textfont2=\twelvesy \scriptfont2=\ninesy \scriptscriptfont2=\sixsy%
\textfont3=\tenex \scriptfont3=\nineex \scriptscriptfont3=\nineex%
\textfont\itfam=\twelveit \def\it{\fam\itfam\twelveit}%
\textfont\slfam=\twelvesl \def\sl{\fam\slfam\twelvesl}%
\textfont\bffam=\twelvebf \scriptfont\bffam=\ninebf%
\scriptscriptfont\bffam=\sixbf \def\bf{\fam\bffam\twelvebf}%
\textfont\msbfam=\twelvemsb \scriptfont\msbfam=\ninemsb%
\scriptscriptfont\msbfam=\sixmsb \def\Bbb{\fam\msbfam\twelvemsb}%
\textfont\eufmfam=\twelveeufm \scriptfont\eufmfam=\nineeufm%
\scriptscriptfont\eufmfam=\sixeufm \def\frak{\fam\eufmfam\twelveeufm}%
\rm}

\def\tenpoint{%
\baselineskip 15pt%
\def\rm{\fam0\tenrm}%
\textfont0=\tenrm \scriptfont0=\sevenrm \scriptscriptfont0=\fiverm%
\textfont1=\tenmi \scriptfont1=\sevenmi \scriptscriptfont1=\fivemi%
\textfont2=\tensy \scriptfont2=\sevensy \scriptscriptfont2=\fivesy%
\textfont3=\tenex \scriptfont3=\nineex \scriptscriptfont3=\nineex%
\textfont\itfam=\tenit \def\it{\fam\itfam\tenit}%
\textfont\slfam=\tensl \def\sl{\fam\slfam\tensl}%
\textfont\bffam=\tenbf \scriptfont\bffam=\sevenbf%
\scriptscriptfont\bffam=\fivebf \def\bf{\fam\bffam\tenbf}%
\textfont\msbfam=\tenmsb \scriptfont\msbfam=\sevenmsb%
\scriptscriptfont\msbfam=\fivemsb \def\Bbb{\fam\msbfam\tenmsb}%
\textfont\eufmfam=\teneufm \scriptfont\eufmfam=\seveneufm%
\scriptscriptfont\eufmfam=\fiveeufm \def\frak{\fam\eufmfam\teneufm}%
\rm}

\def\ninepoint{%
\baselineskip 13.5pt%
\def\rm{\fam0\ninerm}%
\textfont0=\ninerm \scriptfont0=\sixrm \scriptscriptfont0=\fiverm%
\textfont1=\ninemi \scriptfont1=\sixmi \scriptscriptfont1=\fivemi%
\textfont2=\ninesy \scriptfont2=\sixsy \scriptscriptfont2=\fivesy%
\textfont3=\nineex \scriptfont3=\nineex \scriptscriptfont3=\nineex%
\textfont\itfam=\nineit \def\it{\fam\itfam\nineit}%
\textfont\slfam=\ninesl \def\sl{\fam\slfam\ninesl}%
\textfont\bffam=\ninebf \scriptfont\bffam=\sixbf%
\scriptscriptfont\bffam=\fivebf \def\bf{\fam\bffam\ninebf}%
\textfont\msbfam=\ninemsb \scriptfont\msbfam=\sixmsb%
\scriptscriptfont\msbfam=\fivemsb \def\Bbb{\fam\msbfam\ninemsb}%
\textfont\eufmfam=\nineeufm \scriptfont\eufmfam=\sixeufm%
\scriptscriptfont\eufmfam=\fiveeufm \def\frak{\fam\eufmfam\nineeufm}%
\def\smc{\ninesmc}%
\rm}

\def\eightpoint{%
\baselineskip 12pt%
\def\rm{\fam0\eightrm}%
\textfont0=\eightrm \scriptfont0=\sixrm \scriptscriptfont0=\fiverm%
\textfont1=\eightmi \scriptfont1=\sixmi \scriptscriptfont1=\fivemi%
\textfont2=\eightsy \scriptfont2=\sixsy \scriptscriptfont2=\fivesy%
\textfont3=\nineex \scriptfont3=\nineex \scriptscriptfont3=\nineex%
\textfont\itfam=\eightit \def\it{\fam\itfam\eightit}%
\textfont\slfam=\eightsl \def\sl{\fam\slfam\eightsl}%
\textfont\bffam=\eightbf \scriptfont\bffam=\sixbf%
\scriptscriptfont\bffam=\fivebf \def\bf{\fam\bffam\eightbf}%
\textfont\msbfam=\eightmsb \scriptfont\msbfam=\sixmsb%
\scriptscriptfont\msbfam=\fivemsb \def\Bbb{\fam\msbfam\eightmsb}%
\textfont\eufmfam=\eighteufm \scriptfont\eufmfam=\sixeufm%
\scriptscriptfont\eufmfam=\fiveeufm \def\frak{\fam\eufmfam\eighteufm}%
\def\smc{\eightsmc}%
\rm}

\def\sixpoint{%
\def\rm{\fam0\sixrm}%
\textfont0=\sixrm \scriptfont0=\sixrm \scriptscriptfont0=\fiverm%
\textfont1=\sixmi \scriptfont1=\sixmi \scriptscriptfont1=\fivemi%
\textfont2=\sixsy \scriptfont2=\sixsy \scriptscriptfont2=\fivesy%
\textfont\bffam=\sixbf \scriptfont\bffam=\sixbf%
\scriptscriptfont\bffam=\fivebf \def\bf{\fam\bffam\sixbf}%
\textfont\msbfam=\sixmsb \scriptfont\msbfam=\sixmsb%
\scriptscriptfont\msbfam=\fivemsb \def\Bbb{\fam\msbfam\sixmsb}%
\textfont\eufmfam=\sixeufm \scriptfont\eufmfam=\sixeufm%
\scriptscriptfont\eufmfam=\fiveeufm \def\frak{\fam\eufmfam\sixeufm}%
\rm}

\def\titlefont{\twelvepoint\rm
   \baselineskip=15pt\lineskiplimit=2pt\lineskip=2pt}
\def\authorfont{\tensmc}
\def\affiliationfont{\eightpoint\it}
\def\abstractfont{\eightpoint\rm}
\def\headlinefont{\eightpoint\rm}
\def\smc{\tensmc}


\def\colon{{:}\;}
\def\setminus{
   \def\hefreshD{\mathop{\raise1.5pt\hbox{${\smallsetminus}$}}}
   \def\hefreshS{\mathop{\raise0.85pt\hbox{$\scriptstyle\smallsetminus$}}}
   \mathchoice{\hefreshD}{\hefreshD}{\hefreshS}{\hefreshS}}

\def\empty{}
\def\ifempty#1{\ifx\empty#1\empty\relax}
\def\newline{\hfill\break}
\def\setminus{\mathop{\raise1.5pt\hbox{${\smallsetminus}$}}}
\def\shvor{\discretionary{}{}{}}

%
%
\def\item#1{\oldindent=\parindent\parindent=\itemindent%
\par\hangindent\parindent\indent\llap{\rm #1\enspace}%
\parindent=\oldindent\ignorespaces}
\def\itemitem#1{\oldindent=\parindent\parindent=\itemindent%
\par\indent \hangindent2\parindent%
\indent\llap{\rm #1\enspace}\parindent=\itemindent\ignorespaces}
%
\def\condition#1{\item{\rm#1}}
\def\noextra{\extrafalse}

\parindent=1em
\baselineskip 15pt
\hsize=12.3 cm
\def\normalvsize{18.5 cm }

\def\ijvol#1{\def\thevol{#1\unskip}}
\def\ijyear#1{\def\theyear{#1\unskip}}
\def\recdate#1{\setbox\recdatebox\vbox{\ninepoint
\parindent=\footnoteindent
\hangindent\footnoteindent\indent\strut\ignorespaces #1}
  \newvsize=\normalvsize
  \advance\newvsize by -\ht\recdatebox
  \advance\newvsize by -\dp\recdatebox
  \global\vsize=\newvsize
}
\def\shortauthor#1{\def\theshortauthor{#1\unskip}}
\def\shorttitle#1{\def\theshorttitle{#1\unskip}}
\def\rightheadline{\rm Vol.\ \thevol, \theyear\hfill
   \theshorttitle\hfill\the\pageno}
\def\leftheadline{\rm\the\pageno\hfil%
   \uppercase{\theshortauthor}\hfil Isr. J. Math.}
\def\firsthead{\sixpoint\rm \hfil ISRAEL JOURNAL OF MATHEMATICS
   {\sixbf\thevol} (\theyear), \the\firstpage--\the\lastpage\hfil}
%
%
\outer\def\title#1{\def\thetitle{\vglue 36pt plus 12pt minus 12pt
   {\let\\=\cr\tabskip\centering\titlefont
   \halign to \hsize{\hfill\ignorespaces{####}\unskip\hfill\cr
   #1\crcr}}}}

\def\english{english}

\def\language{english}
\def\BY{\ifx\language\english BY\else PAR\fi}
\def\AND{\ifx\language\english AND\else ET\fi}

\outer\def\author#1{\setbox\authorbox=
    \vbox{\vbox{\vglue 21pt plus 6pt minus 3pt}
    \centerline{\sixrm \BY}%
    \vglue 3pt
    {\let\\=\cr\tabskip\centering\authorfont
    \halign to \hsize{\hfill\ignorespaces##\unskip\hfill\cr
    #1\crcr}}}}

\outer\def\affiliation#1{
   \ifempty#1\else\setbox\authorbox
   =\vbox{\unvbox\authorbox
   \vbox{\vglue 4pt plus 1pt minus 1pt}
   {\let\\=\cr\tabskip\centering\affiliationfont
   \halign to\hsize
   {\hfill\ignorespaces##\unskip\hfill\cr #1\crcr}}
    \vglue 21pt plus 6pt minus 3pt\goodbreak}
   \fi}

\outer\def\addauthor#1{\ifempty#1\else
    \setbox\authorbox
    =\vbox{\unvbox\authorbox
    \centerline{\sixrm \AND}%
    \vglue 3pt
    {\let\\=\cr\tabskip\centering
    \halign to \hsize{\authorfont\hfill\ignorespaces##\unskip\hfill\cr
    #1\crcr}}}\fi}

\outer\def\beginabstract{\setbox\abstractbox=
\vbox\bgroup
   \centerline{{\sixrm ABSTRACT}} 
    \nobreak
   \bgroup\baselineskip=12pt\parindent=38pt\narrower
   \parindent=1em
    \abstractfont\noindent}
\def\endabstract{\smallskip\egroup\vskip 12pt plus 6pt minus 6pt
   \egroup}

\def\maketopmatter{%
  \pageno=\firstpage\relax%
  \headline={\ifnum\pageno>\firstpage{\headlinefont\ifodd\pageno\rightheadline\else\leftheadline\fi}\else\firsthead\fi}
  \footline={\ifnum\pageno>\firstpage{}\else{\hfil{\rm\the\pageno}\hfil}\fi}
  \thetitle
  \unvbox\authorbox
  \unvbox\abstractbox
  \bigskip\tenpoint}

\def\footnoteindent{20pt}
\newinsert\footins
\skip\footins=12pt plus 2pt minus 4pt
\dimen\footins=30pc
\count\footins=1000
\def\footnoterule{\kern-3pt\hrule width 2in\kern 2.6pt}
\def\footnote#1#2{\unskip\edef\@sf{\spacefactor\the\spacefactor}%
{$^{#1}$}\@sf
\insert\footins{\parindent=\footnoteindent\ninepoint
\itemindent=\footnoteindent
\interlinepenalty=100 \let\par=\endgraf
\leftskip=0pt \rightskip=0pt
\baselineskip=11pt 
\splittopskip=12pt plus 1pt minus 1pt \floatingpenalty=20000
\smallskip\item{#1}{#2}}}

\def\pagecontents{%
\ifvoid \topins \else \unvbox \topins \fi
\dimen@ =\dp \@cclv \unvbox \@cclv 
\ifnum\pageno=\firstpage
\vskip \skip \footins \footnoterule \unvbox \footins 
\unvbox\recdatebox
\else \ifvoid \footins
\else \vskip \skip \footins \footnoterule \unvbox \footins \fi
 \ifr@ggedbottom \kern -\dimen@ \vfil \fi \fi}

\output={\ijoutput}
\def\ijoutput{%
\ifnum\pageno=\firstpage \global\vsize=\normalvsize\fi
\shipout\vbox{\makeheadline \pagebody \makefootline}
\advancepageno\ifnum\outputpenalty>-20000 \else\dosupereject\fi}

\sectionskipamount=20 pt plus 4 pt minus 4 pt
\def\sectionbreak{\par\ifdim\lastskip<\sectionskipamount
\removelastskip\vskip\sectionskipamount \fi\penalty-200}
%
\def\head#1 #2\endhead\par{\sectionbreak
   \message{#1 #2}
   \setbox\secbox=\hbox{\bf#1\enspace}\ignorespaces
   \vbox{\hangindent=\wd\secbox\noindent
   \hangafter=1\box\secbox{\bf\ignorespaces #2}\smallskip}\ignorespaces
   \nobreak\vskip-\parskip\noindent}

\def\subhead#1\endsubhead{\proclaim@{\smc#1\unskip\ifextra.\fi%
\extratrue}{}{\quad\rm}}
\def\subsubhead#1\endsubsubhead{\proclaim@{\it#1\unskip\ifextra.\fi%
\extratrue}{}{\quad\rm}}

\def\introduction\par{\sectionbreak\noindent{\bf Introduction}
\par\noindent}

%
\def\proclaim@#1#2#3{\ifdim\lastskip<\medskipamount\removelastskip
\medskip\fi\penalty-55\noindent#1#2#3\ignorespaces}

%
\def\proclaim #1{\proclaim@{\smc #1:\enspace}{}\sl}

%
\def\proclaiminfo #1#2{\proclaim@{\smc #1\enspace}{\rm(#2):}%
{\enspace\sl}}

\def\endproclaim{\par\rm
 \ifdim\lastskip<\medskipamount\removelastskip\penalty55\medskip\fi}

%
\def\demo#1{\removelastskip\medskip\par\noindent{\sl #1:}\enspace}
\def\demoinfo#1#2{\removelastskip\medskip\par%
\noindent{\sl #1}{\enspace\rm (#2):}}

\def\bull{\vrule height 6pt width 4pt depth 0.47pt} 
\def\qed{\qquad\bull\vadjust{\medskip}\ifmmode{}\else\par\fi}

%
\def\rem#1{\proclaim@{\sl#1:}{}{\quad\rm}}

%
\def\reminfo#1#2{\proclaim@{\sl#1:}{\enspace\sl#2.}{\quad\rm}}


\def\subdemo#1{\proclaim@{\smc#1\unskip\ifextra:\fi\extratrue}{}%
{\quad\rm}}

%
\def\subdemoinfo#1#2{\proclaim@{\smc#1\ifextra:\fi\extratrue}%
{\enspace\sl #2.}{\quad\rm}}

\def\acknowledgment#1\endacknowledgment%
{\proclaim@{\smc Acknowledgement:\quad}{}{\rm#1}}
\def\thank{\proclaim@{\smc Acknowledgement:\quad}{}{\rm}}
\def\thanks{\thank}

\newbox\keybox
\newdimen\keywidth
\def\references#1#2{
\def\setkeywidth##1{\setbox\keybox=\hbox{[##1]\enspace}
\keywidth=\wd\keybox}
\def\numberkey{\def\makekeybox####1{\setbox\keybox=\hbox
to\keywidth{\hfil[####1\unskip]\enspace}}}
\def\letterkey{\def\makekeybox####1{\setbox\keybox=\hbox
to\keywidth{[####1\unskip]\hfil}}}
\def\endreferences{\filbreak\tenpoint\nonfrenchspacing}
\csname#1\endcsname
\setkeywidth{#2}
\ninepoint\frenchspacing\rm
   \ifdim\lastskip<\bigskipamount\removelastskip\bigskip\fi
   \filbreak\centerline{\bf References}\medskip\nobreak
\def\book{\errmessage{\string\book command misplaced}}
\def\by{\errmessage{\string\by command misplaced}}
\def\editor{\errmessage{\string\editor command misplaced}}
\def\inbook{\errmessage{\string\inbook command misplaced}}
\def\journal{\errmessage{\string\journal command misplaced}}
\def\pages{\errmessage{\string\pages command misplaced}}
\def\seriesandpublisher{\errmessage{\string\seriesandpublisher
command misplaced}}
\def\reftitle{\errmessage{\string\reftitle command misplaced}}
\def\yr{\errmessage{\string\yr command misplaced}} 

\def\refhead##1##2##3{\par\smallbreak
\hangindent=\keywidth\makekeybox{##1}
\noindent\hangafter1\box\keybox\ignorespaces%
{\rm##2\unskip\ifextra, \fi\extratrue}%
{\sl##3\unskip\ifextra, \fi\extratrue}}
\def\endref{\errmessage{Unmatched \string\endref}}
\def\refpaper\key##1\by##2\reftitle##3\journal##4\volume##5\yr##6\pages%
##7\endref{\refhead{##1}{##2}{##3}%
{\rm ##4\unskip\ }%
{\bf ##5\unskip\ }%
{\ifx\noextra##6\else(##6\unskip), \fi\extratrue}%
{##7\unskip\ifextra.\fi\extratrue}%
\filbreak}

\def\refbook\key##1\by##2\reftitle##3\publisher##4\yr##5\endref
{\refhead{##1}{##2}{##3}%
{\ifx\noextra##4\else\rm ##4\unskip, \fi\extratrue}%
{\ifx\noextra##5\else\rm ##5\unskip.\fi\extratrue}\filbreak}

\def\refall\key##1\by##2\reftitle##3\rest##4\endref%
{\refhead{##1}{##2}{##3}{##4}\filbreak}
\extratrue}

\catcode`\@=\active


\input amssym.def
\input amssym.tex 

%

\def\lam{\lambda}		
		


\def\calE{{\cal E}}

\def\calJ{{\cal J}}

\def\calP{{\cal P}}

\def\calS{{\cal S}}
\def\calT{{\cal T}}

\def\calY{{\cal Y}}






\def\fr#1{{\frak #1}}
\def\frA{\fr{A}}	\def\fra{\fr{a}}
\def\frB{\fr{B}}	\def\frb{\fr{b}}
	\def\frc{\fr{c}}
	\def\frd{\fr{d}}
	\def\fre{\fr{e}}



\def\case#1:#2.{\par\medskip \noindent {\smc Case #1:} {\sl #2.}\quad}
\def\cf{}
\def\claim#1:#2.{\par\medskip \noindent {\smc Claim #1:} {\sl #2.}\quad}
\def\Claim:#1.{\par\medskip\noindent{\smc Claim:} {\sl #1.}\quad}
\def\colon{{:}\;}
\def\comment#1\endcomment{}
\def\condition(#1){\item{\rm#1}}

\def\dotcup{
 \def\dotcupD{\cup\hskip-7.5pt\cdot\hskip4.5pt}
 \def\dotcupS{\cup\hskip-4pt\raise0.6pt\hbox{$\cdot$}}
 \mathchoice{\dotcupD}{\dotcupD}{\dotcupS}{}}
\def\dotunion{
\def\dotunionD{\bigcup\kern-9pt\cdot\kern5pt}
\def\dotunionT{\bigcup\kern-7.5pt\cdot\kern3.5pt}
\mathop{\mathchoice{\dotunionD}{\dotunionT}{}{}}}

\def\finishproclaim {\par\rm
 \ifdim\lastskip<\medskipamount\removelastskip\penalty55\medskip\fi}
\def\firstcondition (#1){\hangindent\parindent{\rm(#1)}\enspace
 \ignorespaces}
\def\freeprod{\mathop{\prod\kern-10.9pt{*}\kern5pt}}
\def\freeprodT{\prod\kern-9pt{*}\kern2pt}

\def\setminus{
   \def\hefreshD{\mathop{\raise1.5pt\hbox{${\smallsetminus}$}}}
   \def\hefreshS{\mathop{\raise0.85pt\hbox{$\scriptstyle\smallsetminus$}}}
   \mathchoice{\hefreshD}{\hefreshD}{\hefreshS}{\hefreshS}}

\def\leaderfill{\leaders\hbox to 1em{\hss.\hss}\hfill}

\def\longmapright#1#2{\smash{\mathop
                {\hbox to #1pt{\rightarrowfill}}\limits^{#2}}}

\def\newline{\hfill\break}

\def\part#1:#2.{\par\medskip \noindent {\smc Part #1:} {\sl #2.}\quad}
\def\proclaimauthor #1 (#2):{\medbreak\noindent{\smc#1}
{\rm(#2):}\enspace\sl\ignorespaces}
\def\proofof#1:{\par\medskip\noindent{\sl Proof of {\rm #1:}}}

\def\qed{\qquad\bull\vadjust{\medskip}\ifmmode{}\else\par\fi}
\def\remark#1{}
\def\rmapdown#1
   {\Big\downarrow\rlap{$\vcenter{\hbox{$\scriptstyle#1$}}$}}

\def\shvor{\discretionary{}{}{}}
\def\sub#1{_{\lower2pt\hbox{$\scriptstyle #1$}}}
\def\sdp{\hbox to 7.56942pt{$\times$\hskip-0.25em{\raise 0.3ex
\hbox{$\scriptscriptstyle|$}}}}

\def\und{\ \wedge\ }

\def\pmb#1{\setbox0=\hbox{#1}%
\leavevmode
\kern-.025em\copy0\kern-\wd0
\kern.05em\copy0\kern-\wd0
\kern-.025em\raise.0433em\box0 }


\def\and{{\rm and}}

\def\id{{\rm id}}

\def\mod{\;{\rm mod}\hskip0.2em\relax}



\def\pcf{\mathop{\rm pcf}}
\def\Min{\mathop{\rm Min}}
\def\otp{\mathop{\rm otp}}
\def\nacc{\mathop{\rm nacc}}
\def\sqed#1{\qquad\bull\hbox{$_{#1}$}\vadjust{\medskip}\ifmmode{}%
\else\par\fi}
\def\Max{\mathop{\rm Max}}
\def\Dom{\mathop{\rm Dom}}
\def\cov{\mathop{\rm cov}}
\def\acc{\mathop{\rm acc}}
\def\tlim{\mathop{\rm tlim}}
\def\miss{{\not\subseteq}}
\def\Ord{\mathop{\rm Ord}}
\def\tcf{\mathop{\rm tcf}}
\def\pre{\mathop{\rm pre}}

\def\pp{\mathop{\rm pp}}
\def\Reg{\mathop{\rm Reg}}
\def\cf{\mathop{\rm cf}}
\def\Rang{\mathop{\rm Rang}}
\def\FIL{\mathop{\rm FIL}}
\def\Card{\mathop{\rm Card}}
\def\SCH{\mathop{\rm SCH}}
\def\und{\mathbin{\&}}
\def\Ch{{\rm Ch}}

\def\gl{\mathop{\rm g\ell}}
\def\lg{\mathop{\rm \ell g}}

\newdimen\contindent \contindent=15pt
\def\lbox#1{\hbox to\contindent{#1\hfill}}
\def\rbox#1{\hbox to\contindent{\hfill #1}}














%



\def\IJM{Israel Journal of Mathematics}









\def\PJM{Pacific Journal of Mathematics}





\def\APAL{Annals of Pure and Applied Logic}
%






\def\PJM{Pacific Journal of Mathematics}














\def\DMJ{Duke Mathematical Journal}




\def\BAMS{Bulletin of the American Mathematical Society}

\def\AML{Annals of Mathematical Logic}
\def\NDJFL{Notre Dame Journal of Formal Logic}
\def\ARML{Archive for Mathematical Logic}




 
  


%
\firstpage=61
\lastpage=114
\ijvol{95}
\ijyear{1996}
\recdate{Received January 23, 1992 and in final revised form April 22, 1996}
\shortauthor{S. SHELAH}
\shorttitle{FURTHER CARDINAL ARITHMETIC}
\title{FURTHER CARDINAL ARITHMETIC}
\author{Saharon Shelah\footnote*{Done mainly 1--4/1991.\newline
I thank Alice Leonhardt for typing
and retyping so beautifully and accurately.\newline
Partially supported by the Basic Research Fund, Israel Academy of
Sciences.\newline
Publication number 430.}}
\affiliation{Institute of Mathematics\\
The Hebrew University of Jerusalem, Jerusalem, Israel\\
and\\
Department of Mathematics\\
Rutgers University, New Brunswick, NJ, USA}
\addauthor{}
\affiliation{}
%
\beginabstract		
We continue the investigations in the author's book on cardinal
arithmetic, assuming some knowledge of it.
We deal with the cofinality of
$(\calS_{\le\aleph_0}(\kappa),\subseteq)$ for $\kappa$ real valued
measurable (Section 3), densities of box products (Section 5,3),
prove the equality $\cov (\lambda,\lambda,
\theta^+,2)=\pp(\lambda)$ in more cases even when
$\cf(\lambda)=\aleph_0$ (Section 1), deal with bounds of
$\pp(\lambda)$ for
$\lambda$ limit of inaccessible (Section 4) and give proofs to various
claims I was sure I had already written but did not find (Section 6).
\endabstract		


\maketopmatter	        
\begingroup
\ninepoint

\centerline{\bf Annotated Contents}

\line{\lbox{1.} Equivalence of two
covering properties\leaderfill\rbox{63}}
{\parindent=\contindent\narrower\noindent
[We try to characterize when, say,  $\lambda $  has few
countable subsets; for a given  $\theta \in  (\aleph_{0},\lambda )$, 
we try to translate to expressions with $\pcf$'s  the cardinal 
$$
\eqalign{
\Min
\big\{
|\calP|\colon & \calP \subseteq\calS_{<\mu }(\lambda )
\hbox{ and every }
\ a \in\calS_{\leq\theta} (\lambda )\ \hbox{ is }
\  \bigcup_{n<\omega} a_n,
\hbox{ such that every } \cr
& b\in \bigcup_n S_{\leq\aleph_0}(a_n)
\hbox{ is included in a member of }\ \calP \big\}.
\cr}
$$

This continues and improves [Sh410,{\S}6].]%
}

\line{\lbox{2.}Equality relevant to weak diamond\leaderfill\rbox{70}}
{\parindent=\contindent\narrower\noindent
[We show that if  $\mu  > \lambda \geq  \kappa $,  $\theta = 
\cov(\mu ,\lambda ^+,\lambda ^+,\kappa )$  and $\cov(\lambda ,\kappa,
\kappa ,2)\leq\mu$  (or $\leq\theta)$, then  $\cov(\mu ,\lambda ^+,
\lambda ^+,2) = \cov(\theta ,\kappa ,\kappa ,2)$. 
This is used in [Sh-f, Appendix,{\S}1] to 
clarify the conditions for the holding of versions of the weak
diamond.]

}

\line{\lbox{3.}Cofinality of ${\cal S}_{\leq \aleph_{0}}(\kappa )$
for $\kappa$ real valued measurable and trees\leaderfill\rbox{72}}
{\parindent=\contindent\narrower\noindent
[Dealing with partition theorems on trees, Rubin--Shelah [RuSh117]
arrive at the statement:  $\lambda > \kappa > \aleph_{0}$ are regular,
$a_{\alpha}  \in  {\cal S}_{<\kappa }(\mu )$,  $\mu  < \lambda $; 
can we find unbounded $W \subseteq  \lambda $  such that  $|
\bigcup_{{\alpha} \in  {W}}a_{\alpha} | < \kappa$?  Of course, 
$ \bigwedge_{{\alpha} < {\lambda}} \cov(\alpha ,\kappa ,\kappa ,2)
< \lambda $  suffice, but is it necessary? 
By 3.1, yes.  Then we answer a problem of Fremlin: e.g.\ if  
$\kappa $  is a real valued measurable cardinal then the cofinality of  
$\left( {\cal S}_{\leq \aleph_{0}}(\kappa ),\subseteq \right) $  is 
$\kappa $. 
Lastly we return to the problem of the existence of trees with many
branches (3.3, 3.4).]

}

\line{\lbox{4.}Bounds for $\pp_{\Gamma (\aleph_{1})}$ for limits
of inaccessibles\leaderfill\rbox{79}}
{\parindent=\contindent\narrower\noindent 
[Unfortunately, our results need an assumption: $\pcf(\fra)$ 
does not have an inaccessible accumulation point  $(|\fra| < \Min\fra$, 
$ \fra \subseteq  \Reg$,  of course). 
Our main conclusion (4.3) is that e.g.\ if  $\langle \lambda_{\zeta}
\colon \zeta < \omega_{4}\rangle $  is the list of the first 
$\aleph_{4}$ inaccessibles then  $\pp_{\Gamma (\aleph_{1})}\left(
\bigcup_{{\zeta} < {\omega_{1}}}\lambda_{\zeta} \right)  <
\bigcup_{{\zeta} < {\omega_{4}}}\lambda_{\zeta} $. 
This does not follow from the proof of $\pp\aleph_\omega <
\aleph_{\omega_4}$ [Sh400,{\S}2], nor do we make our life easier by
assuming $`` \bigcup_{{\zeta} < {\omega_{1}}}\lambda_{\zeta} $ is strong
limit".  We indeed in the end quote a 
variant of [Sh400,{\S}2] (= [Sh410,3.5]). 
But the main point now is to 
arrive at the starting point there: show that for 
$\delta < \omega_{4}$,  
$\cf \delta = \aleph_{2}$,  for some club $C$  of  $\delta $,  $\sup  
\pcf_{\aleph_{2}\hbox{-complete}}(\{\lambda_{\zeta} \colon \zeta \in 
C\})$ is $\leq\lambda_{\delta} $.  This is provided by 4.2.]

}

\line{\lbox{5.}Densities of box products\footnote*{There is a paper in
preparation on independence results by Gitik and
Shelah.}\leaderfill\rbox{85}}
{\parindent=\contindent\narrower\noindent
[The behavior of the Tichonov product of topological spaces on densities
is quite well understood for  $^\mu 2:$  it is  $\Min \{\lambda
\colon 2^\lambda \geq  \mu \}$;
but less so for the generalization to box products.  Let  
${\cal T}_{\mu ,\theta ,\kappa }$ be the space with set of points
$^\mu \theta $,  and basis  $\{[f]\colon f$\ \ a\ partial 
function\ from\ \ $\mu $\ \ to\ \ $\theta $\ \ of\ cardinality\ \ $< 
\kappa\}$, where $[f]=\{g \in ^\mu\!\!\theta\colon f \subseteq  g\}$. 
If $\theta \leq  \lambda = \lambda ^{<\kappa }$,
$2^\lambda \geq  \mu $  the 
situation is similar to the Tichonov product. 
Now the characteristic unclear case is  $\mu $  strong limit singular of
cofinality  $< \kappa $,  $\theta = 2$,  $2^\mu  > \mu ^+$. 
We prove that the density is ``usually" large  
$(2^\mu )$, i.e.\ the failure quite limits the cardinal arithmetic
involved (we can prove directly consistency results but what we do
seems more informative).]

}

\line{\lbox{6.}Odds and ends\leaderfill\rbox{90}}

\line{\lbox{}References\leaderfill\rbox{112}}

\line{\lbox{}}

\line{\lbox{}Notation:
{\parindent=\contindent\narrower\noindent
Let $J_{\lambda} [\fra]$ be $\{\frb\subseteq\fra\colon\lambda
\notin \pcf (\frb)\}$, equivalently  $J_{<\lambda }[\fra] +
\frb_{\lambda} [\fra].$}
}

\line{\lbox{}\ 
{\parindent=\contindent\narrower\noindent 
See more in [Sh513], [Sh589].\hfill}
}

\endgroup



\head 1. Equivalence of Two Covering Properties\endhead

\proclaim{1.1 Claim} 
If  $\pp\lambda = \lambda ^+$,  $\lambda > \cf(\lambda)=\kappa>
\aleph_{0}$ then  $\cov(\lambda ,\lambda ,\kappa^+,2) = \lambda ^+.$
\endproclaim

\demo{Proof} 
Let  $\chi  = \beth_{3}(\lambda )^+$; \ choose  $\langle
\frB_{\zeta} \colon \zeta < \lambda ^+\rangle $  
increasing continuous, such that  $\frB_{\zeta}  \prec 
(H(\chi ),\in ,<^\ast_{\chi} )$,  $\lambda + 1 
\subseteq \frB_{\zeta} $,  $\|\frB_{\zeta} \| = \lambda $  and
$\langle \frB_{\xi} \colon \xi \leq\zeta\rangle\in  \frB_{\zeta +1}$. 
Let  $\frB =: \bigcup_{{\zeta} < {\lambda ^+}} \frB_{\zeta} $ and
${\cal P} =: {\cal S}_{<\lambda }(\lambda ) \cap  \frB$. 
Let  $a \in  {\cal S}_{\leq \kappa}(\lambda )$;  it suffices to
prove  $(\exists A \in  {\cal P})[a \subseteq A]$. 
Let  $f_{\xi} $ be the $<^\ast_{\chi} $-first  $f \in  \prod
(\Reg \cap  \lambda )$  such that  $(\forall g)[g \in  \prod
(\Reg\cap\lambda)\und g\in\frB_{\zeta}  \Rightarrow  g < f$\ mod\
$J^{bd}_{\lambda} ]$,  such  $f$  exists as  $\prod (\Reg \cap 
\lambda )/J^{bd}_{\lambda} $ is $\lambda ^+$-directed. 

By [Sh420, 1.5, 1.2] we can find $\langle C_{\alpha} \colon \alpha < 
\lambda ^+\rangle $  such that:  $C_{\alpha} $ is a closed subset of 
$\alpha $, $\otp C_{\alpha}  \leq  \kappa^+$, 
$[\beta \in  \nacc C_{\alpha}  \Rightarrow  C_{\beta}  = C_{\alpha} 
\cap  \beta ]$  and $S =: \{\delta < \lambda ^+\colon  
\cf(\delta)= \kappa^+$ and $\delta = \sup  C_{\delta} \}$  is
stationary.

Without loss of generality  $\bar C \in  \frB_{0}.$

Now we define for every  $\alpha < \lambda ^+$ elementary submodels  
$N^0_{\alpha} $, $N^1_{\alpha} $ of  $\frB$: 

$N^0_{\alpha}$ is the Skolem Hull of $\{f_{\zeta}\colon\zeta\in 
C_{\alpha} \}\cup \{ i\colon i\le\kappa\}$  and $N^1_{\alpha} $ is the
Skolem Hull of\break $a \cup  \{f_{\zeta} \colon \zeta \in  
C_{\alpha} \}\cup \{i: i\leq \kappa\}$,
both in  $(H(\chi ),\in ,<^\ast_{\chi} ).$

Clearly:  
\item{(a)}  $N^0_{\alpha}  \subseteq N^1_{\alpha}  \subseteq
\frB_{\alpha}  \subseteq \frB$ 
[why? as $f_{\zeta} \in \frB_{\zeta +1}$ because $\frB_{\zeta}  \in  
\frB_{\zeta +1}]$,

\item{(b)}  $\|N^\ell_{\alpha} \| \leq  \kappa + \|C_{\alpha} \|$,

\item{(c)}  $N^0_{\alpha}  \in  \frB_{\alpha +1}$.

\noindent
[Why?  As  $\alpha \subseteq \frB_{\alpha} $ (you can prove it by
induction on  $\alpha )$  clearly  $\alpha 
\in  \frB_{\alpha +1}$,  but  $\bar C \in  \frB_{0} 
\subseteq \frB_{\alpha +1}$;  hence  $C_{\alpha}  \in \frB_{\alpha +1}$,
also $\langle\frB_{\gamma}\colon\gamma\leq\alpha\rangle
\in  \frB_{\alpha +1}$ hence  $\langle f_{\gamma}\colon\gamma\leq
\alpha \rangle  \in  \frB_{\alpha +1}$,
hence  $\langle f_{\gamma} \colon \gamma \in  C_{\alpha} \rangle  \in 
\frB_{\alpha +1}$. 
Now  $N^0_{\alpha}  \subseteq \frB_{\alpha}  \in  \frB_{\alpha +1}$
and the Skolem Hull can be computed in  $\frB_{\alpha +1}.]$

\item{(d)} for each $\alpha $  with $\kappa^+>\otp (C_{\alpha})$, for
some $\gamma_{\alpha}<\lambda^+$,  letting  $\fra_{\alpha}  =:$\break
$N^0_{\alpha}\cap\Reg\cap\lambda\backslash\kappa^{++}$ clearly
$\Ch_{\alpha} 
\in  \prod \fra_{\alpha} $ where  $\Ch_{\alpha} (\theta ) =: \sup
(\theta \cap  N^1_{\alpha} )$,  and we have:
$\Ch_{\alpha}  < f_{\gamma_{\alpha} }\upharpoonright
\fra_{\alpha}\mod J^{bd}_{\fra_{\alpha} }$.

\noindent
[Why?  $\fra_{\alpha}  \in  \frB_{\alpha +1}$ as  $N^0_{\alpha}  \in 
\frB_{\alpha +1}$, and $\prod\fra_{\alpha} /J^{bd}_{\fra_{\alpha} }$ is
$\lambda ^+$-directed (trivially) and has cofinality $\leq  \max 
\pcf_{J^{bd}_{\fra_{\alpha} }}(\fra_{\alpha} ) \leq  
\pp(\lambda ) = \lambda ^+$,  so there is  $\langle
f^{\fra_{\alpha} }_{\beta} \colon \beta 
< \lambda ^+\rangle $,  $<_{J^{bd}_{\fra_{\alpha} }}$-increasing cofinal
sequence in $\prod \fra_{\alpha} $,  so without loss of generality 
$\langle f^{\fra_{\alpha}}_{\beta}\colon\beta < \lambda ^+\rangle  \in 
\frB_{\alpha +1}$; 
also by the ``cofinal" above, for some $\beta\in(\alpha,\lambda^+)$, 
$\Ch_{\alpha}  < 
f^{\fra_{\alpha} }_{\beta} $\ mod\ $J^{bd}_{\fra_{\alpha} }$. 
We can use the minimal  $\beta $,  now obviously  $\beta \in 
\frB_{\beta +1}$ so  $f^{\fra_{\alpha} }_{\beta}  
\in  \frB_{\beta +1}$,  hence  $f^{\fra_{\alpha} }_{\beta}  < 
f_{\beta +2}$\ mod\ $J^{bd}_{\lambda} $. 
Together  $\gamma_{\alpha}  =: \beta + 2$  is as required.]

\item{(d)$^+$} for each $\alpha$ with $\otp(C_\alpha)<\kappa^+$
for some $\gamma_\alpha\in(\alpha,\lambda^+)$, for any
$\mu\in\Reg\cap N^0_\alpha$, letting $N^{0,\mu}_\alpha =:
\Ch_{\frB_\alpha} (N^0_\alpha\cup \mu)$,
$\fra_{\alpha,\mu} =
N^{0,\mu}_\alpha\cap\Reg\cap\lambda\setminus\mu^+$ and
$\Ch_{\alpha,\mu}\in\Pi\fra_{\alpha,\mu}$ be
$$
\Ch_{\alpha,\mu}(\theta) =
\cases{
\sup(\theta\cap N^1_\alpha) &  if $\theta\in N^1_\alpha$,\cr 
0 & otherwise,\cr}
$$
we have: $\Ch_\alpha< f_{\gamma_\alpha} \upharpoonright
\fra_{\alpha,\mu}\mod J^{bd}_{\fra_{\alpha,\mu}}$.

\noindent
[Why? Clearly $\Ch_{\frB_\alpha}(N^0_\alpha\cup \mu)\in
\frB_{\alpha+1}$, so $\fra_{\alpha,\mu}\in\frB_{\alpha+1}$,
hence there are in $\frB_{\alpha+1}$ elements $\langle
\frb_\theta [\fra_{\alpha,\mu}]\colon\theta\in
\pcf(\fra_{\alpha,\mu})\rangle$ and
$\langle\langle f^{\fra_{\alpha,\mu},\theta}_\alpha
\colon\alpha<\theta\rangle\colon\theta\in\pcf(\fra_{\alpha,\mu})\rangle$
as in [Sh 371, 2.6, \S1].
So for some $\gamma_{\alpha,\mu}\in(\alpha,\lambda^+)$ we have
$\Ch_\alpha\upharpoonright\frb_{\lambda^+} [\fra_{\alpha,\mu}] <
f_{\gamma_\alpha}$, so it is enough to prove
$\fra_{\alpha,\mu}\hefresh \frb_{\lambda^+}
[\fra_{\alpha,\mu}]$ is bounded below $\mu$ but otherwise
$\pp(\lambda)=\lambda^+$ will be contradicted.
Let $\gamma_\alpha=\sup\{\gamma_{\alpha,\mu}\colon \mu\in
N^0_\alpha\}$.]

\item{(e)} $E^\ast  =: \{\delta < \lambda ^+\colon \alpha < \delta
\und |C_{\alpha} | \leq 
\kappa \Rightarrow  \gamma_{\alpha}  < \delta\hbox{ and }\delta> \lambda \}$
is a club of  $\lambda$. 

Now as  $S$  is stationary, there is  $\delta (\ast ) \in  S \cap 
E^\ast $.  
Remember $\otp C_{\delta (\ast )} = \kappa^+.$ 

Let  $C_{\delta (\ast )} = \{\alpha_{\delta (\ast ),\zeta }\colon
\zeta < \kappa^+\}$  (in increasing order).

\noindent
Let  (for any  $\zeta < \kappa^+$) $M^0_{\zeta} $ be the Skolem Hull
of
$\left\lbrace f_{\alpha_{\delta (\ast ),\xi }}\colon \xi  < \zeta \right
\rbrace \cup \{i: i\leq \kappa\}$,
and let  $M^1_{\zeta} $ be the Skolem Hull of  $a \cup
\left\lbrace f_{\alpha_{\delta (\ast ),\xi }}\colon \xi  < \zeta \right
\rbrace \cup \{i: i\leq \kappa\}$.
Note: for  $\zeta < \kappa^+$ non-limit
$\left\lbrace f_{\alpha_{\delta (\ast ),\xi }}\colon \xi  < \zeta
\right\rbrace  = 
\left\lbrace f_{\xi} \colon \xi  \in  
C_{\alpha_{\delta (\ast ),\zeta }}\right\rbrace $.
Clearly  $\langle M^0_{\zeta} \colon \zeta < \kappa^+\rangle $, 
$\langle M^1_{\zeta} \colon \zeta < \kappa^+\rangle $  are increasing
continuous sequences of countable elementary 
submodels of  $\frB$  and $M^0_{\zeta}  \subseteq 
M^1_{\zeta} $ and for  $\zeta < \kappa^+$ a successor ordinal,  
$N^\ell_{\alpha_{\delta (\ast ),\zeta }} = M^\ell_{\zeta} $.

Now for each successor  $\zeta $,  for some  $\epsilon (\zeta ) \in  
(\zeta ,\omega_{1})$  we have  $\gamma_{\alpha_{\delta (\ast ),\zeta }} 
< \alpha_{\delta (\ast ),\epsilon (\zeta )}$ (by the choice of  
$\delta (\ast ))$  hence
$f_{\gamma_{\alpha_{\delta (\ast ),\zeta }}} < 
f_{\alpha_{\delta (\ast ),\epsilon (\zeta )}}\mod\ J^{bd}_{\lambda} $
hence  $\Ch_{\alpha_{\delta (\ast ),\zeta }} < 
f_{\alpha_{\delta (\ast ),\epsilon (\zeta )}}\mod\ J^{bd}_{\lambda} $.

Let $E =: \{\delta<\omega_{1}$: for every successor  $\zeta < \delta$,
$\epsilon (\zeta ) < \delta \}$,  clearly  $E$  is a club of 
$\kappa^+$.
Let $\lambda=\sum_{i<\kappa}\lambda_i$, $\lambda_i<\lambda$
singular increasing continuous with $i$, wlog $\{\lambda_i\colon
i<\kappa\}\subseteq \Ch_\frB (\{ i\colon i\le
\kappa\}\cup\{\lambda\})$.
So for some  $\mu_{\zeta,i}  < \lambda $,  we have: 
$$
\eqalign{
i<\kappa,
\quad \zeta = & \xi  + 1 < \kappa^+ \und \theta \in  \Reg \cap
\lambda \backslash \mu_{\zeta,i}  \und \theta \in  
N^{0,\lambda_i}_{\alpha_{\delta (\ast ),\zeta }}\cap
N^1_{\alpha_{\delta(\ast),\zeta}}\cr
\Rightarrow &  \sup \left( N^1_{\alpha_{\delta (\ast ),\zeta }}\cap  
\theta \right)  < f_{\alpha_{\delta (\ast ),\epsilon (\zeta )}}(\theta )
\in  \theta \cap  N^{0,\lambda_i}_{\alpha_{\delta (\ast ),\zeta +1}}.
\cr}
\leqno(*)
$$
So for some limit $i(\zeta)<\kappa^+$ we have
$\lambda_{i(\zeta)}=\sup\{\mu_{\zeta,j}\colon j<i(\zeta)\}$.
Now as  $\cf \lambda \leq \kappa^+$ for some  $i(\ast) < \lambda $
$$
W =: \{\zeta < \kappa^+\colon   \zeta \hbox{  successor ordinal 
and } i(\zeta) = i(*) \}
$$
is unbounded in  $\kappa^+$.
So 
$$
\eqalign{
& \hbox{ if } \xi  < \kappa^+,\ \xi\in E,\
\xi  = \sup (\xi  \cap  W)
\hbox{ and } \theta \in M_\xi^1\Reg\cap\lambda\cap
M^{0,\lambda_{i(*)}}_{\xi} \backslash \lambda_{i(*)} \cr
& \hbox{ then  }M^{0,\lambda_i}_{\xi} \cap \theta
\hbox{ is an unbound subset of }
M^1_{\xi} \cap  \theta .\cr}
\leqno{\otimes}
$$
Hence by [Sh400] 5.1A(1), remembering
$M^0_{\zeta+1}= N^0_{\alpha_{\delta(*),\zeta+1}}$,  we
have:  
$M^1_{\xi}  \subseteq $\ Skolem\ Hull$\left[
\bigcup_{{\zeta} < {\xi}} N^0_{\zeta +1} \cup\lambda_{i(*)} \right] 
\subseteq$ Skolem Hull $\left( N^0_{\alpha_{\delta (\ast ),\xi +1}}
\cup  \lambda_{i(*)}\right)$ whenever $\xi\in E$ is an accumulation
point of  $W$.
But $a \subseteq M^1_{\xi}$ and the right side belongs to $\frB$
(as we can take the Skolem Hull in  $\frB_{\delta (\ast )})$. 
So we have finished.
\sqed{1.1}

\rem{Remark}
Alternatively note:
$\cov(\lambda,\lambda,\kappa,2)\le\cov(\theta,\lambda,\sigma,2)$
when $\sigma=\cf(\lambda)<\kappa<\lambda$, $\sigma=>\aleph_0$,
$\theta=\pp_{\Gamma(\kappa,\sigma)}(\lambda)$;
remember $\cf(\lambda)<\kappa<\lambda\ \&\
\pp(\lambda)<\lambda^{+\kappa^+}\Rightarrow
\pp_{<\lambda}(\lambda) = \pp(\lambda)$.

\proclaim{1.2 Claim} 
For  $\lambda > \mu  = \cf(\mu) > \theta > \aleph_{0}$,  we have  
$\lambda (0) \leq  \lambda (1) \leq  \lambda (2) = \lambda (3)$  and if 
$\cov(\theta ,\aleph_{1},\aleph_{1},2) < \mu $  they are all equal,
where:
$$
\eqalign{
\lambda (0) =:&\hbox{ is the minimal }\kappa\ \hbox{ such that: if }\fra
\subseteq \Reg \cap  \lambda ^+\backslash \mu , 
|\fra| \leq  \theta \hbox{ then we    }\cr
& \hbox{ can find } \langle \fra_{\ell} \colon \ell  < \omega \rangle
\hbox{ such that } \fra = \bigcup\nolimits_{{\ell} < {\omega}}
\fra_{\ell} \hbox{ and }\cr
& (\forall \frb)\left[ \frb \in  {\cal S}_{\leq \aleph_{0}}(\fra_{n})
\Rightarrow  \max  \pcf (\frb) \leq  \kappa \right] .\cr}
$$
$$
\eqalign{
\lambda (1) =: &  \Min \left\lbrace |{\cal P}|\colon {\cal P} 
\subseteq 
{\cal S}_{<\mu }(\lambda )\right., \ \hbox{  and for every }
A \subseteq \lambda,  |A| \leq  \theta  \hbox{ there }\cr
& \hbox{ are } A_{n} \subseteq A\ (n < \omega ), A =
\bigcup\nolimits_{{n} < {\omega}} A_{n}, A_{n} \subseteq A_{n+1}
\hbox{ such }\cr
& \hbox{ that: for } n < \omega,\ \hbox{  every } a \in 
{\cal S}_{\leq \aleph_{0}}(A_{n})\ \hbox{  is a subset }\cr          
& \hbox{ of some member of } {\cal P}\left.\!\right\rbrace .\cr}
$$
$\lambda (2)$  is defined similarly to  $\lambda (1)$  as: 
$$
\eqalign{
\Min \Big\lbrace |{\cal P}|\colon &  {\cal P} \subseteq  
{\cal S}_{<\mu }(\lambda )\ \hbox{ and for every } A \in  
{\cal S}_{\leq \theta }(\lambda ) \ \hbox{ for some }
A_{n} \subseteq  A (n < \omega )\cr
& A =\bigcup_{n<\omega} A_n \hbox{ and for each } n < \omega\
\hbox{ for some } {\cal P}_{n} 
\subseteq {\cal P},\ |{\cal P}_{n}| < \mu ,\cr
& \sup_{B\in {\cal P}_{n}}|B| < \mu \hbox{  and every } a \in 
{\cal S}_{\leq \aleph_{0}}(A_{n})\ \hbox{  is a subset of some }\cr
& \hbox{ member of } {\cal P}_{n} \Big\rbrace .\cr}
$$
\item{$\lambda (3)$}  is the minimal  $\kappa $  such that: if  $\fra 
\subseteq 
\Reg \cap  \lambda ^+\backslash \mu $,  $|\fra| \leq  \theta $,  then  
we can find $\langle \fra_{\ell} \colon \ell  < \omega \rangle $, 
$\fra_{\ell}  \subseteq \fra_{\ell +1} \subseteq 
\fra = \bigcup_{{\ell} < {\omega}} \fra_{\ell} $ such that:  
there is  $\{\frb_{\ell ,i}\colon i < i_{\ell}  < \mu \}$,  $\frb_{\ell ,i} 
\subseteq 
\fra_{\ell} $ such that  $\max\pcf  \frb_{\ell ,i} \leq  \kappa $  
and $(\forall \frc)[\frc \subseteq \fra_{\ell}  \und |\frc| \leq  \aleph_{0} \Rightarrow   \bigvee_{i} \frc 
\subseteq 
\frb_{\ell ,i}]$;  equivalently:  ${\cal S}_{\leq \aleph_{0}}(\fra_{n})$  
is included in the ideal generated by  $\{\frb_{\sigma} [\fra_{n}]\colon
\sigma
\in  \frd\}$  for some  $\frd \subseteq 
\kappa ^+ \cap  \pcf \fra_{n}$ of cardinality  $< \mu .$

\rem{1.2A Remark} 
(1) We can get similar results with more parameters: replacing  
$\aleph_{0}$ and/or $\aleph_1$ by higher cardinals.

(2)  Of course, by assumptions as in [Sh410, {\S}6] (e.g.\ 
$|\pcf \fra| \leq  |\fra|)$  we get  $\lambda (0) = \lambda (3).$
This (i.e.\ Claim 1.2) will be continued in [Sh513].

\demo{Proof}  
\subdemo{$\lambda(1)\le\lambda(2)$}
Trivial.

\subdemo{$\lambda(2)\le\lambda(3)$}
Let  $\chi  = \beth _{3}(\lambda (3))^+$ and for\ 
$\zeta \leq  \mu ^+$ we choose  $\frB_{\zeta}  \prec 
(H(\chi ),\in ,<_{\chi} )$,  $\{\lambda ,\mu ,\theta ,\lambda (2),
\lambda (3)\}\in\frB_{\zeta}$,\ $\|\frB_{\zeta}\| =\lambda(3)$  and
$\lambda (3) \subseteq \frB_{\zeta}$, $\frB_{\zeta}$ $(\zeta\leq\mu^+)$ 
increasing continuous and $\langle\frB_{\xi}\colon\xi\leq\zeta\rangle
\in  \frB_{\zeta +1}$ and let  $\frB = \frB_{\mu ^+}$. 
Lastly let  ${\cal P} = \frB \cap  {\cal S}_{<\mu }(\lambda )$.  
Clearly 
\item{$(*)_{0}$} a function $\fra\mapsto\langle\frb_{\sigma}
[\fra]\colon \sigma\in\pcf\fra\rangle$ as in [Sh371, 2.6] is definable
in\break 
$(H(\chi ),\in ,<^\ast_{\chi})$ hence $\frB$ is closed under it.

It suffices to show that  ${\cal P}$  satisfies the requirements in the 
definition of  $\lambda (2).$ 

Let  $A \subseteq \lambda $,  $|A| \leq  \theta $. 
We choose by induction on  $n < \omega $,  
$N^a_{n}$,  (for  $\ell  < \omega )$  and $N^b_{n}$,  $f_{n}$ such that:

\item{(a)} $N^a_{n}$, $N^b_{n}$ are elementary submodels of  
$(H(\chi ),\in ,<^\ast_{\chi} )$  of cardinality  $\theta ,$ 

\item{(b)} $f_{n} \in  \prod \fra_{n}$ where  $\fra_{n} =: N^a_{n}
\cap  \Reg \cap  \lambda ^+\backslash \mu $,  and $f_{n}(\sigma ) >
\sup (N^b_{n} \cap  \sigma )$  (for any  $\sigma \in  \fra_{n})$,  

\item{(c)} $\theta + 1 \subseteq N^a_{n} \subseteq 
N^b_{n} \subseteq  \frB$,

\item{(d)} $N^b_{n}$ is the Skolem Hull of $\bigcup \{\Rang\
f_{\ell}\colon\ell  < n\} \cup  A \cup  (\theta +1)$,

\item{(e)} $N^a_{0}$ is the Skolem Hull of  $\theta + 1$  in  
$(H(\chi ),\in ,<^\ast_{\chi} )$,  

\item{(f)} $N^a_{n+1}$ is the Skolem Hull of $N^a_{n}\cup\Rang f_{n}$,

\item{(g)} there are  ${\cal P}_{n,\ell } \subseteq 
{\cal S}_{<\mu }(\lambda+1)$  and $A_{n,\ell } \subseteq  N^a_{n}$
(for $l<\omega$) such that:

\itemitem{$(\alpha )$} $|{\cal P}_{n,\ell }| < \mu $  and
$\mu_{n,\ell } =:  \sup_{B\in {\cal P}_{n,\ell }}|B| < \mu $  and
${\cal P}_{n,\ell } \subseteq  {\cal P}_{n,\ell +1}$, 

\itemitem{$(\beta )$} $N^a_{n} =  \bigcup_{\ell } A_{n,\ell }$, 
${\cal P}_{n} = \bigcup_{{\ell} < {\omega}} {\cal P}_{n,\ell }
\subseteq  \frB$  and $A_{n,\ell } \subseteq  A_{n,\ell +1}$,  

\itemitem{$(\gamma )$} for every countable  $a \subseteq 
\lambda \cap  A_{n,\ell }$ there is  $b \in  {\cal P}_{n,\ell }$
satisfying  $a \subseteq b$,

\itemitem{$(\delta )$} ${\cal P}_{n,\ell } = {\cal S}_{\leq
\mu_{n,\ell }}(\lambda +1) 
\cap  ($Skolem Hull of  $A_{n,\ell } \cup  {\cal P}_{n,\ell } \cup  
(\theta +1))$.

As in previous proofs, if we succeed to carry out the definition, then 
$ \bigcup_{n} \left( N^a_{n} \cap  \lambda \right)  \shvor = 
\bigcup_{n} N^b_{n} \cap  \lambda $,  but the former is 
$ \bigcup_{{n} , {\ell}} A_{n,\ell } \cap  \lambda $,  hence 
$A \subseteq  \bigcup_{n} \bigcup_{\ell } A_{n,\ell}$,
by (g)$(\alpha ),(\beta )$ the $\calP'_{n,\ell}=\{
a\cap\lambda\colon a\in\calP_{n,\ell}\}$ are of the 
right form and so by (g)$(\gamma )$ we finish. 

Note that without loss of generality: if  $a \in  {\cal P}_{n,\ell }$
then 
$a\cap\Reg\cap(\lambda+1)\backslash \mu  \in  {\cal P}_{n,\ell }.$ 

For\  $n = 0$  we can define  $N^a_{0}$,  $N^b_{0}$,  $A_{n,\ell }$
trivially.  
Suppose  $N^a_{m} $,  $N^b_{m} $,  $A_{m ,\ell }$,  ${\cal P}_{m ,
\ell }$ are defined 
for  $m  \leq  n$,  $\ell  < \omega $  and $f_{m}$ $(m  < n)$  are
defined. 
Now  $\fra_{n}$ is well defined and $ \subseteq 
\Reg \cap  \lambda ^+\backslash \mu  \subseteq 
\frB$  and $|\fra_{n}| \leq  \theta .$
So  $\fra_{n} =  \bigcup_{\ell } \fra_{n,\ell }$ and $\fra_{n,\ell } 
\subseteq \fra_{n,\ell +1}$ where  $\fra_{n,\ell } =: \fra_{n} \cap 
A_{n,\ell }$ and, of course,  $\fra_{n,\ell } \subseteq  \Reg \cap 
\lambda ^+\backslash \mu$ has cardinality  $\leq  \theta $. 
Note that  $\fra_{n,\ell }$ is not necessarily in  $\frB$  but 
\item{$(*)_{1}$} every countable subset of $\fra_{n,\ell }$ is included
in some subset of  $\frB$  which belongs to  ${\cal P}_{n,\ell }$ and
is $\subseteq  \Reg \cap  \lambda ^+\backslash \mu .$ 

By the definition of $\lambda (3)$ (see ``equivalently" there), for each
$n,\ell$ we can find an increase sequence $\langle\fra_{n,\ell ,k}\colon
k <\omega\rangle$ of subsets of $\fra_{n,\ell }$ with union
$\fra_{n,\ell }$ and $\frd_{n,\ell ,k} \subseteq [\mu ,\lambda (3)]\cap
 \pcf (\fra_{n,\ell ,k})$,  $|\frd_{n,\ell ,k}| < \mu $  such that: 
\item{$(*)_{2}$} if $\frb \subseteq \fra_{n,\ell ,k}$ is countable then
$\frb$   is included in a finite union of some members of  
$\{\frb_{\sigma} [\fra_{n,\ell ,k}]\colon \sigma \in  \frd_{n,\ell ,k}\}$   
(hence  $\max  \pcf (\frb) \leq  \lambda (3)).$ 

By the properties of $\pcf$: 
\item{$(*)_{3}$} for each  $\ell ,k < \omega $  and $\frc \subseteq 
\Reg \cap  \lambda ^+\backslash \mu $  such that  $\frc \in 
{\cal P}_{n,\ell }$ we can find $\fre = \fre^{\ell ,k}_{\frc}
\subseteq  \lambda (3)^+ \cap  \pcf
\frc$,  $|\fre| \leq  |\frd_{n,\ell ,k}| < \mu $  such that for every
$\sigma \in  \frd_{n,\ell ,k}$ we have: 
$\frc \cap\frb_{\sigma} [\fra_{n,\ell ,k}]$  is included in a finite
union of members of $\{\frb_{\tau}[\frc]\colon\tau\in\fre_{\frc}\}$.

\noindent
By [Sh371, 1.4] we can find $f_{n} \in \prod_{{\sigma} \in 
{\fra_{n}}}\sigma $  such that:  
\item{$(*)_{4}$} $(\alpha )$ $\sup (N^b_{n} \cap  \sigma ) < f_{n}
(\sigma )$;

\itemitem{$(\beta)$} if $\frc\in\calP_{n,\ell}$,  $\ell,k<\omega$,
$\frc \subseteq  \Reg \cap  \lambda ^+\backslash \mu $  and $\sigma \in  
\fre^{\ell ,k}_{\frc} \subseteq  \pcf (\frc) \cap  [\mu ,\lambda (3)]$\ 
(where  $\fre^{\ell ,k}_{\frc}$ is from $(\ast )_{3})$ {\it then}
for some 
$m  < \omega $,  $\sigma_{p} \in  \sigma ^+ \cap  \pcf (\frc)$ 
and $\alpha_{p} < \sigma_{p}$ ,  (for  $p \leq  
m )$  the function  $f_{n}\upharpoonright (\frb_{\sigma} [\frc])$  is included in 
$ \Max_{p\leq m } f^{\frc,\sigma_{\ell} }_{\alpha_{p}}\upharpoonright
\frb_{\sigma_{p}}[\frc]$  (the $\Max$ taken pointwise).

\noindent
Note 
\item{$(*)_{5}$} if $\frb\subseteq\fra_{n,\ell ,k}$ is countable (where 
$\ell ,k <\omega)$ {\it then} there is  $\frc \in  {\cal P}_{n,\ell }$, 
$|\frc| < \mu $,  $\frc \subseteq  \Reg \cap  \lambda ^+\backslash \mu $ 
such that  $\frb \subseteq  \frc.$

\noindent
By $(*)_{4}:$\  
\item{$(*)_{6}$} if  $\ell ,k < \omega $,  $\frc \in  {\cal P}_{n,\ell }$,
$\frc \subseteq  \Reg \cap  \lambda ^+\backslash \mu $,  and $\sigma \in  
\frd_{n,\ell ,k} \cap  \lambda (3)^+ \cap  \pcf \frc\backslash \mu$
then  $f_{n}\upharpoonright \frb_{\sigma} [\frc] \in  \frB.$

\noindent
You can check that (by $(*)_{2} - (*)_{6}):$ 

\item{$(*)_{7}$} if  $\frb \subseteq \fra_{n,\ell ,k}$ is countable
{\it then} there is  $f^{n,\ell ,k}_{\frb} \in  \frB$,  
$|\Dom \ f^{n,\ell ,k}_{\frb}| < \mu $  such that  $f_{n}
\upharpoonright \frb \subseteq f^{n,\ell ,k}_{\frb}.$

\noindent
Let $\tau_{i}(i < \omega )$  list the Skolem function of  
$(H(\chi ),\in ,<^\ast_{\chi} )$. 
Let
$$
A_{n+1,\ell } = \bigcup \left\lbrace \right. \Rang\left( \tau_{i}
\upharpoonright (A_{n,j} \cup\Rang\ f_{n}\upharpoonright \fra_{n,j,k})
\right) \colon i < \ell ,\quad j < \ell ,\quad k < \ell \left.
\right\rbrace ,
$$
$$
\calP'_{n+1,\ell } = \bigcup_{m\leq\ell}
\calP_{n,m } \cup  \big\{ f_n\upharpoonright \fra'\colon
\fra' \in \bigcup_{m\leq\ell} \calP_{n,m }\ 
\hbox{  and  }
f_{n}\upharpoonright \fra' \in \frB \big\},
$$ 
and ${\cal P}_{n+1,\ell }
= {\cal S}_{<\mu }(\lambda +1) \cap  \left( \right. $Skolem Hull of  
$A_{n+1,\ell } \cup  {\cal P}'_{n+1,\ell} \cup  (\theta + 
1)\left. \right) $.\ 

So  $f_{n}$,  ${\cal P}_{n+1,\ell}$ are as required.

Thus we have carried the induction.

\subdemo{$\lambda(3)\le\lambda(2)$}
Let  ${\cal P}$  exemplify the definition of  $\lambda (2)$. 
Let  $\fra \subseteq 
\Reg \cap  \lambda ^+\backslash \mu $,  $|\fra| \leq  \theta (< \mu )$. 
Let  $J = J_{\leq \lambda (2)}[\fra]$,  and let
\item{$J_{1} =$} $\left\lbrace \frb\colon \frb 
\subseteq \right.  \fra$ and there is  $\langle \frb_{i}\colon i < i^\ast
\rangle $,  satisfying:  $\frb_{i} \subseteq \frb$,  $i^\ast  < \mu ,$  
$\max\pcf  \frb_{i} \leq  \lambda (2)$  and any countable subset of 
$\frb$  is in the ideal    
which  $\{\frb_{i}\colon i < i^\ast \}$  generates$\left. \right\rbrace $.\ 

\noindent
Clearly  $J_{1}$ is an ideal of subsets of $\fra$ extending  $J$. 
Let 
$$
J_2 =\left\lbrace \frb\colon \hbox{ for some } 
\frb_{n} \in  J_{1} \hbox{ (for } n < \omega ),\ \frb\subseteq
\bigcup_{n} \frb_{n} \right\rbrace .
$$ 

\noindent
Clearly  $J_{2}$ is an $\aleph_{1}$-complete ideal extending  $J_{1}$
(and $J)$.  
If  $\fra \in J_{2}$ we have that $\fra$ satisfies the requirement
thus we have finished so we can assume  $\fra \notin  J_{2}$. 
As we can force by Levy $\left( \lambda (2)^+,2^{\lambda (2)}\right)$
(alternatively, replacing  $\fra$ by [Sh355, {\S}1]) without loss of
generality  $\lambda
(2)^+ = \max  \pcf \fra$ and so  $\tcf(\prod \fra/J_{2}) =
\tcf(\prod \fra/J) = \lambda (2)^+$.  Let  $\bar f = 
\langle f_{\alpha} \colon \alpha < \lambda (2)^+\rangle$
be $<_{J}$-increasing,  $f_{\alpha}  \in  \prod \fra$,  cofinal in 
$\prod \fra/J$. 
Let  $\frB \prec (H(\chi ),\in ,<^\ast_{\chi} )$  be of cardinality 
$\lambda (2)$,  $\lambda (2) + 1 \subseteq 
\frB$, $\fra\in\frB$, $\bar f\in\frB$ and ${\cal P}\in\frB$.
Let  ${\cal P}' =: \frB \cap  {\cal S}_{<\mu }(\lambda ).$

For  $B \in  {\cal P}'$  (so  $|B| < \mu )$  let  $g_{B} \in 
\prod \fra$  be  $g_{B}(\sigma ) =: \sup (\sigma \cap  B)$,  so for
some  $\alpha_{B} < \lambda $, $g_{B} <_{J} f_{\alpha_{B}}$. 
Let  $\alpha (\ast ) = \sup \{\alpha_{B}\colon B \in  
{\cal P}\}$,  clearly  $\alpha (\ast ) < \lambda (2)^+$. 
So  $\bigwedge_{{B} \in {\cal P}}g_{B} <_{J} f_{\alpha (\ast )}$.
 Note:  ${\cal P} \subseteq 
{\cal P}'$  (as  ${\cal P} \in  \frB$,  $|{\cal P}| \leq  \lambda (2)$,  
$\lambda (2) + 1 \subseteq  \frB)$  and for each  $B \in  {\cal P}$, 
$\frc_{B} =: \{\sigma \in  \fra\colon g_{B}(\sigma ) \geq 
f_{\alpha(*)}(\sigma )\}$
is in  $J$ and\  $J \subseteq  J_{1} \subseteq  J_{2}$. 
Apply the choice of  ${\cal P}$  (i.e.\ it exemplifies  $\lambda (2))$ 
to  $A =:$\ Rang\ $f_{\alpha (\ast )}$, 
get  $\langle A_{n},{\cal P}_{n}\colon n < \omega \rangle $  as there. 
Let  $\fra_{n} =: \{\sigma \in  \fra\colon f_{\alpha (\ast )}(\sigma )
\in A_{n}\}$,
so  $\fra =  \bigcup_{n} \fra_{n}$,  hence for some  $m $, 
$\fra_{m}  \notin  J_{2}$ (as  $\fra \notin  J_{2}$,  $J_{2}$ is 
$\aleph_{1}$-complete) hence  $\fra_{m}  \notin  J_{1}$. 
As  $\fra\in\frB$, ${\cal P}\in\frB$ clearly
 ${\cal P}_{m}  \subseteq  \frB$. 
So $\{\frc_{B}\colon B\in {\cal P}_{m}\}$ is a family of $< \mu $ 
subsets of $\fra$, each in  $J$  and every countable  $\frb \subseteq 
\fra_{m} $ is included in at least one of them (as for some
$B\in \calP_m$, $\Rang ( f_{\alpha(\ast)}\upharpoonright
\frb)\subseteq B$, hence $\frb\subseteq \frc_B$). 
Easy contradiction.

\subdemo{$\lambda(3)\le\lambda(0)\hbox{ if }\cov(\theta, \aleph_1,
\aleph_1, 2)<\mu$}
Let $\fra \subseteq  \Reg \cap  \lambda ^+\backslash \mu $, $|\fra|\leq
\kappa$, let
$\langle\fra_{\ell}\colon\ell<\omega \rangle $  be as guaranteed by the
definition of  $\lambda (0)$, let ${\cal P}_{\ell} \subseteq  
\calS_{<\aleph_{1}}(\fra_{\ell} )$  exemplify  $\cov(\theta ,\aleph_{1},
\aleph_{1},2)\shvor<\mu$, for each $\frb\in \calP_{\ell} $ we can
find a finite  $\fre_{\frb} \subseteq  (\pcf \fra_{\ell} ) \cap 
\lambda ^+\backslash \mu $  such that  $\frb 
\subseteq  \bigcup \{\frb_{\sigma} [\fra_{\ell} ]\colon \sigma \in 
\fre_{\frb}\}$  and $\{\frb_{\ell ,i}\colon i < i^\ast \}$ 
enumerates  $\{\fre_{\frb}\colon \frb\in{\cal P}_{\ell} \}.$

\subdemo{$\lambda(0)\le\lambda(1)$}
Similar to the proof of  $\lambda (3) \leq  \lambda (2)$.              
\sqed{1.2}

\proclaim{1.3 Claim} 
Assume  $\aleph_{0} < \cf \lambda \leq  \theta < \lambda < 
\lambda ^\ast$,  $\pp(\lambda ) \leq \lambda ^\ast $ and
$$\cov(\lambda^* ,\lambda ^+,\theta^+, 2) < \lambda ^\ast. $$
Then  $\cov(\lambda ,\lambda ,\theta ^+,2) < \lambda ^\ast .$
\endproclaim

\demo{Proof} 
Easy.

\rem{1.3A Definition} 
Assume  $\lambda \geq  \theta = \cf \theta > \kappa = 
\cf \kappa > \aleph_{0}$.

(1) $(\bar C,\bar {\cal P}) \in  T^\oplus [\theta ,\kappa ]$  if 
$(\bar C,\bar {\cal P}) \in  {\cal T}^\ast [\theta ,\kappa ]$  (see
[Sh420, Def 2.1(1)]), and $\delta \in  S(\bar C) \Rightarrow  \delta = 
\sup (\acc \ C_{\delta} )$  (note: $\acc C_{\delta}  \subseteq 
C_{\delta} )$,  and we do not allow (viii$)^-$ (in [Sh420, Definition
2.1(1)]), or replace it by:\  

{
\itemindent=20pt
\item{(viii$)^*$} for some list  $\langle a_{i}\colon i < \theta
\rangle $
of  $ \bigcup_{{\alpha \in } S {(\bar C)}}{\cal P}_{\alpha} $, 
we have:  $\delta \in  S(\bar C)$,  $\alpha \in  \acc \ C_{\delta} $
implies  $\{a \cap  \beta \colon a \in  {\cal P}_{\delta} ,\beta 
\in  a \cap  \alpha \} \subseteq  \{a_{i}\colon i < \alpha \}.$

}

(2) For  $(\bar C,\bar {\cal P}) \in  {\cal T}^\oplus [\theta,
\kappa ]$  we define a\ filter  ${\cal D}^\oplus_{(\bar C,
\bar\calP)}(\lambda)$
on 
$[{\cal S}_{<\kappa }(\lambda )]^{<\kappa }$ (rather than on
${\cal S}_{<\kappa }(\lambda )$  as in [Sh420, 2.4]) (let  $\chi  = 
\beth _{\omega +1}(\lambda )):$

$Y \in  {\cal D}^\oplus_{(\bar C, \bar\calP)}(\lambda)$
iff $Y\subseteq\left( {\cal S}_{<\kappa }(\lambda )\right)^{<\kappa}$
and for some  $x \in  H(\chi )$  for every  $\langle N_{\alpha},
N^\ast_{a}\colon\shvor\alpha  < \theta$, $a \in  
 \bigcup_{{\delta} \in  {S}}{\cal P}_{\delta} \rangle $ 
satisfying condition  $\otimes$ from [Sh420, 2.4], and also 
$[a \in  {\cal P}_{\delta}  \und \delta \in  S \und \alpha 
< \theta \Rightarrow  x \in  N^\ast_{a} \und x \in  N_{\alpha} ]$  there
is  $A \in  \id^a(\bar C)$  such that  $\delta \in  S(\bar C)\backslash
A \Rightarrow \langle \bigcup_{{a} \in  {{\cal P}}_{\delta} }
N^\ast_{a} \cap  \lambda \cap  N_{\alpha} \colon \alpha \in  
\acc \ C_{\delta} \rangle  \in  Y.$

\rem{Remark} 
For 1.3B below, see Definition of  ${\cal T}^\ell (\theta ,\kappa )$  
and compare with [Sh420, Definition 2.1(2), (3)].

\proclaim{1.3B Claim}
\item{(1)} If  $(\bar C,\bar {\cal P}) \in  T^\oplus [\theta ,\kappa]$
(so $\lambda > \kappa $  are regular uncountable) then 
$D_{(\bar C, \bar\calP)}^\oplus(\lambda )$ is a non-trivial ideal on  
$[ {\cal S}_{<\kappa }(\lambda )] ^{<\kappa }$.

\item{(2)}  If  $\bar C \in  {\cal T}^0[\theta ,\kappa ]$,\  $[\delta
\in  S(\bar C) \Rightarrow  \delta = \sup  \acc \ C_{\delta} ]$, 
${\cal P}_{\delta}  = \{C_{\delta} \cap  \alpha \colon \alpha \in 
C_{\delta} \}$  then  $(\bar C,\bar {\cal P}) \in  
{\cal T}^\oplus [\theta ,\kappa ]$.  If  $\bar C \in  
{\cal T}^1[\theta ,\kappa ]$,  $[\delta \in  S(\bar C) \Rightarrow 
\delta = \sup  \acc \ C_{\delta} ]$  and ${\cal P}_{\delta}  = 
{\cal S}_{<\aleph_{0}}(C_{\delta} )$  then  $(\bar C,\bar {\cal P}) \in 
\calT^\oplus [\theta ,\kappa ].$

\item{(3)}  If  $\theta $  is successor of regular,  $\sigma = \cf
\sigma < \kappa $,  there is  $\bar C \in  {\cal T}^0[\theta ,\kappa]
\cap  {\cal T}^1[\theta ,\kappa ]$  with: for  $\delta \in  S(\bar C)$,
$C_{\delta} $ is closed,  $\cf \delta = \sigma $  and $\otp
C_{\delta} $ divisible by  $\omega ^2$
(hence  $\delta = \sup  \acc \ C_{\delta} ).$

\item{(4)}  Instead of $``\theta $  successor of regular", it suffices
to demand 
$$
\theta > \kappa \hbox{  regular uncountable, and }
 \bigwedge_{\alpha < \theta} 
\bigvee_{\kappa_1\in[\kappa,\theta)} \cov(\alpha ,\kappa_{1},\kappa ,2)
< \theta .
\leqno(*)
$$
Replacing 2 by  $\sigma $,  $``C_{\delta} $ closed"  is weakened to  
$``\{\otp (\alpha \cap  C_{\delta} )\colon \alpha \in  C_{\delta} \}$ 
is stationary''.
\endproclaim

\demo{Proof} 
Check.

\proclaim{1.3C Claim} 
Let  $\lambda > \kappa = \cf \kappa > \aleph_{0}$,  $\theta = 
\kappa ^+$,  $(\bar C,\bar {\cal P}) \in  {\cal T}^\oplus [\theta,
\kappa ]$  then the following cardinals are equal: 
{\itemindent=50pt
\item{$\mu(0)=$} $\cf \left( {\cal S}_{<\kappa }(\lambda ),\subseteq
\right)$, 

\item{$\mu(4)=$} $\Min \left\lbrace |Y|\colon Y \in 
{\cal D}^\oplus_{(\bar C, \bar\calP)}(\lambda )\right\rbrace .$
}
\endproclaim

\demo{Proof} 
Included in the proof of [Sh420, 2.6].

\proclaim{1.3D Claim} 
Let  $\lambda_{1} \geq  \lambda_{0} > \kappa = \cf \kappa > 
\aleph_{0}$,  $\theta = \kappa ^+$ and $(\bar C,\bar {\cal P}) \in  
{\cal T}^\oplus [\theta ,\kappa ]$.  Let  $\frB_{\lambda_{1}}$ be a
rich enough model with universe  $\lambda_{1}$ and countable vocabulary
which is rich enough $(e.g$. all functions (from  $\lambda_{1}$ to 
$\lambda_{1})$ \ definable in  $(H( \beth
_{\omega} (\lambda_{1})^+),\in ,<^\ast )$  with any finite number of
places). 
Then the following cardinals are equal: 

{\itemindent=50pt
\item{$\mu^*(0)=$} $\cov(\lambda_{1},\lambda ^+_{0},\kappa ,2)$, 

\item{$\mu^+(4)=$} $\Min \left\lbrace
|Y/\approx ^{\lambda_{0}}_{\frB_{\lambda_{1}}}|\colon
Y \in  {\cal D}^\oplus_{(\bar C, \bar\calP)}
(\lambda_{1}) \right\rbrace $  where    
$\langle a'_{i}\colon i \in  \acc \ C_{\delta} \rangle$
$\approx^{\lambda_{0}}_{\frB}$\break
$\langle a^{''}_i \colon i\in  \acc
C_{\delta}\rangle$ iff $\bigwedge_{{i\in ac} c {C_{\delta}} }$
Skolem Hull $_{\frB_{\lambda_{1}}}(a^{'}_i \cup  \lambda_{0}) =$\break
Skolem 
Hull $_{\frB_{\lambda_{1}}}({a'}_i \shvor \cup  \lambda_{0}).$

}
\endproclaim

\demo{Proof} 
Like the proof of [Sh420], 2.6, but using [Sh400, 3.3A].



\head 2. Equality Relevant to Weak Diamond \endhead

It is well known that:
$$
\kappa = \cf \kappa \und \theta > 2^{<\kappa } \Rightarrow  
\cov(\theta ,\kappa ,\kappa ,2) = \theta ^{<\kappa } = 
\cov(\theta ,\kappa ,\kappa ,2)^{<\kappa }.
$$
Now we have

\proclaim{2.1 Claim} 
\item{(1)} If  $\mu  > \lambda \geq  \kappa $,  $\theta = 
\cov(\mu ,\lambda ^+,\lambda ^+,\kappa )$,  $\cov(\lambda ,\kappa,
\kappa ,2) \leq  \mu$  (or  $\leq  \theta$)  then
$$
\cov(\mu ,\lambda ^+,\lambda ^+,2) = \cov(\theta ,\kappa ,\kappa ,2).
$$

\item{(2)}  If in addition  $\lambda \geq  2^{<\kappa }$ (or just 
$\theta \geq  2^{<\kappa })$  then
$$
\cov(\mu ,\lambda ^+,\lambda ^+,2)^{<\kappa } = 
\cov(\mu ,\lambda ^+,\lambda ^+,2).
$$
\endproclaim

\rem{2.1A Remark} 
\item{(1)} A most interesting case is  $\kappa = \aleph_{1}.$

\item{(2)}  This clarifies things in [Sh-f,AP1.17].

\demo{Proof} 
(1) Note that  $\theta \geq  \mu $  (because  $\mu  > \lambda \geq  
\kappa )$. 
First we prove  $``\leq "$. Let  ${\cal P}_{0}$ be a family of  
$\theta $  subsets of  $\mu $  each of cardinality  $\leq  \lambda $, 
such that every subset of  $\mu $  of cardinality  $\leq  \lambda $ 
is included in the union of $< \kappa $  of them (exists by the
definition of  $\theta = 
\cov(\mu ,\lambda ^+,\lambda ^+,\kappa ))$. 
Let  ${\cal P}_{0} = \{A_{i}\colon i < \theta \}$. 
Let  ${\cal P}_{1}$ be a family of  $\cov(\theta ,\kappa ,\kappa ,2)$
subsets of  $\theta $,  each of cardinality  $< \kappa $  such that any
subset of  $\theta $  of cardinality  $< \kappa $  is included in\ one
of them.

Let  ${\cal P} =: \left\{ \bigcup_{{i} \in  {a}}A_{i}\colon a
\in  {\cal P}_{1}\right\}$; clearly  ${\cal P}$  is a family of 
subsets of  $\mu $  each of cardinality  $\leq  \lambda $,  $|{\cal P}| 
\leq  |{\cal P}_{1}| = \cov(\theta ,\kappa ,\kappa ,2)$,  and every  $A 
\subseteq \mu $,  $|A| \leq  \lambda $  is included in some union of
$< \kappa $ 
members of  ${\cal P}_{0}$ (by the choice of  ${\cal P}_{0})$,  say 
$ \bigcup_{{i} \in  {b}}A_{i}$,  $b \subseteq  \theta $, 
$|b| < \kappa $;  by the choice of  
${\cal P}_{1}$,  for some  $a \in  {\cal P}_{1}$ we have 
$b \subseteq  a$,  hence $A \subseteq \bigcup_{{i} \in  {b}}
A_{i} \subseteq \bigcup_{{i} \in  {a}}A_{i} \in  {\cal P}$. 
So  ${\cal P}$  exemplify  $\cov(\mu ,\lambda ^+,\lambda ^+,2) \leq 
\cov(\theta ,\kappa ,\kappa ,2).$ 

Second we prove the inequality  $\geq.$ 
If  $\kappa \leq  \aleph_{0}$ then  
$\cov(\mu ,\lambda ^+,\lambda ^+,2) = \theta $  and 
$\cov(\theta,\kappa,\kappa,2)=\theta$ so $\geq$ trivially holds;
so assume  $\kappa > \aleph_{0}$. 
Obviously  $\cov(\mu ,\lambda ^+,\lambda ^+,2) \geq  \theta $. 
Note, if  $\kappa $  is singular then, as  $\cf \lambda ^+ > 
\lambda \geq  \kappa $  for some  $\kappa_{1} < \kappa$, we have 
$\theta = \cov(\mu ,\lambda ^+,\lambda ^+,\kappa ) = 
\cov(\mu ,\lambda ^+,\lambda ^+,\kappa ')$  whenever  $\kappa ' \in  
[\kappa_{1},\kappa ]$  is a successor (by [Sh355, 5.2(8)]); also  
$\cov(\theta ,\kappa ,\kappa ,2) \leq  
\sup \{\cov(\theta ,\kappa ,\kappa ',2)\colon \kappa ' \in  
[\kappa_{1},\kappa ]$\ \ is\ a\ successor\ cardinal$\}$  and 
$\cov(\theta ,\kappa ,\kappa ',2) \leq  \cov(\theta ,\kappa ',\kappa ',2)$
when $\kappa ' < \kappa$, so without loss of generality  $\kappa $ 
is regular uncountable.  Hence for any  $\theta_{1} < \theta $  we have 

\item{$(*)_{\theta_{1}}$} we can find a family 
${\cal P} = \{A_{i}\colon i
< \theta_{1}\}$, $A_{i} \subseteq  \mu $,  $|A_{i}| \leq  \lambda $, 
such that any subfamily of 
cardinality  $\leq  \lambda ^+$ has a transversal.    
[Why?  By [Sh355, 5.4], $(=^+)$ and [Sh355,1.5A] even for $\leq \mu$.]

\medskip
Hence if  $\theta_{1} \leq  \theta $,  $\cf \theta_{1} < \lambda ^+$
(or even  $\cf \theta_{1} \leq  \mu )$  then
$(\ast )_{\theta_{1}}$.  Now we shall prove below 
$$
(*)_{\theta_{1}} \Rightarrow  \cov(\theta_{1},\kappa ,\kappa ,2) \leq 
\cov(\mu ,\lambda ^+,\lambda ^+,2)
\leqno(\otimes_{1})
$$
and obviously 
$$
\hbox{ if }  \cf \theta \geq  \kappa \hbox{  then  }
\cov(\theta ,\kappa ,\kappa ,2) = 
\sum_{\alpha <\theta }\cov(\alpha ,\kappa ,\kappa ,2)
\leqno(\otimes_{2})
$$
together; (as  $\theta \leq  \cov(\theta ,\lambda ^+,\lambda ^+,2)$
which holds as $\lambda<\mu\le\theta$) we are done.

\demo{Proof of $\otimes_1$}
Let  $\{A_{i}\colon i < \theta_{1}\}$  exemplify $(\ast )_{\theta_{1}}$ and
${\cal P}_{2}$ exemplify the value of  $\cov(\mu ,\lambda ^+,\lambda ^+,
2)$.  Now for every  $a \subseteq 
\theta_{1}$,  $|a| < \kappa $,  let  $B_{a} =: \bigcup_{{i} \in 
{a}}A_{i}$;  so  $B_{a} \subseteq \mu $,  $|B_{a}| \leq  \lambda$
hence there is\  $A_{a} \in  {\cal P}_{2}$ such that:  $B_{a} 
\subseteq A_{a}.$
Now for  $A \in  {\cal P}_{2}$ define  $b[A] =: \{i < \theta_{1}\colon
A_{i} 
\subseteq A\}$;  it has cardinality  $\leq  \lambda $  (as any
subfamily of  $\{A_{i}\colon A_{i} \subseteq 
A\}$  of cardinality  $\leq  \lambda ^+$ has a transversal).  Note  $a 
\subseteq b[A_{a}]$  (just read the definitions of  $b[A]$  and
$A_{a}$;  note  $a \in  {\cal S}_{<\kappa }(\theta_{1}))$. 
For  $A \in  {\cal P}_{2}$ let  ${\cal P}_{A}$ 
be a family of $\leq  \cov(\lambda ,\kappa ,\kappa ,2)$  subsets of 
$b[A]$  each of cardinality  $< \kappa $  such that any such set is
included in one of 
them (exists as  $|b[A]| \leq  \lambda $  by the definition of  
$\cov(\lambda ,\kappa ,\kappa ,2))$.  So for any  $a \in  
{\cal S}_{<\kappa }(\theta_{1})$  for some  $c \in  {\cal P}_{A_{a}}$, 
$a \subseteq c$. 
We can conclude that  $\bigcup \{\calP_A\colon A \in 
\calP_2\}$  is a
family exemplifying  $\cov(\theta_{1},\kappa ,\kappa ,2) \leq  
\cov(\mu ,\lambda ^+,\lambda ^+,2) + \cov(\lambda ,\kappa ,\kappa ,2)$ 
but the last term is $\leq\mu$  (by an assumption) whereas the
first is $\ge\mu$ (as $\mu>\lambda$) hence the second term is redundant.

(2)  By the first part it is enough to prove  
$\cov(\theta ,\kappa ,\kappa ,2)^{<\kappa } = \cov(\theta ,\kappa,
\kappa ,2)$,  which is easy and well known (as $\theta\ge\mu>\lam
\ge 2^{<\kappa}$).
\sqed{2.1}

\rem{2.1B Remark} 
So actually if  $\mu  > \lambda \geq  \kappa $,  $\theta = 
\cov(\mu ,\lambda ^+,\lambda ^+,\kappa )$  then ($\theta\ge\mu>
\lam\ge \kappa$ and)
$$
\eqalign{
\cov(\mu ,\lambda ^+,\lambda ^+,2) & \leq  
\cov(\mu ,\lambda ^+,\lambda ^+,\kappa ) + \cov(\theta ,\kappa ,\kappa,
2) \cr
& = \theta + \cov(\theta ,\kappa ,\kappa ,2) = \cov(\theta ,\kappa,
\kappa ,2)\cr}
$$
and
$$ 
\cov(\theta ,\kappa ,\kappa ,2)  \leq  \cov(\mu ,\lambda ^+,\lambda^+,
2) + \cov(\lambda ,\kappa ,\kappa ,2),
$$

\noindent
hence, $\cov(\theta,\kappa,\kappa ,2) = \cov(\mu,\lambda ^+,\lambda^+,
2) + \cov(\lambda ,\kappa ,\kappa ,2).$



\head 3. Cofinality of $S_{\le\aleph_0}(\kappa)$ for $\kappa$
Real Valued Measurable and Trees \endhead

In Rubin--Shelah [RuSh117] two covering properties were discussed
concerning partition theorems on trees, the stronger one was sufficient,
the weaker one necessary so it was asked whether they are equivalent. 
[Sh371, 6.1, 6.2] gave a partial positive answer (for  $\lambda $ 
successor of regular, but then it 
gives a stronger theorem); here we prove the equivalence. 

In Gitik--Shelah [GiSh412] cardinal arithmetic, e.g.\ near a real valued
measurable cardinal $\kappa $,  was investigated, e.g.\ 
$\{2^\sigma\colon
\sigma < \kappa \}$  is finite (and more); this section continues it. 
In particular we answer a problem of Fremlin: 
for $\kappa$ real valued measurable, do we have
$\cf ({\cal S}_{<\aleph_{1}}(\kappa ),\subseteq)=\kappa$?
Then we deal with trees with many branches; on earlier theorems 
see [Sh355, {\S}0], and later [Sh410, 4.3].

\proclaim{3.1 Theorem} 
Assume  $\lambda $, $\theta $, $\kappa $  are regular cardinals 
and $\lambda > \theta = \kappa > \aleph_{0}$. 
Then the following conditions are equivalent: 

\item{(A)} for every  $\mu  < \lambda $  we have  $\cov(\mu ,\theta,
\kappa ,2) < \lambda $, 

\item{(B)} if  $\mu  < \lambda $  and $a_{\alpha}  \in 
{\cal S}_{<\kappa }(\mu )$  for  $\alpha < \lambda $  then for some  $W 
\subseteq \lambda $  of cardinality  $\lambda $  we have 
$| \bigcup_{{\alpha} \in  {W}}a_{\alpha} | < \theta .$
\endproclaim

\rem{3.1A Remark} 
(1) Note that (B) is equivalent to: if  $a_{\alpha}  \in  
{\cal S}_{<\kappa }(\lambda )$  for  $\alpha < \lambda $,  then for some
unbounded  $S \subseteq  \{\delta < \lambda \colon \cf (\delta ) \geq 
\kappa \}$  and $b \in  {\cal S}_{<\theta }(\lambda )$,  for  $\alpha
\neq  \beta $  in  $S$,  $a_{\alpha}  \cap  a_{\beta}  \subseteq  b$ 
(we can start with any stationary  $S_{0} \subseteq  \{\delta <
\lambda \colon \cf \delta \geq  \kappa \}$,  and use Fodour Lemma).

\noindent
(2)  We can replace everywhere  $\theta $  by  $\kappa $,  but want to
prepare for a possible generalization.
By the proof we can strengthen ``$W\subseteq\lambda$ of
cardinality $\lambda$'' to ``$W\subseteq\lambda$ is stationary''
(for $\lnot({\rm A})\to \lnot({\rm B})$ this is trivial, for
$({\rm A})\rightarrow({\rm B})$ real), so these two versions of
$({\rm B})$ are equivalent.

\demo{Proof} 

\subdemo{(A)$\Rightarrow$(B)}

Trivial [for  $\mu  < \lambda $  let  ${\cal P}_{\mu}  
\subseteq 
{\cal S}_{<\theta }(\mu )$  exemplify  $\cov(\mu ,\theta ,\kappa ,2) < 
\lambda $;  suppose  $\mu  < \lambda $  and $a_{\alpha}  \in  
{\cal S}_{<\kappa }(\mu )$  for  $\alpha < \lambda $  are given,  for
each  $\alpha $  for some  $A_{\alpha}  \in  {\cal P}_{\mu} $ we have 
$a_{\alpha}  \subseteq A_{\alpha} $;  as  $|{\cal P}_{\mu} | <
\lambda = \cf \lambda $  for some  $A^\ast $
we have  $W =: \{\alpha < \lambda \colon A_{\alpha}  = A^\ast \}$  has
cardinality  $\lambda $,  so  $S$  is as required in (B)]. 

\subdemo{$\lnot {\rm (A)} \Rightarrow  \lnot {\rm (B)}$}

\subdemoinfo{First Case}{For some  $\mu  \in  [\theta ,\lambda )$, 
$\cf \mu  < \kappa < 
\mu $  and $pp^+_{<\kappa }(\mu ) > \lambda$} 
Then we can find $\fra \subseteq 
\Reg\cap\mu\backslash\theta$, $|\fra| < \kappa$, $\sup\fra=\mu$
and $\max  \pcf_{J^{bd}_{\fra}}\fra \geq  \lambda $. 
So by [Sh355, 2.3] without loss of generality  $\lambda=\max\pcf\fra$;
let  $\langle f_{\alpha} \colon \alpha < \lambda \rangle$ be
$<_{J_{<\lambda }[\fra]}$-increasing cofinal in  $\prod \fra$. 

Let  $a_{\alpha}  =\Rang(f_{\alpha} )$,  so for  $\alpha < \lambda $, 
$a_{\alpha} $ is a subset of  $\mu  < \lambda $  of cardinality 
$< \kappa $.  Suppose  $W \subseteq 
\lambda $  has cardinality  $\lambda $,  hence is unbounded, and we
shall show 
that  $\mu  = | \bigcup_{{\alpha} \in  {W}}a_{\alpha} |$;  as 
$\mu  \geq  \theta $  this is enough. 
Clearly  $a_{\alpha}  = \Rang\ f_{\alpha}  
\subseteq 
\sup  \fra = \mu $,  hence \ $ \bigcup_{{\alpha} \in  {W}}a_{\alpha} 
\subseteq \mu $. 
If  $| \bigcup_{{\alpha} \in  {W}}a_{\alpha} | < \mu $  define 
$g \in  \prod \fra$  by:  $g(\sigma )$  is  
$\sup \left( \sigma \cap \bigcup_{{\alpha} \in  {W}}a_{\alpha}
\right) $  if  $\sigma >\vert\bigcup_{{\alpha} \in  {W}}a_{\alpha}\vert$
and 0  otherwise.  So  $g \in  \prod \fra$  hence for\ some  $\beta <
\lambda \ \ g < f_{\beta} $\ mod\ $J_{<\lambda }[\fra]$. 
As the  $f_{\beta}$'s are $<_{J_{<\lambda }[\fra]}$-increasing
and $W \subseteq \lambda $  unbounded, without loss of generality 
$\beta \in  W$, hence by  
$g$'s  choice  $\left[ \sigma \in  \fra\backslash |
\bigcup_{{\alpha} \in  {W}}a_{\beta} |^+ \Rightarrow  f_{\beta}
(\sigma ) \leq  g(\sigma )\right] $  but  
$\left\lbrace \sigma \colon \sigma\in\fra,\sigma > | \bigcup_{{\theta}
\in  {W}}a_{\alpha} |^+\right\rbrace  \notin  J_{<\lambda }[\fra]$ 
(as $\mu$ is a limit cardinal and $\max\pcf_{J^{bd}_{\fra}}(\fra)\geq 
\lambda )$,  contradiction.

The main case is:

\subdemoinfo{Second Case}{For no  $\mu  \in  [\theta ,\lambda )$ 
is  $\cf \mu  < \kappa < \mu $,  $pp^+_{<\kappa }(\mu ) > \lambda$}
Let  $\chi  =: \beth
_{2}(\lambda )^+$,  $\frB$ be the model with universe $\lambda$ and the
relations and functions definable in  $(H(\chi ),\in ,<^\ast_{\chi} )$ 
possibly with the parameters  $\kappa ,\theta ,\lambda $. 
We know that  $\lambda > 
\theta ^+$ (otherwise  $\lambda = \theta ^+$ and (A) holds). 
Let  $S \subseteq \{\delta < \lambda \colon \cf \delta = \theta \}$ 
be stationary and in  $I[\lambda ]$
(see [Sh420, 1.5]) and let  $S \subseteq 
S^+$,  $\bar C = \langle C_{\alpha} \colon \alpha \in  S^+\rangle $ 
be such that:  
$C_{\alpha} $ closed, $\otp C_{\alpha}  \leq  \theta $,  $[\beta \in 
\nacc \ C_{\alpha}  \Rightarrow  C_{\beta}  = C_{\alpha}  \cap  \beta]$,
$[\otp \ C_{\alpha}  = \kappa \Leftrightarrow  \alpha \in  S]$ 
and for $\alpha\in S^+$ limit, $C_\alpha$ is unbounded in $\alpha$
(see [Sh420, 1.2]).

Without loss of generality  $\bar C$  is definable in  
$(\frB,\kappa,\theta,\lambda)$. 
Let  $\mu_{0} \in  [\theta ,\lambda )$ be minimal such that 
$\cov(\mu_{0},\theta ,\kappa ,2) \geq  \lambda $,  
so  $\mu_{0} > \theta $,  $\kappa > \cf \mu_{0}$. 
We choose by induction on  $\alpha < 
\lambda $,  $\frA_{\alpha} $,  $a_{\alpha} $ such that:

\item{$(\alpha)$} $\frA_{\alpha} \prec (H(\chi ),\in ,<^\ast_{\chi} )$, 
$\|\frA_{\alpha} \| <
\lambda $  and $\frA_{\alpha}  \cap  \lambda $  is an ordinal and\break
$\{\lambda ,\mu_{0},\theta ,\kappa ,\frB,\shvor\bar C\} \in 
\frA_{\alpha} .$ 

\item{$(\beta )$} $\frA_{\alpha} (\alpha < \lambda )$  is increasing
continuous and $\langle\frA_{\beta}\colon\beta\leq\alpha\rangle\in 
\frA_{\alpha +1}.$ 

\item{$(\gamma )$} $a_{\alpha}  \in  {\cal S}_{<\kappa }(\mu_{0})$ 
is such that for no $A \in  {\cal S}_{<\theta }(\mu_{0}) \cap 
\frA_{\alpha} $ is  $a_{\alpha}  \subseteq  A$. 

\item{$(\delta )$} $\langle a_{\beta} \colon \beta \leq  \alpha \rangle 
\in  \frA_{\alpha +1}$.

There is no problem to carry the definition and let  $\frA =
\bigcup_{\alpha< \lambda }\frA_{\alpha} $. 
Clearly it is enough to show that  $\bar a = 
\langle a_{\alpha} \colon \alpha < \lambda \rangle $  contradict (B). 
Clearly  $\mu_{0} \in  (\theta ,\lambda )$  and $a_{\alpha}  \in  
{\cal S}_{<\kappa }(\mu_{0})$.  So let  $W 
\subseteq 
\lambda $,  $|W| = \lambda $  and we shall prove that  $|
\bigcup_{{\alpha} \in  {W}}a_{\alpha} | \geq  \theta $.  Note: 
\item{$(*)$} if $\fra\subseteq [\theta,\lambda)$, $|\fra| < \kappa $,
$\fra \in  \frA_{\alpha} $ (and $\fra \subseteq  \Reg$, of course) then 
$(\prod \fra) \cap  \frA_{\alpha} $
is cofinal in  $\prod \fra$  (as $\max\pcf  \fra < \lambda )$.

\noindent
Let $R=\{(\alpha,\beta)\colon\beta\in a_{\alpha},\alpha<\lambda\}$
and 
$$
E = : \left\lbrace \delta < \lambda \colon (\frA_{\delta} ,R
\upharpoonright \delta ,W \cap
\delta ,\mu_{0}) \prec (\frA,R,W,\mu_{0})\right. \ \hbox{ and }
\ \frA_{\delta}  \cap  \lambda = \delta \left. \right\rbrace .
$$
Clearly $E$ is a club of $\lambda $, hence we can find $\delta (\ast ) 
\in  S \cap  \acc (E)$. 
Let  $C_{\delta (\ast )} =$\break
$\{\gamma_{i}\colon i < \theta \}$  (in increasing order). 
We now define by induction on $n<\omega$, $M_{n}$,\break
$\langle N^n_{\zeta}\colon\zeta<\theta \rangle $,  $f_{n}$ such that: 

\item{(a)} $M_{n}$ is an elementary submodel of  $(\frA,R,W)$, 
$\|M_{n}\| = \theta$, 

\item{(b)} $\langle N^n_{\zeta} \colon \zeta < \theta \rangle $  is an
increasing continuous sequence of elementary submodels of  $\frB$,

\item{(c)} $\|N^n_{\zeta} \| < \theta $,

\item{(d)} $N^n_{\zeta}  \in  \frA_{\delta (\ast )}$, 

\item{(e)} $ \bigcup_{{\zeta} < {\kappa}} |N^n_{\zeta} | 
\subseteq |M_{n}|$,  

\item{(f)} $f_{n} \in  \prod (\Reg \cap  M_{n})$,

\item{(g)} $f_{n}(\sigma ) > \sup (M_{n} \cap  \sigma )$  for 
$\sigma \in  \Dom(f_{n})\backslash \theta^+$,

\item{(h)} for every $\zeta<\theta$, $f_{n}\upharpoonright (\Reg \cap  
N^n_{\zeta} \backslash \theta ^+) \in  \frA_{\delta (\ast )}$, 

\item{(i)} $N^0_{\zeta} $ is the Skolem Hull in  $\frB$  of  $\{\gamma_{i},
i\colon i < \zeta \}$,

\item{(j)} $N^{n+1}_{\zeta} $ is the Skolem Hull in $\frB$ of
$N^n_{\zeta}  \cup  \{f_{n}(\sigma )\colon \sigma \in  \Reg \cap 
N^n_{\zeta} \backslash \theta ^+\}$,

\item{(k)} $M_{n}$ is the Skolem Hull in  $(\frA,R,W)$  of 
$ \bigcup_{{\ell} < {n}}M_{\ell}  \cup \bigcup_{{\zeta} <
{\theta}} N^n_{\zeta} .$ 

There is no problem to carry the definition: for  $n = 0$  define 
$N^0_{\zeta} $ by (i) [trivially (b) holds and also (c), as for
(d), note that  $\bar C \in  \frA_{0} \prec \frA_{\delta (\ast )}$
and\break $\{\gamma_{i}\colon i < \zeta \} \in 
\frA_{\delta (\ast )}$ as  $\bar C$  is definable in  $\frB$  hence  
$\{\langle \alpha ,\gamma ,\zeta \rangle \colon \alpha \in  S^+,\zeta < 
\theta $,  and  $\gamma $  is the $\zeta$-th member of
$C_{\alpha} \}$  is a relation of  $\frB$  hence each  $C_{\gamma_{\zeta
+1}}(\zeta < \theta )$  is in  $\frA_{\delta (\ast )}$ hence each 
$\{\gamma_{i}\colon i < \zeta \}$  is and we can compute the Skolem
Hull in $\frA_{\gamma_{j}}$ for  $j < \theta $  large enough].

Next, choose  $M_{n}$ by  (k),  it satisfies (e) $+$ (a).  If  
$\langle N^n_{\zeta} \colon \zeta < \theta \rangle $,  $M_{n}$ are
defined, we can find $f_{n}$ satisfying  (f) $+$ (g) $+$ (h) by
[Sh371,1.4] (remember $(\ast ))$. 
For $n + 1$ define $N^n_{\zeta}$ by (j) and then $M_{n+1}$ by (k). 

Next by [Sh400, 3.3A or 5.1A(1)] we have 
$$
\bigcup_{{n} < {\omega}} M_{n} \cap  \delta (\ast ) =
\bigcup_{{n<\omega}\atop{\zeta<\theta}}N^n_{\zeta}  
\cap \delta(*)
\quad\hbox{ hence } \bigcup_{{n<\omega}\atop{\zeta<\theta}}
N^n_\zeta\cap W\ \hbox{ is unbounded in } \delta(*),
\leqno(*)
$$
hence for some  $n$ 
$$
\bigcup_{\zeta<\theta} N^n_{\zeta}  \cap  W \ \hbox{ is unbounded
in  }\delta(*).
\leqno{(*)_n}
$$

Remember  $N^n_{\zeta}  \in  \frA_{\delta (\ast )} =
\bigcup_{{\alpha <} \delta {(\ast )}}\frA_{\alpha}  = \bigcup_{{i}
< {\theta}} \frA_{\gamma_{i}}$.
So for some club $e$ of $\theta$ we have: 
$$
\hbox{ if } \zeta \in  e,\quad  \xi  < \zeta \ \hbox{ then: }
N^n_{\xi} 
\in  \frA_{\gamma_{\zeta} }, \hbox{ and }
\gamma_{\zeta}  \in  E \cap  C_{\delta (*)}
\leqno(\otimes) 
$$
(remember  $\delta (\ast ) \in  \acc (E)).$

Hence, for  $\zeta \in  e$,  we have:
$\frA_{\gamma_{\zeta} } \cap  \lambda = \gamma_{\zeta} $,  and $W \cap  
N^n_{\zeta} \backslash \sup  N^n_{\xi}  \neq  \emptyset $  for every 
$\xi  < \zeta $.  Let  $e = \{\zeta (\epsilon )\colon \epsilon  < \theta \}$, 
$\zeta (\epsilon )$  strictly increasing continuous in  $\epsilon $. 
Now for every  $\epsilon  < \theta $,  $N^n_{\zeta (\epsilon )} \in  
\frA_{\gamma_{\zeta (\epsilon +1)}}$ (and $\langle a_{\beta} \colon
\beta \leq  
\sup (\lambda \cap  N^n_{\zeta (\epsilon )}\rangle  \in  
\frA_{\gamma_{\zeta (\epsilon +1)}})$  hence  $A_{1} =: \bigcup
\{a_{\beta} \colon
\beta \in W \cap  N^n_{\zeta (\epsilon )}\} \subseteq  A_{2} =: \bigcup
\{a_{\beta} \colon \beta \in N^n_{\zeta (\epsilon +1)}\} \cap  \mu_{0} \in  
\frA_{\gamma_{\zeta (\epsilon +1)}}$ and $A_{2}$ is a subset of  $\mu_{0}$
of cardinality  $< \theta $  hence (by the choice of the
$a_\gamma$'s above)
$a_{\gamma_{\zeta(\epsilon+1)}}\not\subseteq A_{2}$ hence  
$a_{\gamma_{\zeta (\epsilon +1)}} \not\subseteq \bigcup
\{a_{\beta} \colon \beta \in 
W \cap  N^n_{\zeta (\epsilon )}\}$;  moreover, similarly  
$\gamma_{\zeta (\epsilon +1)} \leq  \gamma < \lambda \Rightarrow 
a_{\gamma}  
{\miss} \bigcup \{a_{\beta} \colon \beta \in  W \cap 
N^n_{\zeta (\epsilon )}\}.$

But  $W \cap  N^n_{\zeta (\epsilon +2)}\backslash \gamma_{\zeta
(\epsilon +1)}
\neq  \emptyset $,  hence  $\langle \bigcup \{ a_{\beta} \colon \beta \in  W \cap 
N^n_{\zeta (\epsilon )}\}\colon \epsilon  < \theta \rangle $  is
not eventually constant, hence 
$$
\bigcup \left\{ a_{\beta} \colon \beta \in
 W \cap \bigcup_{{\epsilon} < {\theta}} N^n_{\zeta (\epsilon )}
\right\}  = \bigcup \left\{ a_{\beta} \colon \beta \in  W \cap
\bigcup_{{\zeta} < {\theta}} N^n_{\zeta} \right\}
$$
has cardinality  $\theta$. 
Hence  $ \bigcup_{{\beta} \in  {W}}a_{\beta} $ has cardinality 
$\geq  \theta $,  as required. 
\sqed{3.1}

\rem{3.2 Conclusion} 
(1) If  $\lambda $  is real valued measurable then  $\kappa = 
\cf \left[ {\cal S}_{<\aleph_{1}}(\lambda ), \subseteq 
\right] $  (equivalently,  $\cov(\lambda ,\aleph_{1},\aleph_{1},2) =
\lambda ).$

\noindent
(2)  Suppose  $\lambda $  is regular  $> \kappa = \cf \kappa >
\aleph_{0}$, $I$ is a $\lambda $-complete ideal on  $\lambda$
extending $J^{bd}_{\lambda} $ and
is $\kappa$-saturated (i.e.\ we cannot partition  $\lambda $  to 
$\kappa $  sets not in  $I$).  Then for  $\alpha < \lambda $,  
$\cf ({\cal S}_{<\kappa }(\alpha ),\subseteq ) < \lambda $, 
equivalently  $\cov(\alpha ,\kappa ,\kappa ,2) < \lambda .$

\rem{3.2A Remark} 
(1) So for regular  $\theta \in  (\kappa ,\lambda )$  (in the 
above situation) we have  $ \bigwedge_{{\alpha} < {\lambda}}
\cov(\alpha ,\theta ,\theta ,2) < \lambda $; actually  $\kappa \leq  
\cf \theta \leq  \theta < \lambda $  suffices by the proof.

\demo{Proof} 
(1) Follows by (2).

\noindent
(2)  The conclusion is (A) of Theorem 3.1, hence it suffices to
prove (B).
Let  $\mu  < \lambda $  and $a_{\alpha}  \in  {\cal S}_{<\kappa }(\mu )$
for  $\alpha < \lambda $  be given. 
As  $\kappa < \lambda = \cf \lambda $  without loss of generality for
some  $\sigma < \kappa $,  $ \bigwedge_{{\alpha} < {\lambda}}
|a_{\alpha} | = \sigma $. 
Let  $f_{\alpha} $ be a function from  
$\sigma $  onto  $a_{\alpha} $,  so  $\Rang f_{\alpha}  
\subseteq \mu $. 
Now for each  $i < \sigma$, $\langle\{\alpha <\lambda\colon f_{\alpha}
(i) =
\gamma \}\colon \gamma <\mu \rangle$ is a partition of $\lambda $ to $\mu $  
sets; as $I$ is $\kappa $-saturated, $b_{i} =: \left\lbrace \gamma < 
\mu \colon \{\alpha < \lambda \colon f_{\alpha} (i) = \gamma \} \notin 
I\right\rbrace $  has cardinality  $< \kappa $,  hence 
$b =: \bigcup_{{i} < {\sigma}} b_{i}$ has cardinality 
$< \kappa + \sigma ^+ \leq  \kappa $  (remember 
$\sigma < \kappa = \cf \kappa )$.  For each  $i < \sigma $,  $\gamma \in
\mu \backslash b_{i}$ the set  $\{\alpha < \lambda \colon f_{\alpha} (i) =
\gamma \}$  is in  $I$;  so as  $I$  is $\lambda $-complete, 
$\lambda > \mu $  we have:  
$\{\alpha < \lambda \colon f_{\alpha} (i) \notin  b_{i}\}$  is in $I$. 
Now let
$$
W =: \{\alpha < \lambda \colon \hbox{ for some } i < \sigma ,f_{\alpha}
(i) \notin  b_i\} \subseteq  \bigcup_{{i} < {\sigma}}
\{\alpha < \lambda \colon f_{\alpha} (i) \notin  b_{i}\}.
$$
This is the union of $\leq\sigma<\lambda$  sets each in $I$,
hence is in $I$,  so  $|\lambda \backslash W| = \lambda $,  and clearly 
$$
\bigcup_{{\alpha \in \lambda } \backslash  {W}}a_{\alpha}  =
\{f_{\alpha} (i)\colon \alpha \in  \lambda \backslash W,i < \sigma \} 
\subseteq  \{f_{\alpha} (i)\colon \alpha < \lambda ,\neg f_{\alpha} (i)
\notin  b_{i},i < \sigma \} \subseteq  b,
$$ 
and $|b| < \kappa $  so 
$\lambda \backslash W$  is 
as required in (B) of Theorem 3.1. 
\sqed{3.2}

\proclaim{3.3 Lemma}
For every  $\lambda $  there is  $\mu $,  $\lambda \leq  \mu  < 
2^\lambda $ such that {\rm (A)} or {\rm (B)} or {\rm (C)} below holds
(letting 
$\kappa = \Min \{\theta \colon 2^\theta = 2^\lambda \}):$\  

\item{(A)} $\mu=\lambda$ and for every regular $\chi\leq 2^\lambda $
there is a tree  $T$  of cardinality  $\leq  \lambda $  with  $\geq  
\chi\ \cf(\kappa)$-branches (hence there is a linear order of
cardinality $\geq  \chi $  and density  $\leq  \lambda ).$ 

\item{(B)} $\mu  > \lambda $  is singular,  and:  

\itemitem{$(\alpha)$} $\pp (\mu)  = 2^\lambda $ (even  $\lambda =
\kappa \Rightarrow \pp^+(\mu ) =
\left( 2^\lambda \right) ^+)$,  $\cf \mu  \leq  \lambda $,  
$(\forall \theta )[\cf \theta \leq  \lambda < \theta < \mu  \Rightarrow 
\pp_{\lambda} \theta < \mu ]$  (and $\mu  \leq  2^{<\kappa })$

\noindent
hence  

\itemitem{$(\alpha)'$} for every successor\footnote*{If $\lambda =
\kappa$, just regular, and we can change $\lambda$ for
this.} $\chi  \leq  2^\lambda $
there is a tree from [Sh355, 3.5]:  $\cf \mu $  levels, every level of
cardinality  $< \mu $  and $\chi$  $(\cf \mu )$-branches,  

\itemitem{$(\beta)$} for every  $\chi  \in  (\lambda ,\mu )$,  there
is a tree  $T$  of cardinality  $\lambda $  with  $\geq  \chi $
branches of the same height, 

\itemitem{$(\gamma)$} $\cf\mu\geq\cf\kappa$ and even $\cf\kappa>
\aleph_{0} 
\Rightarrow  pp_{\Gamma (\cf \mu )}(\mu ) =^+ 2^\lambda .$ 

\item{(C)} Like (B) but we omit $(\alpha)$ and retain $(\alpha)'.$
\endproclaim

\demo{Proof}
\subdemoinfo{First Case}{$\kappa = \aleph_{0}$} 
Trivially (A) holds.

\subdemoinfo{Second Case}{$\kappa$ is regular uncountable} 
So  $\kappa \leq  \lambda $  and $2^\kappa = 2^\lambda $ and $[\theta < 
\kappa \Rightarrow  2^\theta < 2^\kappa ]$  hence  $2^{<\kappa } <
2^\kappa $ (remember  $\cf (2^\kappa ) > \kappa ).$
Try to apply [Sh410, 4.3], its assumptions (i) + (ii) hold (with 
$\kappa $  here standing for  $\lambda $  there) and if possibility (A)
here fails then the assumption (iii) there holds, too; so there is 
$\mu $  as there;  so $(\alpha )$, $(\gamma )$ of (B) of 3.3
holds\footnote{**}{Alternatively to quoting [Sh410, 4.3], we can get this
directly, if\break $\cov(2^{<\kappa },\lambda ^+,(\cf \kappa )^+,$
$\cf \kappa) < 2^\lambda $ we can get 
(A); otherwise by [Sh355, 5.4]  for some 
$\mu_{0} \in 
(\lambda ,2^{<\kappa }]$,  $\cf (\mu_{0}) = \cf \kappa $  and
$\pp(\mu_{0}) = (2^\lambda )$. 
Let  $\mu  \in  (\lambda ,2^{<\kappa }]$  be minimal such 
that $\cf\mu\leq\lambda\und\pp_{\lambda} (\mu ) > 2^{<\kappa }$. 
Necessarily ([Sh355, 2.3] and [Sh371, 1.6(2), (3), (5)]) 
$\pp_{\lambda} (\mu ) = \pp\mu  =
\pp(\mu_{0})=(2^\lambda)$ and (again using [Sh355, 2.3]) we have 
$(\forall \theta )[\cf \theta \leq  
\lambda < \theta < \mu  \Rightarrow  \pp_{\lambda} (\theta ) < \mu ]$; 
together $(\alpha )$ of (B) holds. 
Also  $\mu  \leq  2^{\leq \kappa }$,  hence  
$\cf (\mu ) < \kappa \Rightarrow  \pp\mu  \leq  \mu ^{<\kappa } \leq  
2^{<\kappa }$,  contradiction, so $(\gamma )$ of (B) follows
from $(\alpha)$.
Note that if we replace $\lambda$ by $\kappa$ (changing the
conclusion  a little; or $\lambda=\kappa$) then by [Sh355, 5.4(2)] if
$2^\lambda$ is regular the conclusion holds for $\chi =
2^\lambda$ too.
}
and let us prove $(\beta )$, so assume  $\chi  \in  (\lambda ,\mu )$, 
without loss of generality, is regular, and we shall prove the
statement in $(\beta )$ of 3.3(B). 
Without loss of generality  $\chi$ is regular and
$\mu ' \in  (\lambda ,\chi ) \und \cf \mu '
\leq  \lambda \Rightarrow  \pp_{\lambda} (\mu ') < 
\chi $;  i.e.\  $\chi $  is $(\lambda ,\lambda ^+,2)$-inaccessible.
[ Why?  If   $\chi $  is not as required, we shall show
how to 
replace  $\chi $  by an appropriate regular  $\chi ' \in [\chi ,\mu)$.]

Let  $\mu ' \in  (\lambda ,\chi )$  be minimal such that  $\pp_{\lambda}
(\mu ') \geq  \chi $,  (so  $\cf \mu ' \leq  \lambda )$  now 
$\pp(\mu ') < \mu $  (by the
choice of  $\mu )$  and $\chi ' =: \pp(\mu ')^+$,  by [Sh355, 2.3] is as 
required$\left. \right] .$

Let  $\theta $  be minimal such that  $2^\theta \geq  \chi $. 
So trivially  $\theta \leq  \kappa \leq  \lambda < \chi $  and 
$\left( 2^{<\kappa }\right) ^\kappa = 2^\kappa $ hence  $\mu  \leq  
2^{<\kappa }$ hence  $\chi  < 2^{<\kappa }$;  as  $\chi $  is regular 
$< 2^{<\kappa }$ but $> \lambda \geq  \kappa $,  clearly  $\theta <
\kappa \leq \lambda $;  also trivially  $2^{<\theta } \leq  \chi 
\leq  2^\theta $ but  
$\chi $  is regular  $> \lambda \geq  \kappa > \theta $  and $[\sigma < 
\theta \Rightarrow  2^\sigma < \chi ]$,  so  $2^{<\theta } < \chi  \leq 
2^\theta $.  Try to apply [Sh410, 4.3] with  $\theta $  here standing
for $\lambda$ there; assumptions (i), (ii) there hold, and if assumption
(iii) fails we get a tree with  $\leq  \theta $  nodes {\sl and} $\geq  
\chi$ $\theta $-branches as required. 
So assume (iii) holds and we get there $\mu '$;  if  $\mu '\leq\lambda$
we have a tree as required; if  $\mu ' \in  (\lambda ,2^{<\theta }] 
\subseteq 
(\lambda ,\chi )$  we get contradiction to  $``\chi $  is 
$(\lambda ,\lambda ^+,2)$-inaccessible" which, without loss of
generality, we have assumed above.

\subdemoinfo{Third Case}{$\kappa$ is singular (hence  $2^{<\kappa }$
is singular,  $\cf (2^{<\kappa }) = \cf \kappa )$}
Let  $\mu  =: 2^{<\kappa }$ and we shall prove (C); easily
(B)$(\gamma)$ holds. 
Now  $^{\kappa >}2$  is a tree with  $2^{<\kappa } = \mu $  nodes and 
$2^\kappa = 2^\lambda$ $\kappa $-branches, so $(\alpha )'$ of (C) holds.
As for $(\beta )$ of (B), if $\kappa$ is strong limit checking the
conclusion is immediate, otherwise it follows from 3.4 part (3) below.

Clearly if  $\cf \kappa > \aleph_{0}$,  also (B) holds. 
\sqed{3.3}

\proclaim{3.4 Claim} 
\item{(1)} Assume $\theta_{n+1} = \Min \left\lbrace \theta \colon
2^\theta > 
2^{\theta_{n}}\right\rbrace $  for  $n < \omega $  and 
$\sum_{n<\omega }\theta_{n} < 2^{\theta_{0}}$ (so  $\theta_{n+1}$ is
regular,  $\theta_{n+1} > \theta_{n})$. 
Then: for infinitely many  $n < \omega$, for some
$\mu_{n} \in [\theta_{n},\theta_{n+1})$ (so  $2^{\mu_{n}} =
2^{\theta_{n}})$  we have: 
{\itemindent=50pt
\item{$(*)_{\mu_{n},\theta_{n}}$} for every regular $\chi\leq 
2^{\theta_{n}}$ there is a tree of cardinality  $\mu_{n}$ with  $\geq  
\chi$ $\theta_{n}$-branches; if $\mu_n>\theta_n$ then
$\cf(\mu_n)=\theta_n$, $\mu_n$ is $(\theta_n, \theta_n^+,
2)$-inaccessible.

}

\item{(2)}  Moreover 

\itemitem{$(\alpha)$} for every $n < \omega $ large enough for some 
$\mu_{n}:$    
$$
\eqalign{
& \theta_{n} \leq  \mu_{n} < \sum_{ m <\omega }\theta_{m} \quad
\hbox{ and } (*)_{\mu_{n},\theta_{n}}\quad\hbox{ and }
\cf (\mu_{n}) = \theta_{n}, \cr
& [\mu_{n} > \theta_{n} \Rightarrow  \mu_{n}\ \hbox{ is }
[(\theta_{n},\theta ^+_{n},2)\hbox{-inaccessible},  
\pp(\mu_{n}) = 2^{\theta_{n}}]. \cr}
$$  

\itemitem{$(\beta)$} Moreover, for infinitely many  $m$  we can demand:
for every  $n < m $,  $\chi  = \cf \chi  \leq  2^{\theta_{n}}$ the
tree  $T^n_{\chi}$ (witnessing $(*)_{\mu_n,\theta_n}$ for
$\chi$) has cardinality  $< \theta_{m +1}$     
(i.e.\  $\mu_{m}  < \theta_{m +1}).$

\item{(3)} If $\kappa$ is singular, $\kappa< 2^{<\kappa} < 2^\kappa$
then for every regular  $\chi  \in  (\kappa ,2^{<\kappa })$,  there
is a tree with  $< \kappa $  nodes $and \geq  \chi $  branches (of
same height). 
Also for some  $\theta ^\ast  \in  (\kappa ,\pp^+(\kappa )) \cap  \Reg$,
for every regular  $\chi  \leq  2^\kappa $ there is a tree  $T$, 
$|T| \leq  \kappa ^{\cf \kappa }$,
with  $\geq  \chi \ \ \theta ^\ast $-branches.
\endproclaim

\demo{Proof} 
Clearly (2) implies (1) and (3) (for (3) second sentence use ultraproduct).
Let
$\theta =: \sum_{n<\omega }\theta_{n}$.  Let  $S_{0} =: \left\lbrace n <
\omega \colon (\ast )_{\theta_{n},\theta_{n}}\right. $\ fails$\left.
\right\rbrace $.  
Let for  $n \in  \omega \backslash S_{0}$,  $\mu_{n} = \theta_{n}$ and
note that $(\alpha )$ of 3.4(2) holds and if $S_0$ is co-infinite, also
$(\beta)$ of 3.4(2) holds.
We can assume that  $S_{0}$ is infinite 
(otherwise the conclusion of 3.4(2) holds). 
By [Sh355, 5.11], fully [Sh410, 4.3] for  $n \in  S_{0}$ there is 
$\mu_{n}$ such that: 

\item{$(\alpha)_{n}$} $\theta_{n} = \cf \mu_{n} < \mu_{n} \leq 
2^{<\theta_{n}}$,

\item{$(\beta)_{n}$}  $\pp_{\Gamma (\theta_{n})}(\mu_{n}) \geq 
2^{\theta_{n}}$ (hence 
equality holds and really  $\pp^+_{\Gamma (\theta_{n})}(\mu_{n}) = 
\left( 2^{\theta_{n}}\right) ^+)$  and 

\item{$(\gamma)_{n}$}  $\theta_{n} < \mu ' < \mu_{n} \und \cf \mu ' \leq 
\theta_{n} \Rightarrow  \pp_{\leq \theta_{n}}(\mu ') < \mu_{n}$ hence  
$\pp^+_{\theta_{n}}(\mu_{n})  =
\pp^+_{\Gamma (\theta_{n})}(\mu_{n}) \shvor = 
\left( 2^{\theta_{n}}\right).$

\noindent
Note that  $2^{<\theta_{n}} = 2^{\theta_{n-1}}$ so  $\mu_{n} \leq  
2^{\theta_{n-1}}$. 
By [Sh355, 5.11] for  $n \in  S_{0}$,  part $(\alpha )$ (of 
3.4(2)) holds except possibly  $\mu_{n} < \theta .$ 

Remember  $\cf (\mu_{n}) = \theta_{n}.$

Let  $n < m$ be in $S_0$ and $\mu_n>\theta_m$, so $\Max\{\cf \mu_{n},\cf
\mu_{m} \} = \Max \{\theta_{n},\theta_{m} \} < \Min \{\mu_{n},
\mu_{m} \}$  so by $(\gamma )_{n}$ (and 
[Sh355, 2.3(2)]) we have  $\mu_{n} \geq  \mu_{m} $. 
Note  $\cf \mu_{n} = \theta_{n}$, 
$\cf \mu_{m}  = \theta_{m} $ (which holds by  $(\alpha )_{n}$,
$(\alpha )_{m} )$  hence  $\mu_{n} > \mu_{m} $. 
As the class of cardinals is well ordered we get
$S_1=: \{n<\omega\colon n\in S_0, \mu_n\ge\theta_{n+1}\}$ is
co-infinite and $S=:\{ n\colon \mu_n\ge\theta\}$ 
is finite (so $(\alpha )$ of 3.4(2)(b) holds).

So for some  $n(\ast ) < \omega $,  $S \subseteq  n(\ast )$  hence for
every  $n \in  [n(\ast ),\omega )$  for some  $m  \in  (n,\omega )$, 
$\mu_{n} < \theta_{m} $. 
Note:  $n \neq  m  \Rightarrow  \mu_{n} \neq  \mu_{m} $ (as their 
cofinalities are distinct) and $[n \notin  S_{0} \Rightarrow  \mu_{n}
\notin  \{\theta_{m} \colon m < \omega \}]$. 
Assume  $n \geq  n(\ast )$,  if  $\mu_{n} > 
\theta_{n+1}$,  let  $m  = m_{n} = \Min \{m \colon \mu_{m +1} > \mu_{n}$
and $m  \geq  n\}$  (it is well defined as  $ \bigvee_{k}
\mu_{n} < \theta_{k}$ and $\theta_{k} < \mu_{k} < \theta =
\bigcup_{{\ell} < {\omega}} \theta_{\ell} )$  and we shall show 
$\mu_{m}  < \theta_{m +1}$;  assume 
not, hence  $m  \in  S_{0}$; so  $\mu_{m +1} \leq  2^{\theta_{m} } = 
\pp_{\Gamma (\theta_{m} )}(\mu_{m} ) \leq  \pp_{\theta_{m +1}}
(\mu_{m} )$  but  $\mu_{m} 
\leq  \mu_{n}$ (by the choice of  $m$)  so
as $\cf(\mu_m)=\theta_m\not= \theta_{m +1}$,  
necessarily  $\mu_{m}  > \theta_{m +1}$ and if  $m  + 1 \notin 
S_{0}$ trivially and
if  $m  + 1 \in  S_{0}$ by one of the demands on  $\mu_{m +1}$ (in
its choice) and [Sh355, 2.3] we have  $\mu_{m +1} \leq  \mu_{m} $; 
but  $\mu_{m}  < \mu_{n}$,  so  
$\mu_{m +1} < \mu_{n}$ contradicting the choice of  $m $. 
So by the last sentence,  $n \geq  n(\ast ) \Rightarrow  \mu_{m_{n}}
< \theta_{m_{n}+1}$. 
By [Sh355, 5.11] we get the desired conclusion (i.e.\ also part
$(\beta )$ of 3.4(2)). 
\sqed{3.4}

\rem{Remark}
It seemed that we cannot get more as we can get an appropriate
product of a forcing notion as in Gitik and  Shelah [GiSh344].



\head 4. Bounds for $\pp_{\Gamma(\aleph_1)}$ for Limits of
Inaccessibles$^*$\endhead

\rem{4.1 Convention} 
For\footnote{}{\hskip-10pt $*$
In previous versions these sections
have been in [Sh410], [Sh420] hence we use  ${\cal Y}$,  etc. (and not
the context of [Sh386]); see 4.2B below.}
any cardinal  $\mu $,  $\mu  > \cf \mu  = \aleph_{1}$ we let 
${\cal Y}_{\mu} $,  $Eq_{\mu} $ be as in [Sh420, 3.1],  $\bar \mu $  is
a strictly 
increasing continuous sequence of singular cardinals of cofinality  
$\aleph_0$ of length $\omega_1$,  $\mu  = \sum_{i<\aleph_{1}}\mu_{i}.$
\endproclaim

So  $\mu $  stands here for  $\mu ^\ast $ in [Sh420, {\S}3, {\S}4,
{\S}5]. 
(Of course,  $\aleph_{1}$ can be replaced by ``regular uncountable".)

\proclaiminfo{4.2 Theorem}{Hypothesis [Sh420, 6.1C]\footnote{*}{I.e.:
if $\fra\subset\Reg$, $|\fra|<\min(\fra)$, $\lambda$
inaccessible then $\lambda>\sup(\lambda\cap\pcf\fra)$.}
}
\item{(1)}  Assume 

\itemitem{(a)} $\mu  > \cf \mu  = \aleph_{1}$,  ${\cal Y} = {\cal Y}_{\mu}$,
$Eq'_{\mu}  \subseteq  Eq_{\mu}$, 

\itemitem{(b)} every $D\in\FIL({\calY})$ is nice (see [Sh420, 3.5]), 
$E = \FIL({\cal Y})$  (or at least there is a nice  ${\cal E}$  (see
[Sh420, 5.2--5], $E=\bigcup{\calE}=\Min\ {\calE}$,  ${\calE}$  is
$\mu$-divisible having weak $\mu$-sums, but we concentrate on the
first case), 

\itemitem{(c)} $\mu <\lambda< \pp^+_{E}(\mu)$, $\lambda$ inaccessible.

\noindent
Then there are  $e \in  Eq_\mu$  and $\langle \lambda_{x}\colon x \in
{\cal Y}/e\rangle $,  a sequence of inaccessibles  $< \mu $  and a
$D \in  \FIL(e,{\cal Y}) \cap  E$  nice to  $\mu $,  $D \in 
\FIL(e,{\cal Y}_{\mu} )$  such that:

\itemitem{$(\alpha)$} $\prod_{x\in \calY_\mu /e}\lambda_{x}/D$  has
true cofinality $\lambda$,

\itemitem{$(\beta)$} $\mu  =\tlim_{D}\langle \lambda_{x}\colon x \in 
{\cal Y}_{\mu} \rangle $.

\item{(2)} We can weaken ``{(\rm b)}" to  $``E \subseteq 
\FIL(Eq,{\cal Y})$
and for  $D \in  E$,  in the game  $wG(\mu ,D,e,{\cal Y})$  the second
player wins choosing filters only from  $E.$

\item{(3)}  Moreover, for given  $e_{0}$,  $D_{0}$, 
$\langle \lambda ^0_{x}\colon x \in  {\cal Y}/e_{0}\rangle $,  if 
$\prod_{x\in {\cal Y}/e_0}\lambda ^0_{x}/D^e_{0}$ is
$\lambda $-directed, then without loss of 
generality  $e_{0} \leq  e$,  $D_{0} \leq  D$  and $\lambda_{x} \leq  
\lambda_{x^{[e_{0}]}}.$
\endproclaim

\rem{4.2A Remark} 
(1) We could have separated the two roles of  $\mu $  (in the 
definition of  ${\cal Y}$,  etc. and in  $\lambda \in  (\mu ,\pp^+_{E}
(\mu )))$  but the result is less useful; except for the unique
possible cardinal appearing later.

(2)  Compare with a conclusion of [Sh386] (see in particular 5.8
there):

\proclaim{Theorem} 
Suppose  $\lambda > 2^{\aleph_{1}}$,  $\lambda $  (weakly) 
inaccessible.
\item{(1)}  If  $\aleph_{1} < \lambda_{i} = \cf \lambda_{i} < \lambda $
for  $i < \omega_{1}$,  $D$  is a normal filter on  $\omega_{1}$, 
$ \prod_{{i} < {\omega_{1}}}\lambda_{i}/D$  is
$\lambda $-directed, then for some  $\lambda '_{i}$,
$\aleph_{1} < \lambda '_{i} = \cf \lambda '_{i} \leq  \lambda_{i}$ and
normal filter  
$D'$  extending  $D$,  $\lambda = \tcf\left( \prod_{{i} <
{\omega_{1}}}\lambda '_{i}/D'\right) $  and $\{i\colon \lambda_{i}$
inaccessible$\} \in  D'.$

\item{(2)}  If  $\aleph_{1} = \cf \mu  < \mu  < \lambda $, 
$\pp_{\Gamma (\aleph_{1})}(\mu ) \geq  \lambda $  then for some 
$\langle \lambda_{i}\colon i < \omega_{1}\rangle $,  
$\aleph_{1} < \lambda_{i} = \cf \lambda_{i} < \mu $,  each 
$\lambda_{i}$ inaccessible and $\lambda \in  \pcf_{\Gamma
(\aleph_{1})}\{\lambda_{i}\colon i < \omega_{1}\}.$
\endproclaim

\demo{Proof of 4.2} 
(1) By the definition of  $\pp^+_{E}(\mu )$  (and assumption $(c)$,
and [Sh355, 2.3 (1) + (3)]) there are  $D \in  E$  and $f \in\  
^{{\cal Y}_{\mu} /e}\mu$ such that:
{\itemindent=32pt
\item{(A)$_{f}$} $\mu  > f(x) = \cf [f(x)] > \mu_{\iota (x)}$, 

}
{
\itemindent=40pt
\item{(B)$_{f,D}$} $\lambda = \tcf\left[ 
\prod_{x\in\calY/e}f(x)/D\right] $.

}
\noindent
Let  $K_{0} =: \left\lbrace (f,D)\colon D \in  E,f
\in\ ^{{\cal Y}_{\mu} /e} \mu \right.$
and conditions (A)$_{f}$ and (B)$_{f,D}$ hold$\left. \right\rbrace $,
so  $K_{0} \neq  \emptyset $. 
Now if  $(f,D) \in  K_{0}$,  for some  $\gamma $\  

{\itemindent=50pt
\item{(C)$_{f,D,\gamma }$} in  $G^\gamma (D,f,e,{\cal Y})$  the second
player wins (see [Sh420, 3.4(2)])

hence  $K_{1} \neq  \emptyset $  where  $K_{1} =: \{(f,D,\gamma ) \in  
K_{0}$\ condition\ (C)$_{f,D,\gamma }$\ holds$\}$.

}

\noindent
Choose  $(f^1,D_{1},\gamma_{\langle \rangle }) \in  K_{1}$ with 
$\gamma_{\langle\rangle}$ minimal. 
By the definition of the game 
$$
\hbox{ for every } A \neq  \emptyset \mod\ D_{1}\hbox{ we have }
(f^1,D_{1} + A,\gamma_{\langle \rangle }) \in  K_{1}.
\leqno(*)
$$
Let  $e_{1} = e(D_{1}).$

\subdemoinfo{Case A}{$\{x\colon f^1(x)$  inaccessible$\} \neq  \emptyset\mod
D_{1}$}
We can get the desired conclusion (by increasing  $D_{1}).$

\subdemoinfo{Case B}{$\{x\colon f^1(x)$  successor cardinal$\} \neq 
\emptyset\mod  D_{1}$}
By $(\ast )$, without loss of generality  $f^1(x) = g(x)^+$,  $g(x)$  a 
cardinal (so $\geq  \mu_{\iota (x)})$  for every  $x \in 
{\cal Y}_{\mu} /e$.  
By [Sh355, 1.3] for every regular $\kappa\in(\mu ,\lambda)$ there is  
$f_{\kappa}  \in\ ^{({\cal Y}/e)}\Ord$  satisfying:

\item{(a)} $f_{\kappa}  < f^1$,  each  $f_{\kappa} (x)$  regular,

\item{(b)} tlim$_{D_{1}}f_{\kappa}  = \mu ,$ 

\item{(c)} $ \prod_{x} f_{\kappa} (x)/D_{1}$ has true cofinality 
$\kappa$.

\noindent
By (a) we get 

\item{(d)} $f_{\kappa}  \leq  g$.

\noindent
By (b) we get, by the normality of  $D_{1}$, that for the
$D_{1}$-majority of  $x 
\in  {\cal Y}/e$,  $f_{\kappa} (x) \geq  \mu_{{\iota }(x)}$;  as  
$f_{\kappa} (x)$  is regular (by (a)) and $\mu_{{\iota }(x)}$ singular 
(see 4.1) we get
\item{(e)} for the $D_{1}$-majority of  $x \in  {\cal Y}/e$,  we have 
$f_{\kappa} (x) > \mu_{{\iota }(x)}$.

Let $\chi$ be large enough, let $N$ be an elementary submodel of  
$(H(\chi),\in ,<^\ast_{\chi})$,\break  $\lambda \in N$, $D_{1}\in N$, 
$N \cap  
\lambda $  is the ordinal  $\|N\|$  (singular for simplicity) and\break
$\{\mu ,\langle f^1,g,f_{\kappa} \colon \kappa \in  \Reg \cap  
(\mu ,\lambda )\rangle \}$  belongs to  $N$. 
Choose $\kappa\in\Reg\cap\lambda\backslash(\sup\lambda\cap N)$,\break 
now in  $ \prod_{x\in\calY/e_{1}}f_{\kappa} (x)/D_{1}$, 
there is a cofinal sequence  
$\langle f_{\kappa ,\zeta }\colon \zeta < \kappa \rangle $;  as 
$\kappa > \sup (\lambda \cap  N)$,  so for some $\zeta(\ast)< \kappa$: 
$$
h \in  N \cap\   ^{{\cal Y}/e_{1}}\Ord \Rightarrow 
\left\lbrace x \in  {\cal Y}/e_{1}\colon f_{\kappa ,\zeta (\ast )}(x)
\leq h(x) < f_{\kappa} (x)\right\rbrace  = \emptyset\mod  D_{1}.
\leqno{\otimes}
$$
[Why?  For any such  $h$  define  $h' \in$ $^{{\cal Y}/e_{1}}\!\!\Ord$
by: $h'(x)$ is $h(x)$ if $h(x)<f_{\kappa}(x)$ and zero otherwise,
so for some  
$\zeta_{h} < \kappa $,  $h' < f_{\kappa ,\zeta_h}$\ mod\ $D_{1}$. 
Let
$\zeta (\ast ) = \sup \left\lbrace \zeta_{h}\colon h \in  N \cap  
^{{\cal Y}/e_{1}}N\right\rbrace $;  it is $< \kappa $  as  $\|N\| <
\kappa $,  and it is as required.]

Let  $f_{\ast}  = f_{\kappa ,\zeta (\ast )}$.  The continuation imitates
[Sh371, {\S}4], [Sh410, {\S}5].

Let 
$$
\eqalign{
K_{2} = \Big\{ (D,\bar B,\langle j_{x}\colon & x  \in 
{\cal Y}/e_{1}\rangle )\colon D_{1}
\subseteq D \in  E, \quad\hbox{ player II wins  }
G^{\gamma_{\langle \rangle }}_{E}(f^1,D),                      \cr
& e_{1} = e(D),\bar B = \langle <B_{x,j}\colon j < j^0_{x} \leq 
\mu_{\iota (x)}>\colon x \in  {\cal Y}/e_{1}\rangle  \in  N,     \cr    
& |B_{x,j_{x}}| \leq  g(x)\ \hbox{ and } j_{x} < j^0_{x} \leq 
\mu_{\iota (x)},                                                \cr
& \{x \in  {\cal Y}/e_{1}\colon f_{\ast} (x)\ \hbox{ is in }
B_{x,j_{x}}\} \in  D \Big\} .                           \cr}
$$  
Clearly $K_2 \not= \emptyset$.
For each  $(D,\bar B,\langle j_{x}\colon x \in  {\cal Y}/e_{1}\rangle ) \in 
K_{2}:$ 

\item{$(*)_{1}$} letting  $h \in\  ^{{\cal Y}/e_{1}}\Ord $,  $h(x) =
|B_{x,j_{x}}|$,  for some  $\bar h = \left\langle
\langle \langle \rangle ,f^1\rangle ,\langle \langle 0\rangle ,h
\rangle\right\rangle$, 
for some  $\gamma_{<0>} < \gamma_{<>}$ and $D$  player II\    
wins in  
$G^{\langle \gamma_{<>},\gamma_{<0>}\rangle }_{E}(D,\bar h,e_{1},
{\cal Y}_{\mu} ).$

\noindent
So choose  $(D,\bar B,\langle j_{x}\colon x \in  {\cal Y}/e_{1}\rangle ,
\gamma_{\langle 0\rangle })$  such that: 
\item{$(*)_{2}$} $(D,\bar B,\langle j_{x}\colon x \in 
{\cal Y}/e_{1}\rangle) \in  K_{2}$,  $(\ast )_{1}$ for 
$\gamma_{\langle 0\rangle }$ holds and (under those restrictions)  
$\gamma_{\langle 0\rangle }$ is minimal.

\noindent
So (as player I can ``move twice"), for every  $A \in  D^+$,  if we
replace  $D$  by  $D + A$, then $(\ast )_{2}$ still holds. 

So without loss of generality (for the first and third members use
normality): 
\item{$(*)_{3}$} one of the following sets belongs to  $D$:    
$$
\eqalign{
A_{0,\zeta } = & \left\lbrace x \in  {\cal Y}/e_{1}\colon \cf
|B_{x,j_{x}}| > \mu_{\iota (x)}\right.\hbox{  and  }
j^0_{x} < \mu_{\zeta} \left.  \right\rbrace \cr
& \quad\quad \quad\hbox{ (for some } 
 \zeta < \omega_{1}\quad\hbox{ such that }
|{\cal Y}/e_{1}| < \mu_{\zeta}), \cr
A_{1} = & \left\lbrace x \in  {\cal Y}/e_{1}\colon \cf |B_{x,j_{x}}| <
\mu_{\iota (x)} \leq  |B_{x,j_{x}}|\right\rbrace,  \cr  
A_{2,\zeta } = & \left\lbrace x \in  {\cal Y}/e_{1}\colon |B_{x,j_{x}}|
\leq  \mu_{\zeta} \right. \hbox{  and  }
j_{x} < \mu_{\zeta} \left. \right\rbrace \quad\hbox{ (for some }
\zeta < \omega_{1}). \cr}
$$

If  $A_{2,\zeta } \in  D$  then  (for  $x \in  {\cal Y}/e_{1})$
$$
B^\ast_{x} =: \bigcup \left\lbrace B_{x,j}\colon x \in 
{\cal Y}/e_{1},j < 
j^0_{x} \hbox{ and }\ |B_{x,j_{x}}| < \mu_{\zeta} 
\hbox{ and } j < \mu_{\zeta} \right\rbrace
$$ 
is a set of $\leq\mu_{\zeta}$
ordinals and
$$
\left\lbrace x\in \calY/e_{1}\colon f_{\ast} (x) \in  B^\ast_{x}\right
\rbrace  \in  D
$$
and $\langle B^\ast_{x}\colon x\in\calY/e_{1}\rangle$  belongs to 
$N$ (as $(D,\bar B,\langle j_{x}\colon x \in\calY/e_{1}\rangle) \in
K_{2}$ and the definition of  $K_{2})$,  contradiction to the choice of
$f_{\ast}$ (see $\otimes$, remember $D_{1}\subseteq D$ by the definition
of $K_{2})$.

If  $A_{1} \in  D$,  we can find $\bar B^1 \in  N$,  $\bar B^1 = 
\langle \langle B^1_{x,j}\colon j < j^1_{x} \leq  \mu_{\iota (x)}\rangle
\colon x \in  {\cal Y}/e_{1}\rangle $,  $|B^1_{x,j}| \leq  g(x)$  and $
\bigwedge_{{j} < {j^1_{x}}}\left[ \cf |B^1_{x,j}| \geq  \mu_{\iota (x)}
\vee  |B^1_{x,j}| = 1\right] $ 
and each  $B_{x,j}$ satisfying $\cf|B_{x,j}|<\mu_{i(x)}$
is a union of  $\cf |B_{x,j}|$  sets\ of the form  
$B^1_{x,j^1}$ of smaller cardinality and so for some 
$j^2_{x} < j^1_{x}$,  
$f_{\ast} (x) \in  B_{x,j_{x}} \Rightarrow  f_{\ast} (x) \in 
B_{x,j^2_{x}} \und 
|B_{x,j^2_{x}}| < |B_{x,j_{x}}|$.  Now playing one move in  
$G^{\langle \gamma_{<>},\gamma_{<0>}\rangle }_{E}(D,\bar h,e,{\cal Y})$ 
we get contradiction to choice of  $\gamma_{\langle 0\rangle }.$

We are left with the case  $A_{0,\zeta } \in  D$,  so without loss of 
generality\break  $ \bigwedge_{{x} , {j}}\cf|B_{x,j}|$
$> \mu_{\iota (x)}$.
Let
$$
\fra = \left\lbrace \cf |B_{x,j}|\colon \cf |B_{x,j}| >
\mu_{\iota (x)},x
\in  {\cal Y}/e_{1}, j < j^0_{x}, j<\mu_\zeta
\hbox{ and } {\iota }(x) > \zeta \right\rbrace,
$$ 
so $\fra$ is a set of regular cardinals, and
(remember  $|\calY/e_{1}| < \mu_{\zeta} $)  we have  $|\fra| < \Min\fra$,
so let $\bar \frb =
\langle \frb_{\theta} [\fra]\colon \theta \in  \pcf \fra\rangle $  be
as in [Sh371, 2.6]. 
So as (by the Definition of  $K_{2})$,  $\langle \langle B_{x,j}\colon
j < j^0_{x}\rangle \colon x \in  
{\cal Y}/e_1\rangle \in N$, clearly $\fra\in N$  hence without loss of 
generality  $\bar \frb \in  N$. 
Let  $\lambda ^\ast  = \sup [\lambda \cap  
\pcf\fra]$,  so by Hypothesis $[420,6.1(C)]$, 
$\lambda ^\ast  < \lambda $,  but 
$\lambda ^\ast  \in  N$,  so  $\lambda ^\ast  + 1 \subseteq N.$ 

By the minimality of the rank we have for every  $\theta \in 
\lambda^\ast  
\cap  \pcf\fra$,\break
$\left\lbrace x \in  y/e_{1}\colon\cf  |B_{x,j_{x}}|
\in  \frb_{\theta} \right\rbrace  = \emptyset $\ mod\ $D$  hence 
$ \prod_{x} \cf |B_{x,j_{x}}|/D$  is $\lambda $-directed, hence
we get contradiction to the minimality of the rank of  $f_{1}.$

(2), (3)  Proof left to the reader. 
\sqed{4.2}

\rem{4.2B Remark} 
\item{(1)}
The proof of 4.3 below shows that in [Sh386] the assumption of
the existence of nice filters is very weak, removing it will cost
a little for at most one place.

\item{(2)}  We could have used the framework of [Sh386] but 
not for 4.3 (or use forcing).

\proclaiminfo{4.3 Claim}{Hypothesis 6.1(C) of [Sh420] even in any 
$K[A]$} 
Assume  $\mu  > \cf \mu  = \aleph_{1}$,  $\mu  > \theta > \aleph_1$,  
$\pp_{\Gamma (\theta ,\aleph_{1})}(\mu ) \geq  \lambda > \mu $, 
$\lambda $  inaccessible. 
Then for some  $e \in  Eq_{\mu} $,  $D \in  \FIL(e,{\cal Y}_{\mu} )$ 
and sequence of inaccessibles  $\langle \lambda_{x}\colon x \in\calY_\mu/e
\rangle $, we have  $\tlim_{D}\lambda_{x} = \mu $  and $\lambda =
\tcf(\prod \lambda_{x}/D)$  
except perhaps for a unique $\lambda$  in  $V$  (not depending on
$\mu )$  
and then  $\pp^+_{\Gamma (\theta ,\aleph_{1})}(\mu ) \leq  \lambda ^+.$
\endproclaim

\demo{Proof}
By the Hyp. (see [Sh513, 6.12]) 
for some $\fra\subseteq \Reg\cap\mu$, $|\fra|<\Min(\fra)$,
$\lambda=\max\pcf(\fra)$, and
$$
(\forall
\lambda'<\lambda)(\exists\frb)[\frb\subseteq\fra\ \&\
|\frb|<\theta>] \&\lambda>\sup\pcf_{\aleph_1\scriptstyle{\rm -complete}}
(\frb)>\lambda'],
$$
$J=J_{<\lambda}[\fra]$.
First assume ``in  $K[A]$  there is a Ramsey cardinal 
$> \lambda ^{\theta}$ when  $A \subseteq 
\lambda ^{\theta}$''.  Choose  $A 
\subseteq 
\lambda ^{\theta}$ such that  $^{\theta}\lambda  
\subseteq 
L[A]$  and for every  $\alpha < \lambda^{\theta}$,  there is a one to
one 
function  $f_{\alpha} $ from  $|\alpha|$ (i.e.\  $|\alpha |^V)$  onto  
$\alpha $,  $f_{\alpha}  \in  L[A]$,  so $\Card^{L[A]} \cap  
\left( \lambda^\theta+1\right)  = \Card^V$,\  and apply 4.2 to the
universe  $K[A]$  (its assumption holds by [Sh420, 5.6]). 

Second assume $(\ast )_{\lambda} $ ``in  $K[A]$  there is a Ramsey
cardinal  $> \lambda $  when  $A 
\subseteq 
\lambda ^+"$  and assume our desired conclusion fails. 
Let  $S \subseteq  
\lambda $  be stationary  $[\delta \in  S \Rightarrow  \cf \delta = 
\theta ^+]$,  $\langle a_{\alpha} \colon \alpha < \lambda \rangle $, 
exemplify  $S 
\in  I[\lambda ]$ (exist by [Sh420, \S1]).
We can find $\fra$, $J$ as described above.
Let  $\langle f_{\alpha} \colon \alpha < 
\lambda \rangle $  exemplify  $\lambda = \tcf(\prod \fra/J)$, 
now by [Sh355, 1.3] 
without loss of generality  $\lambda = \max  \pcf\fra$.  Let  $A_{0} 
\subseteq 
\lambda$ be such that $\fra$, $\langle f_{\alpha} \colon \alpha <
\lambda \rangle $,  
$\langle \frb_{\sigma} [\fra]\colon \sigma \in  \pcf \fra\rangle $  are in 
$L[A_{0}]$.  Hence in
$L[A_{0}]$  for suitable  $J$,  $\langle f_{\alpha} /J\colon \alpha <
\lambda \rangle $ 
is increasing, and without loss of generality for some  
$\langle\langle\frc^\delta_{\alpha}\colon\alpha\in  a_{\delta} \rangle
\colon \delta \in S\rangle  \in 
L[A_{0}]$,  we have: for  $\delta \in  S$,  $\cf \delta = |\fra|^+$,  
$a_{\delta} $ a club of  $\delta $  and 
$\langle f_{\alpha} \upharpoonright (\fra\backslash
\frc^\delta_{\alpha} )\colon \alpha \in  
a_{\delta} \rangle$ is $<$-increasing (see [Sh345b, 2.5]
(``good point'')) and 
$\frc^\delta_{\alpha}  \in  J$  and $S$  is stationary in  $V$,  so the
assumption
of 4.3 holds in  $V^1$ whenever  $L[A_{0}] 
\subseteq 
V^1 
\subseteq 
V$;  hence for  $A 
\subseteq 
\lambda ^+$,  in  $K[A_{0},A]$  the conclusion of 4.2 holds as
we are assuming $(*)_\lambda$. 

Note: if  $A 
\subseteq 
\lambda $,  in  $K[A]$,  $\lambda ^{<\lambda } = \lambda $  hence if 
$\alpha < \lambda ^+$,  $A \subseteq 
\alpha $  then  $K[A] \models  ``\lambda ^{<\lambda } <
(\lambda ^+)^V".$ 

Choose by induction on  $\alpha < \lambda ^+$ a set  $A_{\alpha}  
\subseteq 
[\lambda \alpha ,\lambda (\alpha + 1))$  such that:  $A_{0}$ is as
above and for $\alpha > 0$:  if  $\langle \lambda_{x}\colon x
\in\calY/e\rangle $,  $J$
exemplify the conclusion of 4.2 in  $K\left[ \bigcup_{{\beta} <
{\alpha}} A_{\beta} \right] $,  and $\langle f_i \colon i <
\lambda \rangle $
exemplify the  $\lambda = \tcf\left(\prod_{{x} \in  {{\cal Y}}/e}
\lambda_{x}/J\right) $,
without loss of generality  $J$  canonical  (all in 
$K\left[\bigcup_{{\beta} < {\alpha}} A_{\beta} \right]$, 
canonical means: 
the normal ideal generated by  $\{x\colon \lambda_{x} \in  
\frb_{<\lambda }[\{\lambda_{y}\colon y \in  {\cal Y}/e\}]\})$,  then in 
$K\left[ \bigcup_{{\beta} \leq  {\alpha}} A_{\beta} \right] $ 
we can find $f$,  $ \bigwedge_{{\alpha} < {\lambda}} f <_{J}
\langle \lambda_{x}\colon x \in  {\cal Y}/e\rangle $,  $ 
\bigwedge_{\alpha } f \not<_{J} f_{\alpha} $ (as they cannot exemplify
the conclusion of 4.5 in  $V$  --- otherwise we have finished). 

Let  $A = \bigcup_{{\alpha} < {\lambda ^+}} A_{\alpha} .$ 

Now in $K[A]$ there are $e$, $\langle \lambda_{x}\colon \lambda \in  
{\cal Y}/e\rangle$, $\langle f_i\colon i <\lambda\rangle$ 
(and $J)$\break exemplifying the conclusion of 4.2 (by $(*)$ and
[Sh513, 6.12(3)]).
By 4.5 below, for some  $\delta < 
\lambda^+$, $e$, $\langle\lambda_{x}\colon x \in  {\cal Y}/e\rangle $,  
$\langle \frb_{\sigma} [\{\lambda_{x}\colon x \in  {\cal Y}/e\}]\colon
\sigma \in  
\pcf \{\lambda_{x}\colon x \in  {\cal Y}/e\}\rangle $,  $f_{\alpha} (\alpha
< \lambda )$ all belongs to  $K\left[ \bigcup_{{\gamma} <
{\delta}} A_{\gamma} \right] $,  and in  $K\left[
\bigcup_{{\gamma} \leq  {\delta}} A_{\gamma} \right] $  we get a
contradiction.

If $(\ast )_{\lambda} $ holds for every  $\lambda $  we are done. 
If not, let  $\lambda_{0}$ be minimal such that $(\ast )_{\lambda_{0}}$
fails; so if $\lambda < \lambda_{0}$ the conclusion holds, and if 
$\lambda > \lambda_{0}$ then let  $A \subseteq 
\lambda ^+_{0}$ be such that in  $K[A]$  there is no Ramsey, hence
([DoJ]) for  $\mu  \geq  \lambda ^+_{0}$ in  $V$,  $\cov(\mu ,\theta,
\theta ,2) \leq  \mu $,  so the assumptions of 4.3 fail. 
Similarly $\mu>\theta$, $\cf(\mu)=\aleph_1$,
$\pp_{\Gamma(\theta,\aleph_1)}(\mu)>\lambda_0^+$ bring a contradiction.
\sqed{4.3}

\reminfo{4.4 Conclusion}{Hypothesis [Sh420, 6.1(C)] in any  $K[A]$}
(1) Assume  $\mu  > \cf \mu  = \aleph_{1}$,  $\mu_{0} < \mu $, \
$\sigma \geq  |\{\lambda \colon \mu_{0} < 
\lambda < \mu $,\ \ $\lambda $\ \ inaccessible$\}| < \mu $.  Then  
$$
\sigma^{+4} > |\big\{ \lambda \colon \mu  < \lambda
< \shvor \pp_{\Gamma (\sigma ,\aleph_{1})}(\mu )
\hbox{ and } \lambda \ \hbox{ is inaccessible} \big\} |.
$$

\noindent
(2) The parallel of [Sh400, 4.3].

\demo{Proof} 
See [Sh410, 3.5] and use 4.2(3).
\qed

By [DoJe]

\proclaim{4.5 Theorem} 
If  $\lambda $  is regular  $(> \aleph_{1})$ $A \subseteq 
\lambda $,  $Z \in  K[A]$  a bounded subset of  $\lambda $  then for
some  $\alpha < \lambda $,\  $Z \in \bigcup_{{\alpha} <
{\lambda}} K[A \cap  \alpha ].$
\endproclaim

We shall return to this elsewhere.



\head 5. Densities of Box Products \endhead

\rem{5.1 Definition} 
$d_{<\kappa }(\lambda ,\theta )$  is the density of the 
topological space  $^\lambda \theta $  where the topology is generated
by the following family of clopen sets: 
$$
\{[f]\colon f\in\ {}^a\theta\hbox{ for some } a\subseteq\lambda, |a|
< \kappa \}
$$
where
$$
[f] = \{g \in\ {}^\lambda\theta \colon g \subseteq  f\}\hbox{.}
$$
So 
$$
\eqalign{
d_{<\kappa } & (\lambda ,\theta ) = \cr
& \Min \left\lbrace |F|\colon F \subseteq\ 
{}^\lambda\theta \hbox{ and if }a \in  {\cal S}_{<\kappa }(\lambda )
\hbox{ and }g \in\  {}^a\theta \hbox{ then }
(\exists f \in  F)g \subseteq  f\right\rbrace. \cr}
$$

\noindent
If  $\theta = 2$  we may omit it, if  $\kappa = \aleph_{0}$ we may
omit it (i.e.\ $d(\lambda ,\theta ) = d_{<\aleph_{0}}(\lambda,
\theta)$). 
Always we assume  $\lambda \geq  \aleph_{0}$,  $\kappa \ge\aleph_0,
\theta> 1$ and $\lambda^+\ge\kappa$. 
We write $d_{\kappa}(\lambda,\theta)$ for $d_{<\kappa^+}(\lambda,\theta).$

\rem{5.1A Discussion} 
Note: for  $\kappa = \aleph_{0}$ this is the Tichonov 
product, for higher  $\kappa $  those are called box products and $d$ 
has obvious monotonicity properties. 

$d\left( 2^{\aleph_{0}}\right)  = \aleph_{0}$ by the classical 
Hewitt--Marczewski--Pondiczery theorem [H], [Ma], [P].  This has been
generalized by Engelking--Karlowicz [EK] and by Com\-fort--Negrepontis
[CN1], [CN2] to show, for example, that $d_{<\kappa}(2^\alpha
,\alpha) = \alpha $  if and 
only if  $\alpha = \alpha ^{<\kappa }$ ([CN1] (Theorem 3.1)).  
Cater--Erd\H os--Galvin [CEG] show that every non-degenerate space $X$  
satisfies  $\cf (d_{<\kappa}(\lambda ,X))\geq\cf(\kappa )$  when 
$\kappa \leq  \lambda ^+$,  and they note (in our notation) that 
$``d_{<\kappa} (\lambda )$  is
usually (if not always) equal to the well-known upper bound $(\log  
\lambda )^{<\kappa }"$. 
It is known ($\cf$. [CEG], [CR]) that  $\SCH \Rightarrow  
d_{<\aleph_{1}}(\lambda ) = (\log  \lambda )^{\aleph_{0}}$,  but it is
not known whether  $d_{<\aleph_{1}}(\lambda ) = (\log 
\lambda)^{\aleph_{0}}$ is a theorem of ZFC.

The point in those theorems is the upper bound, as, of course,  
$d_{<\kappa }(\mu,\theta) >\chi$ if $\mu  > 2^\chi  \und \theta > 2$
[why? because if $F = \{f_{i}\colon i< \chi \}$ exemplify  $d_{<\kappa}
(\mu ,\theta ) \leq  \chi $,  the number of possible sequences  
$\langle\Min\{1,f_{i}(\zeta)\}\colon i< \chi\rangle$ (where  $\zeta <
\mu )$ is $\leq  2^\chi$, so for some  $\zeta \neq  \xi $  they are
equal and we get contradiction by  $g$,  $g(\zeta ) = 0$,  $g(\xi ) =1$,
$\Dom g = \{\zeta ,\xi \}].$

Also trivial is: for  $\kappa $  limit,  $d_{<\kappa} (\lambda
,\theta) = \kappa + \sup_{\sigma <\kappa}d_{<\sigma}(\lambda,\theta )$,
so we only use  $\kappa $  regular;  
$d_{<\kappa }(\lambda ,\theta ) \geq  \sigma ^\theta $ for  $\sigma < 
\kappa .$

Also if  $\cf(\lambda)< \kappa $,  $\lambda $  strong limit then  
$d_{<\kappa }(\lambda ) > \lambda $. 
The general case (say  $2^{<\mu } < 
\lambda < 2^\mu $,  $\cf \mu  \leq  \theta $)  is similar; we ignore
it in order to make the discussion simpler.

So the main problem is:

\noindent
\proclaim{5.2 Problem} 
Assume  $\lambda $  is strong limit singular, $\lambda>\kappa > 
\cf(\lambda)$,  what is  $d_{<\kappa }(\lambda )?$  Is it always 
$2^\lambda ?$  Is it always  $> \lambda ^+$ when  $2^\lambda >
\lambda ^+?$ 
\endproclaim

In [Sh93] this question was raised (later and independently) for
model theoretic reasons. 
I thank Comfort for 
asking me about it in the Fall of '90.

\proclaim{5.3 Lemma} 
Suppose  $\lambda $  is singular strong limit,  $\cf(\lambda)= 
\cf (\delta ^\ast ) \leq  \delta ^\ast  < \cf(\kappa)\leq  \kappa <
\lambda $,  
$2 \leq  \theta < \lambda $,  $\lambda \leq  \chi  < 2^\lambda $ and 
$\langle \lambda_{\alpha} ,\mu_{\alpha} ,\chi_{\alpha},
\chi ^\ast_{\alpha} :%
\alpha < \delta ^\ast \rangle $  is such that: 
\item{} $\chi_{\alpha}=\theta^{\mu_{\alpha}}, \chi ^\ast_{\alpha}  = 
\cov(\chi_{\alpha} ,\lambda_{\alpha} ,\lambda_{\alpha} ,2)$,

\item{} $\alpha  < \beta \Rightarrow  \mu_{\alpha}  < \mu_{\beta} $,

\item{} $\lambda= \bigcup_{{\alpha} < {\delta ^\ast }}\mu_{\alpha}  =
 \tlim_{\alpha <\delta } \lambda_{\alpha}$, $\theta<\mu_\alpha$,

\item{} $d_{<\kappa }(\mu_{\alpha} ,\theta ) \geq  \lambda_{\alpha}$
(this holds e.g.\ if $(\forall \lambda'<\lambda_\alpha)
[2^{\lambda'}<\mu_\alpha]$),

\item{} $A_{\alpha}   = [\mu_{\alpha} ,\mu_{\alpha}  + \mu_{\alpha} ]$,

\item{} $G_{\alpha} = \{g\colon g \hbox{ a partial function from some }
a \in  {\cal S}_{<\kappa }(A_{\alpha} ) \hbox{ to } \theta \}$,

\item{} $\phantom{G_{\alpha}=}$ for $g \in  G_{\alpha}$,

\item{} $[g] = \{f \in  X_{\alpha} \colon g \subseteq  f\}
\hbox{ where } X_{\alpha}  =:\ ^{(A_{\alpha} )}\!\theta,\hbox{  so  }
|X_{\alpha} | = \chi_{\alpha}$,

\item{} $h_{\alpha}$  is a function from ${\cal S}_{<\lambda_{\alpha} }
\!\left( ^{(A_{\alpha} )}\theta \right)
\hbox{  to }  G_{\alpha} \hbox{ such that } h_{\alpha} (a)$
``exemplifies" 
\itemitem{} that $a$ is not dense in 
$^{(A_{\alpha} )}\theta$, i.e. 
$[f \in  a \und g = h_{\alpha} (a) \Rightarrow  g {\miss} f]$.

\noindent
Then  (F)$\Rightarrow$(E)$\Rightarrow$(D)$\Leftrightarrow$(C)
$\Rightarrow$(B)$\Leftrightarrow$(A); and (E)$^\sigma $
decrease with  $\sigma $  and (E)$^\sigma\Rightarrow$(G) when
$\chi^*_\alpha=\chi_\alpha$; and
if every $\lambda_\alpha$ is regular
(G)$\Rightarrow$(F) and if in addition
$\bigwedge_{\alpha<\delta^*}\chi^*_\alpha=\chi_\alpha$ then
(G)$\Leftrightarrow$(F)$\Leftrightarrow$(E), and if
$\{\alpha<\delta^*\colon\sigma\le\lambda_\alpha\}\not=\emptyset\mod
J$ and $\sigma<\lambda$ then (E)$\Leftrightarrow$(E)$^\sigma$ (fixing
$J$), where

\item{(A)} $d_{<\kappa }(\lambda ,\theta ) > \chi $; 

\item{(B)} if  $x_{\zeta}  \in \prod_{{\alpha} < {\delta ^\ast
}} X_{\alpha} $ for  $\zeta < \chi $  then there is  $\bar g \in  
 \prod_{{\alpha} < {\delta ^\ast }}G_{\alpha} $ such that:\   
for every  $\zeta < \chi,   \{\alpha < \delta ^\ast \colon x_{\zeta}
(\alpha) \notin [g_{\zeta} ]\} \neq  \emptyset $;

\item{(C)} if  $x_{\zeta}  \in \prod_{{\alpha} < {\delta ^\ast
}} X_{\alpha} $ for  $\zeta < \chi $  then for some  $w_{\alpha}  
\in \calS_{<\lambda_{\alpha}}(X_{\alpha})$ $(\alpha< \delta^\ast)$
for every  $\zeta < \chi,   \{\alpha < \delta ^\ast \colon x_{\zeta}
(\alpha) \in  w_{\alpha} \} \neq  \emptyset $; 

\item{(D)} for every  $x_{\zeta}  \in \prod_{{\alpha} <
{\delta^\ast }} \chi_{\alpha} $ for  $\zeta < \chi $  there is 
$\bar w \in   \prod_{{\alpha} < {\delta ^\ast }} {\cal S}_{
<\lambda_{\alpha} }(\chi_{\alpha} )$  such that:  
for each  $\zeta < \chi $,  $ \bigvee_{{\alpha} < {\delta ^\ast
}}x_{\zeta} (\alpha ) \in  w_{\alpha} $;

\item{(E)$^\sigma$} for some ideal  $J$  on  $\delta ^\ast $ extending
$J^{bd}_{\delta ^\ast }$ for every  $x_{\zeta}  \in
\prod_{{\alpha} < {\delta ^\ast }}\chi_{\alpha} $\ (for
$\zeta < \chi )$  there are  
$\epsilon (\ast ) < \sigma $  and $\bar w^\epsilon  \in
\prod_{{\alpha} < {\delta ^\ast }} {\cal S}_{<\lambda_{\alpha} }
(\chi_{\alpha} )$  for  
$\epsilon  < \epsilon (\ast )$  such that for each  $\zeta $  we have  
$\bigvee_{\epsilon } \{\alpha < \delta ^\ast \colon x_\zeta(\alpha )
\notin  w^\epsilon_{\alpha} \} = \emptyset $\ mod\ $J.$ 

If  $\sigma = 2$  we may omit it;

\item{(F)} for some non-trivial ideal  $J$  on  $\delta ^\ast $
extending $J^{bd}_{\delta^*}$ we have
$$ \prod_{{\alpha} < {\delta ^\ast }}\left(
{\cal S}_{<\lambda_{\alpha} }(\chi_{\alpha} ),\subseteq %
\right)/J\hbox{ is }\chi ^+\hbox{-directed};
$$

\item{(G)} for some non-trivial ideal  $J$  on  $\delta ^\ast$
extending $J^{bd}_{\delta^*}$, for any
$\langle {\cal P}_{\alpha} \colon \alpha < \delta ^\ast \rangle $, 
${\cal P}_{\alpha} $
a $\lambda_{\alpha} $-directed partial order of cardinality  $\leq  
\chi ^\ast_{\alpha} $,  we have: $ \prod_{{\alpha} <
{\delta ^\ast }}{\cal P}_{\alpha} /J$  is $\chi ^+$-directed.

\rem{5.3A Remark} 
\item{(1)} Note that the desired conclusion is 5.2(A).

\item{(2)} The interesting case of 5.3 is when  $\{\mu_{\alpha}\colon
\alpha < \delta ^\ast \}$  does not contain a club\break of  $\lambda$.

\item{(3)}  Note that with notational changes we can arrange 
$``\lambda $  is the 
disjoint union of  $A_{\alpha} (\alpha < \delta ^\ast )$,  hence  
$\lambda_\theta = \prod_{{\alpha} < {\delta ^\ast }}
X_{\alpha}".$

\demo{Proof} 
Check.  Clearly  $(E)^\sigma $ decreases with  $\sigma$, i.e.\ if
$\sigma_1<\sigma_2$ then $(E)^{\sigma_1}\Rightarrow
(E)^{\sigma_2}$.

\subdemo{(E)$\Rightarrow$(D)} 
Just for  $J$  varying on non-trivial ideals, we have 
monotonicity in  $J$;  and for  $J = \{\emptyset \}$  we get (D).

\subdemo{(D)$\Leftrightarrow$(C)}
(C) is a translation of (D).

\subdemo{(C)$\Rightarrow$(B)} 
If $x_{\zeta}\in\prod_{{\alpha}<{\delta ^\ast }}X_{\alpha}$
for $\zeta<\chi$,  let  $\langle w_{\alpha}\colon\alpha< \delta^\ast
\rangle $  be as in (C); for
each  $\alpha $  we know that  $w_{\alpha} $ is not a dense subset of 
$X_{\alpha} $  (as\  $d_{<\kappa }(\mu_{\alpha} ,\theta ) \geq 
\lambda_{\alpha}>|w_\alpha|)$  so there is
$g_{\alpha}  \in  G_{\alpha} $ for which $[g_\alpha]\cap
w_\alpha=\emptyset$, so  $\bar g =: \langle g_{\alpha}
\colon \alpha < 
\delta ^\ast \rangle $  is as required in (B).

\subdemo{(B)$\Leftrightarrow$(A)}
They say the same (see 5.3A(3)).

\subdemo{(F)$\Rightarrow$(E)} 
Note that (E) just says that\ in  $ \prod_{{\alpha} <
{\delta ^\ast }}\left( {\cal S}_{<\lambda_{\alpha} }(\chi_{\alpha} ),
\subseteq \right) $,  any subset of  $\left\lbrace f\colon  f \in
\prod_{{\alpha} < {\delta ^\ast }} 
{\cal S}_{<\lambda_{\alpha} }(\chi_{\alpha} )\right. $, such that each
$f (\alpha )$\ \ is\ a\ singleton$\left. \right\rbrace $  has
a $\leq_{J}$-upper bounded. 
In this form it is clearly a specific case of (F).

\subdemo{(E)$^\sigma\Rightarrow$(G) when
$\chi_\alpha=\chi_\alpha^*$}
where $\{\alpha<\delta^*\colon\sigma\le\lambda_\alpha\}\not=
\emptyset\mod J$: Easy too.               

\smallskip

\noindent
Next assume every  $\lambda_{\alpha} $ is regular,  $J$  an ideal on  
$\delta ^\ast .$

\subdemo{(G)$\Rightarrow$(F)}
(F) is a particular case of (G), because  
$({\cal S}_{<\lambda_{\alpha} }(\chi_{\alpha} ) \subseteq )$  is  
$\lambda_{\alpha}$-directed as  $\lambda_{\alpha} $ is regular and 
${\cal S}_{<\lambda_{\alpha} }(\chi_{\alpha} )$  can be replaced by
any cofinal subset and there is one of cardinality  $\chi^\ast_\alpha$
by its definition.

\noindent
The rest should be clear.
\sqed{5.3}

\proclaim{5.4 Claim} 
Assume  $\lambda $  is strong limit, $\theta<\lambda_0$,  
$\langle \lambda_{\alpha} \colon \alpha < \delta ^\ast \rangle $,  
$\langle \chi ^\ast_{\alpha} \colon \alpha < \delta ^\ast \rangle $  are
(strictly) increasing with limit  $\lambda $,  $\delta ^\ast <
\kappa\le\cf(\lambda)< \lambda$,
$\lambda < \chi  < 2^\lambda $ and $\lambda_\alpha\le\chi_\alpha^*$,
$\lambda_\alpha$ regular for each $\alpha<\delta^*$. 
Then (G) of {\rm 5.3} holds (hence  
$d_{<\kappa}(\lambda ,\theta) > \chi)$ in any of the following cases:

\item{(a)} for some  $\mu_{\alpha} $ strong limit,  $\cf (\mu_{\alpha})<
\kappa $,\  $2^{\mu_{\alpha} } = \mu ^+_{\alpha}$, $\lambda_{\alpha} =
\mu ^+_{\alpha} $,  
$\chi ^\ast_{\alpha}  = \mu ^+_{\alpha} $ and
$\prod_{\alpha<\delta^*}\mu^+_\alpha/J$ is $\chi^+$-directed,

\item{(b)} $k < \omega $  and for every  $\alpha $, 
$\chi^\ast_{\alpha}  \leq  \lambda ^{+k}_{\alpha} $ and for some ideal
$J$  on  $\delta ^\ast $,  for  
$\ell  \leq  k$,  $\prod \lambda ^{+\ell }_{\alpha} /J$  is
$\chi ^+$-directed, and
$d_{<\kappa}(\chi^*_\alpha,\theta)\ge\lambda_\alpha$,

\item{(c)}  for some  $\gamma < \cf(\lambda)$  for every 
$\alpha<\delta^*$, $\chi^\ast_{\alpha} \leq\lambda^{+\gamma }_{\alpha}$
and for some ideal  $J$ on 
$\delta ^\ast $ for every  $\zeta < \gamma $,  $
\prod_{{\alpha} < {\delta ^\ast }}$,  $\lambda ^{+(\zeta +1)}_{\alpha}
/J$  is $\chi ^+$-directed, and
$d_{<\kappa}(\chi^*_\alpha, \theta)\ge\lambda_\alpha$,

\item{(d)} for some ideal  $J$  on  $\delta ^\ast $ extending  
$J^{bd}_{\delta ^\ast }$ for every regular  $\lambda '_{\alpha}  \in  
[\lambda_{\alpha} ,\chi ^\ast_{\alpha} ]$ satisfying  $\tlim_J 
 (\cf \lambda'_{\alpha} ) = \lambda $,
we have  $ \prod_{{\alpha} < {\delta ^\ast }}\lambda '_{\alpha}
/J$  is $\chi ^+$-directed and
$d_{<\kappa}(\chi^*_\alpha,\theta)\ge\lambda_\alpha$.
\endproclaim

\demo{Proof} 
Clearly (a)$\Rightarrow$(b)$\Rightarrow$(c)$\Rightarrow$(d).

Now the statements follow from the following observations
5.4A--5.7.

\rem{5.4A Observation} 
Assume that for  $\alpha < \delta $,  ${\cal P}_{\alpha} $ is
a (non-empty) $\lambda_{\alpha} $-directed partial order of cardinality 
$\chi_{\alpha} $,  $|\delta|^+<\lambda_{\alpha} =\cf(\lambda_{\alpha})
\leq  \chi_{\alpha} $,  $J$  an ideal on  $\delta $,  $\theta ^\ast  = 
\Min\{\theta\colon$ for some  A and $\bar f\colon \bar f = \langle
f_{i}\colon i < \theta \rangle$,  $f_{i} \in \prod_{{\alpha} <
{\delta}} {\cal P}_{\alpha}$ is $<_{J+A}$-increasing,  $A \subseteq 
\delta $,  
$\delta \backslash A \notin  J$  but for no  $g \in
\prod_{{\alpha} < {\delta}} {\cal P}_{\alpha} $,  $
\bigwedge_{{i} < {\theta}} \{\alpha \colon\calP_\alpha\models
f_{i}(\alpha ) \leq  g(\alpha )\} \neq  \emptyset $\ mod\ $(J+A)\}$. 
{\sl Then} $ \prod_{{\alpha} < {\delta}} {\cal P}_{\alpha} /J$  is
$\theta ^\ast $-directed.
\endproclaim

\demo{Proof} 
Without loss of generality no  ${\cal P}_{\alpha} $ has a maximal 
element.  If the conclusion of 5.4A fails, let  $F$  be a subset of
$ \prod_{{\alpha} < {\delta}} {\cal P}_{\alpha} $ with
no $<_{J}$-upper bound, of minimal cardinality.  
Let  $\theta = |F|$,  so let  $F = \{f_{i}\colon i < \theta \}$;  by the
choice of  $F$  without loss of generality  $\alpha < \beta
\Rightarrow  f_{\alpha}  <_{J} 
f_{\beta} $ hence  $\theta $  is necessarily regular. 
If  $\{\alpha < \delta \colon \lambda_{\alpha} \leq  \theta \} \in  J$  we
can find an upper bound:  
$g(\alpha )$  is a ${\cal P}_{\alpha} $-upper bound of 
$\{f_{i}(\alpha )\colon i < \theta \}$  when  $\lambda_{\alpha}  > \theta $,
and arbitrarily otherwise. 
So without loss of generality $\bigwedge_{\alpha}\lambda_{\alpha}\leq 
\theta$. 
Now, remember $|\delta|^+< \lambda_{\alpha} $,  and so 
 $|\delta |^+ < \theta $. 
By [Sh420, \S1] we can find ${\bar C}=\langle C_i\colon i<\theta\rangle$,
$C_i\subset i$, $j\in C_i\Rightarrow C_j = j\cap C_i$, 
$\otp (C_i)\le |\delta|^+$ and $S=:\{ i<\lambda\colon \cf (i)=
|\delta|^+, \delta=\sup(C_i)\}$ stationary: so wlog $j\in C_i
\Rightarrow \bigwedge_{\alpha<\delta} \calP_\alpha\models
f_j(\alpha) <f_i(\alpha)$.
Now we repeat the proof from [Sh282, 14]; better see [Sh345a, 2.6]
or here 6.1.\footnote*{In the main case here, $\bigwedge_\alpha
2^{|\delta^*|} < \lambda_\alpha$ and then trying all the possible $A$'s,
using their $g$'s, the proof is very simple.}
\sqed{5.4A}

\rem{5.5 Observation} 
In 5.4A, if  $A$, $\bar f$  exemplify  $\theta ^\ast  =
\theta$ then
$$
\theta^\ast\ge\min\{\pre_{J+A}(\bar\chi,\bar\lambda)\colon
A\subseteq\delta\hbox{ and }\delta\setminus A\not\in J\}
$$
where
\endproclaim

\rem{5.6 Definition} 
For ideal $I$ on  $\delta $ and $\bar\chi =\langle\chi_\alpha
\colon \alpha<\delta\rangle$, $\bar\lambda=\langle\lambda_\alpha
\colon\alpha<\delta\rangle$, $\lambda_\alpha=\cf(\lambda_\alpha)\le
\chi_\alpha$ we let
$\pre_{I}(\bar \chi ,\bar \lambda ) =: \Min\left\lbrace |{\cal P}|
\colon  {\cal P}\right. $\ \ is a family of sequences of the form  
$\langle B_{\alpha} \colon \alpha < \delta \rangle ,$      
$B_{\alpha}  \subseteq  \chi_{\alpha} $,  $|B_{\alpha} | <
\lambda_{\alpha} $ such 
that for every  $g \in \prod_{{\alpha} < {\delta}} \chi_{\alpha}$
for some  $\bar B \in  {\cal P}$,  $\{\alpha < 
\delta \colon g(\alpha ) \in  B_{\alpha} \} \neq  
\emptyset $\ mod\ $I\left. \right\rbrace .$

\demo{Proof} 
Check.

\rem{5.6A Remark} 
We use other parts of 5.3.

\rem{5.7 Observation} 
Let $I$ be an ideal on  $\delta ^\ast $,  
$\chi_{\alpha}  \geq  \lambda_{\alpha}  > \delta ^\ast $.

\item{(1)} Define $\calJ[I] = \{I + A\colon A \subseteq  \delta ,\delta
\backslash A \notin  I\}$. 

\item{(2)}  If  $I_{1} \subseteq  I_{2}$,  $\lambda ^1_{\alpha}  \geq 
\lambda ^2_{\alpha} $,  $\chi_\alpha^1 \leq  \chi_\alpha^2$ for
$\alpha<\delta$ then $\pre_{I_{1}}(\bar \chi ^1,\bar \lambda ^1) 
\leq  \pre_{I_2}(\bar \chi ^2,\bar \lambda ^2).$

\item{(3)}  If  $\delta ^\ast $ is the disjoint union of  $A_{1}$, 
$A_{2}$,  $A_{\ell}  \notin  I$  and $I_{\ell}  =: I + A_{\ell} $ then  
$\pre_{I}(\bar \chi ,\bar \lambda ) = \Min \left\lbrace \pre_{I_{1}}
(\bar \chi ,\bar \lambda ),\pre_{I_{2}}(\bar \chi ,\bar %
\lambda )\right\rbrace .$

\item{(4)}  $\pre_{I}(\bar \chi ^+,\bar \lambda ) \leq  \pre_{I}(\bar\chi,
\bar \lambda ) + \sup \{\tcf(\prod\chi^+_{\alpha} /I + A)\colon A
\subseteq  \delta ,\delta \backslash A \notin
I\}.$\footnote{**}{Of course, $\bar\chi^+=\langle \chi^+_\alpha\colon
\alpha<\delta\rangle$.}

\item{} Moreover $\pre_I(\bar\chi^+,\bar\lambda)\le\Min\{\pre_{I+A}
(\bar\chi,\bar\lambda)+\tcf(\prod_{\alpha<\delta} \chi_\alpha^+
/ (I+A))\colon A\subseteq \delta, \delta\hefresh A\not\in I$
(and the tcf is well defined)$\}$.

\item{(5)}  If each  $\chi_{\alpha} $ is a limit cardinal, 
$\cf \chi_{\alpha}  > \delta ^\ast$, then $\sup_{J\in\calJ[I]}
\pre_J(\bar\chi, \bar\lambda)=$\break $\sup_{\bar\chi'<\bar\chi}
\shvor\sup_{J\in\calJ[I]}\pre_J(\bar\chi',
\bar\lambda)+\sup_{J\in\calJ[I]}\tcf(\Pi\chi_\alpha /I)$.

\item{(6)} 
$2^{|\delta^\ast|} + \sup_{J\in\calJ[I]}\sup \{ \tcf
(\Pi_{\alpha<\delta}\chi'_\alpha/J)$: $\lambda_\alpha\le
{\chi'}_\alpha = \cf (\chi'_\alpha) \le \chi_\alpha$ and the true
cofinality is well defined$\}$ $\le 2^{|\delta^\ast|} + 
\sup_{J\in\calJ[I]} \pre_J (\bar\chi, \bar\lambda)\le
2^{|\delta^\ast|} + \sup_{J\in\calJ[I]} \sup\{\tcf
(\Pi_{\alpha<\delta} \chi'_\alpha / J)\colon
|\delta^\ast| < \cf (\chi'_\alpha)$ and
$\lambda_\alpha \le \chi'_\alpha\le \chi_\alpha\}$.

\item{(7)}
In part (6), if $I$ is a precipitous ideal then the first
inequality is equality.
\endproclaim

\demo{Proof}
Straightforward.

\rem{5.9 Observation} 
In several of the models of set theory in which we know  $``\lambda $  
strong, singular, limit,  $2^\lambda > \lambda ^+"$  our sufficient
conditions for  $d_{\cf \lambda }(\lambda ,2) = 2^\lambda $ usually
hold by the sufficient condition 5.4(a) (simplest: if
$GCH$ holds below $\lambda$, $\cf\lambda=\aleph_0$).
\endproclaim

\rem{Remark} 
We could prove this consistency by looking more at the consistency 
proofs, adding many Cohen subsets to  $\lambda $  in preliminary
forcing; but the present way looks more informative.\footnote{**}{See
much more on independence in a paper of Gitik and Shelah.}



\head 6. Odds and Ends \endhead

\proclaim{6.1 Lemma} 
Suppose  $\cf (\delta ) > \kappa ^+$, $I$ an ideal on  $\kappa $,  
$f_{\alpha}  \in$  $^\kappa\!\!\Ord$  for  $\alpha < \delta$ is\break
$\leq_{I}$-increasing.  Then there are  $J_{\alpha} $, $\bar s$, 
$f'_{\alpha} (\alpha < \delta )$  such that:  
\item{(A)} $\bar s = \langle s_{i}\colon i < \kappa \rangle $, 
each  $s_{i}$ a set $of  \leq \kappa $  ordinals, 

\item{(B)} $ \bigwedge_{i < \kappa}
\bigwedge_{\alpha<\delta} \bigvee_{\beta\in {s_{i}}}f_{\alpha}
(i)\leq\beta$,

\item{(C)} $f'_{\alpha}  \in  \prod_{i<\kappa} s_{i}$ is defined by 
$f'_{\alpha} (i) = \Min [s_{i}\backslash f_{\alpha} (i)]$,

\item{(D)} $\cf[{f'_\alpha}(i)]\le\kappa$ (e.g.\
${f'_\alpha}(i)$ is a successor ordinal) implies
${f'_\alpha}(i)=f_\alpha(i)$,

\noindent
such that:

\item{(E)} $J_{\alpha} $ is an ideal on  $\kappa $  extending $I$
(for  $\alpha < \lambda )$,  decreasing with  $\alpha $  (in fact for
some  $a_{\alpha ,\beta }
\subseteq  \kappa $  (for  $\alpha < \beta < \kappa )$,  
$a_{\alpha ,\beta }/I$  decreases with  $\beta $,  increases with 
$\alpha $  and $J_{\alpha} $ is the ideal generated by  $I \cup  
\{a_{\alpha ,\beta }\colon \alpha < \beta < \lambda \}$) so
possibly $J_\alpha=\calP(\kappa)$ and possibly $J_\alpha =I$, 

\item{(F)} if  $D$  is an ultrafilter on  $\kappa $  disjoint to 
$J_{\alpha} $ then  
$f'_{\alpha} /D$  is a $<_{D}$-l.u.b  of  \break
$\langle f_{\beta} /D\colon
\beta < \delta \rangle $  and $\{i < \kappa \colon \cf [f'_{\alpha}
(i))] > \kappa \} \in  D$,

\item{(G)} if  $D$  is an ultrafilter on  $\kappa $  disjoint to $I$
but for every  
$\alpha $  not disjoint to  $J_{\alpha} $ then  $\bar s$  exemplifies  
$\langle f_{\alpha} \colon \alpha < \delta \rangle $  is chaotic for 
$D$, i.e.\ for some club $E$  of  $\delta $,  $\beta < \gamma \in  E
\Rightarrow  f_{\beta}  \leq_{D} f'_{\beta}  <_{D} f_{\gamma} $,

\item{(H)} if $\cf(\delta) > 2^\kappa$ then  $\langle f_{\alpha} \colon
\alpha < 
\delta \rangle $  has a $\leq_{I}$-l.u.b. and even $\leq_I$-e.u.b,

\item{(I)} if $b_\alpha=:\{ i\colon {f'_\alpha}(i)$ has cofinality
$\le\kappa$ (e.g.\ is a successor)$\}\not\in J_\alpha$ then: for
every $\beta\in(\alpha, \delta)$ we have
$f^{'}_\alpha\upharpoonright b_\alpha=
f_\beta\upharpoonright b_\alpha\mod J_\alpha$.

\noindent
Moreover

\item{(F)$^+$} if  $\kappa \notin  J_{\alpha} $ then  $f'_{\alpha}$ 
is an $<_{J_{\alpha}}$-e.u.b (= exact upper bound) of  $\langle
f_{\beta} \colon \beta < \delta \rangle .$
\endproclaim

\demo{Proof} 
Let  $S = \{j\colon j \leq  \sup \bigcup_{{\alpha} < {\delta}}
$Rang$(f_{\alpha} )$  has cofinality  $\leq  \kappa \}$,  $\bar e = 
\langle e_{j}\colon j \in  S\rangle $  be such that  $[j = i + 1
\Rightarrow  e_{j} = \{i\}]$,  $[j$ limit $\und j' \in  S \cap  e_{j}
\Rightarrow  e_{j'} \subseteq e_{j}]$,  $e_{j} \subseteq  j$\break
$[j$ limit $\Rightarrow  j = \sup  e_{j}]$ and  $|e_{j}| \leq 
\kappa$.

For a set  $a \subseteq  \sup \bigcup_{{\alpha} < {\delta}}
$Rang\ $(f_{\alpha})$ let  $\bar e[a] = a \cup \bigcup_{{j} \in 
{a\cap S}}e_{j}$ hence  $\bar e[\bar e[a]] = \bar e[a]$ and
$\left[ a\subseteq b\Rightarrow\bar e[a] \subseteq  \bar e[b]\right]$ 
and $|\bar e[a]|\le |a|+\kappa$.
We try to choose by 
induction on  $\zeta < \kappa ^+$,  the following:  $\alpha_{\zeta} $, 
$D_{\zeta} $, $g_{\zeta} $, $\bar s_{\zeta} 
= \langle s_{\zeta ,i}\colon i < 
\kappa \rangle $,  $\left\langle f_{\zeta ,\alpha }\colon \alpha < 
\delta \right\rangle$ such that:

\item{(a)} $g_{\zeta}  \in\  ^\kappa\!\Ord$,  

\item{(b)} $s_{\zeta ,i} = \bar e\left[ \{g_{\epsilon} (i)\colon
\epsilon 
< \zeta \} \cup  \{ \sup_{\alpha <\delta} f_{\alpha} (i) + 1\}
\right] $  so it is a set $of  \leq  \kappa $  ordinals, 
increasing with  $\zeta $,  $ \sup_{\alpha <\delta } f_{\alpha}
(i) + 1 \in  s_{\zeta ,i}$,

\item{(c)} $f_{\zeta ,\alpha } \in\  ^\kappa\!\Ord $,  $f_{\zeta,\alpha}
(i) = \Min [s_{\zeta ,i}\backslash f_{\alpha} (i)]$,

\item{(d)} $D_{\zeta} $ is an ultrafilter on  $\kappa$ disjoint to $I$,

\item{(e)} for  $\alpha < \delta $,  $f_{\alpha}  \leq_{D_{\zeta}}
g_{\zeta}$,

\item{(f)} $\alpha_{\zeta} $ is an ordinal  $< \delta$,

\item{(g)} $\alpha_{\zeta}  \leq  \alpha < \lambda \Rightarrow 
g_{\zeta}  <_{D_{\zeta} }f_{\zeta ,\alpha }$.

If we succeed, let  $\alpha (\ast ) =  \sup_{\zeta <\kappa^+}
\alpha_{\zeta} $,  so as  $\cf(\delta)> \kappa ^+$ clearly 
$\alpha (\ast ) < \delta $.  Now let  $i < \kappa $  and look at  
$\langle f_{\zeta ,\alpha (\ast )}(i)\colon \zeta < \kappa ^+\rangle$;
by
its definition (see (c)),  $f_{\zeta ,\alpha (\ast )}(i)$  is the
minimal member of the set  $s_{\zeta ,i}\backslash f_{\alpha(\ast)}(i)$.
This set increases with  $\zeta $,  so  $f_{\zeta ,\alpha (\ast )}(i)$ 
decreases with  $\zeta $  (though not necessarily strictly), hence is
eventually constant; so for some  $\zeta_{i} < \kappa ^+$ we have 
$\zeta \in  [\zeta_{i},\kappa ^+) \Rightarrow  
f_{\zeta ,\alpha (\ast )}(i) = f_{\zeta_{i},\alpha (\ast )}(i)$. 
Let  $\zeta (\ast ) =  \sup_{i<\kappa } \zeta_{i}$,  so  $\zeta
(\ast ) < \kappa ^+$,  hence 
$$
\zeta \in  [\zeta (\ast ),\kappa ^+) \Rightarrow   
\bigwedge_{i} f_{\zeta ,\alpha (\ast )}(i) = f_{\zeta (\ast ),\alpha
(\ast )}(i) \Rightarrow 
f_{\zeta ,\alpha (\ast )} = f_{\zeta (\ast ),\alpha (\ast )}.
\leqno(*)
$$
We know that  $f_{\alpha (\ast )} \leq_{D_{\zeta (\ast )}}g_{\zeta
(\ast )} <_{D_{\zeta (\ast )}}f_{\zeta (\ast ),\alpha (\ast )}$ hence
for some  $i$,  $f_{\alpha (\ast )}(i) \leq  g_{\zeta (\ast )}(i) < 
f_{\zeta (\ast ),\alpha (\ast )}(i)$,  but  $g_{\zeta (\ast )}(i) \in  
s_{\zeta (\ast )+1,i}$ hence  $f_{\zeta (\ast )+1,\alpha (\ast )}(i)
\leq  g_{\zeta (\ast )}(i) < f_{\zeta (\ast ),\alpha (\ast )}(i)$,
contradicting the choice of  $\zeta (\ast ).$ 

So necessarily for some  $\zeta < \kappa ^+$ we are stuck, 
and clearly  $s_{\zeta ,i}(i < \kappa )$,  
$f_{\zeta ,\alpha }(\alpha < \lambda )$  are well defined. 

Let  $s_{i} =: s_{\zeta ,i}$ (for  $i < \kappa )$  and $f'_{\alpha}  = 
f_{\zeta ,\alpha }$ (for  $\alpha < \lambda )$.
Clearly  $s_{i}$ is a set of $\leq\kappa$ ordinals; now clearly: 
\item{$(\ast )_{1}$}  $f_{\alpha}  \leq  f'_{\alpha} $

\item{$(\ast )_{2}$}  $\alpha < \beta \Rightarrow  f'_{\alpha} 
\leq_{I} f'_{\beta} $,

\item{$(\ast)_{3}$} if $b=\{i\colon f'_{\alpha} (\alpha ) < f'_{\beta}
(i)\} \notin  I$,  $\alpha < \beta < \delta $  then 
$f'_{\alpha} \upharpoonright b <_{I} f_{\beta} \upharpoonright b$.

We let for  $\alpha < \delta $ 
$$
\eqalign{
J_{\alpha}  = \Big\{ b \subseteq  \kappa \colon b & \in  I
\hbox{  or  } b \notin  I \quad\hbox{ and for some } \beta \
\hbox{ we have: } \alpha < \beta < \delta \ \hbox{ and } \cr
&  f'_{\alpha} \upharpoonright (\kappa\setminus b) =_{I} f'_{\beta}
\upharpoonright (\kappa\setminus b)
\Big\} .\cr}
$$

We let for $\alpha<\beta<\delta$, $a_{\alpha,\beta}=:\{
i<\kappa\colon {f'_\alpha}(i)<{f'_\beta}(i)\}$.
Then

\item{$(\ast )_{4}$} $J_{\alpha} $ is an ideal on  $\kappa $ 
extending $I$, in fact is the ideal generated by $I\cup\{
a_{\alpha,\beta}\colon\beta\in(\alpha,\delta)\}$.

As $\langle f'_{\alpha}\colon \alpha<\delta \rangle$ is
$\leq_{I}$-increasing (i.e.\ $(*)_1$): 

\item{$(\ast )_{5}$}  $J_{\alpha} $ decreases with  $\alpha$, in fact 
$a_{\alpha,\beta}/I$ increases with $\beta$, decreases with
$\alpha$,

\item{$(\ast )_{6}$}  if  $D$  is an ultrafilter on  $\kappa $ 
disjoint to  $J_{\alpha}$, then  $f'_{\alpha} /D$  is a $<_{D}$-lub
of\break  $\{f_{\beta} /D\colon
\beta\shvor < \delta \}$.

[Why?
We know that $\beta\in(\alpha,\delta)\Rightarrow
a_{\alpha,\beta}=\emptyset\mod D$,
 so  $f_{\beta}  \leq 
f'_{\beta}  =_{D} f'_{\alpha} $  for  $\beta \in  (\alpha ,\delta )$,
so  $f'_{\alpha} /D$  is an $\leq_{D}$-upper bound.
If it is not a least upper bound then for some  $g \in\  
^\kappa\!\Ord $,  $ \bigwedge_{\beta } f_{\beta}  \leq_{D}
g <_{D} f'_{\alpha} $  and we can get a contradiction to the 
choice of  $\zeta $, $\bar s$, $f'_{\beta} $  as:  $(D,g)$  could
serve as  $D_{\zeta} $,\ $g_{\zeta} .]$ 

\item{$(\ast )_{7}$} If  $D$  is an ultrafilter on  $\kappa $  disjoint
to $I$ but not to 
$J_{\alpha} $ (for\ every\ \ $\alpha < \lambda )$  then  $\bar s$ 
exemplifies  
$\langle f_{\alpha} \colon \alpha < \delta \rangle $  is chaotic 
for  $D.$   

[Why?  For every  $\alpha < \delta $  for some  $\beta \in  
(\alpha ,\delta )$  we have $a_{\alpha,\beta}\in D$, i.e.\break
$\{i < \kappa \colon f'_{\alpha} (i) < f'_{\beta}
(i)\} \in D$,  so  $\langle f'_{\alpha} /D\colon \alpha < \delta \rangle $ 
is not eventually constant, so if  $\alpha < \beta $,  $f'_{\alpha}
<_{D} f'_{\beta} $  then  $f'_{\alpha}  <_{D} f_{\beta} $ (by
$(\ast )_{3})$ and $f_{\beta}  \leq_{D} f'_{\beta} $  
(by (c)) as required.]

\item{$(\ast )_{8}$} if $\kappa\not\in J_\alpha$ then $f^\prime_\alpha$
is an $\leq_{J_\alpha}$-e.u.b. of $\langle f_\beta\colon
\beta<\delta\rangle$.

[Why?
By $(\ast)_6$, $f^\prime_\alpha$ is a $\leq_{J_\alpha}$-upper bound
of $\langle f_\beta\colon \beta<\delta\rangle$;
so assume that it is not a $\leq_{J_\alpha}$-e.u.b. of $\langle
f_\beta\colon \beta<\delta\rangle$, hence there is a function
$g$ with domain $\kappa$, such that $g(i)<\Max\{1,f^\prime_\alpha(i)\}$, but
for no $\beta<\delta$ do we have
$$
C_\beta =:\{ i<\kappa \colon
g(i) <\Max\{1,f_\beta(i)\}=\kappa\mod J_\alpha.
$$
Clearly $\langle C_\beta\colon\beta<\delta\rangle$ is increasing
modulo $J_\alpha$ so there is an ultrafilter $D$ on $\kappa$
disjoint to $J_\alpha\cup\{ C_\beta\colon\beta<\delta\}$.
So $f_\beta\leq_D g\leq_D f^\prime_\alpha$, so we get a contradiction
to $(\ast)_6$ except when $g=_D f^\prime_\alpha $ and then
$f^\prime_\alpha =_D O_\kappa$ (as $g(i)< 1\ \vee\ g(i)< f'_\alpha(i)$).
If we can demand $b^*=\{i\colon f^\prime_\alpha(i) =0\}\notin D$ we are
done, but easily $b^*\hefresh C_\beta\in J_\alpha$ so we finish.]

\item{$(\ast )_{9}$} If $\cf [f^\prime_\alpha(i)]\le\kappa$ then
$f^\prime_\alpha(i) = f_\alpha(i)$.

[Why? By the definition of $s_\zeta =\bar e[\dots]$ and the choice of
$\bar e$, and $f^\prime_\alpha(i)$.]

\item{$(\ast )_{10}$} Clause (I) of the conclusion holds.

[Why?
As $f_\alpha\le_{J_\alpha} f_\beta\le_{J_\alpha} f^\prime_\alpha$ and
$f_\alpha\upharpoonright b =_{J_\alpha}
f^\prime_\alpha\upharpoonright b$ by $(\ast)_9$.]

The reader can check the rest.  
\sqed{6.1}

\rem{6.1A Example} 
We show that l.u.b and e.u.b are not the same.  Let $I$ be
an ideal on  $\kappa $,  $\kappa ^+ < \lambda = \cf(\lambda)$, 
$\bar a = \langle a_{\alpha} \colon \alpha < \lambda \rangle $  be
a sequence of subsets of  $\kappa $,  (strictly) increasing modulo $I$, 
$\kappa \backslash a_{\alpha}  \notin  I$  but there is no 
$b \in  {\cal P}(\kappa )\backslash I$  such that 
$ \bigwedge_{\alpha } b \cap  a_{\alpha}  \in  I$. 
[Does this occur? E.g.\ for  $I = {\cal S}_{<\aleph_{0}}(\omega )$, 
the existence of such  $\bar a$  is known to be consistent; e.g.\
MA\ $\und \kappa =\aleph_0\und\lambda = 2^{\aleph_{0}}$. 
Moreover, for any  $\kappa $  and $\kappa ^+ < \lambda = \cf \lambda
\leq  2^\kappa $ we can find $a_{\alpha}  \subseteq  
\kappa $  for  $\alpha < \lambda $  such that, e.g., any Boolean 
combination of the $a_{\alpha}$'s  has cardinality  $\kappa $ 
(less needed). 
Let  $I_{0}$ be the ideal on  $\kappa $  generated by  
${\cal S}_{<\kappa }(\kappa ) \cup  \{a_{\alpha} \backslash a_{\beta}
\colon \alpha < \beta < \lambda \}$,  and let $I$ be maximal in  
$\{J\colon J$\ \ an\ ideal\ on\ \ $\kappa $,\ \ $I_{0} \subseteq  
J$\ \ and\ \ $[\alpha < \beta < \lambda \Rightarrow  
a_{\beta} \backslash a_{\alpha}  \notin  J]\}$. 
So if G.C.H. fails, we have examples.] 
For  $\alpha < \lambda $,  we let  $f_{\alpha} \colon \kappa
\rightarrow \Ord$  be:     
$$
f_{\alpha} (i) = \cases{ \alpha & if $\alpha\in\kappa\setminus a_i$,\cr
\lambda+\alpha & if $\alpha\in a_i$.}
$$

Now the constant function $f\in$ ${}^\kappa\!\Ord$,
$f(i) = \lambda + \lambda $  is a l.u.b of 
$\langle f_{\alpha} \colon \alpha < \lambda \rangle $  but not an
e.u.b.\ (both $\mod J$) (not e.u.b. is exemplified by $g\in
{}^\kappa\!\Ord$ which is constantly $\lambda$).

\proclaim{6.2 Claim} 
Suppose  $\mu  > \kappa = \cf \mu $,  $\mu  = \tlim_{J}\lambda_{i}$,  
$\delta < \mu $,  $\lambda_{i} = \cf(\lambda_{i})> \delta $ for
$i<\delta$,  $J$  a 
$\sigma $-complete ideal on  $\delta $  and $\lambda =
\tcf\left( \prod_{{i} < {\delta}} \lambda_{i}/J\right) $, 
and $\langle f_{\alpha}\colon\alpha <\lambda\rangle$ exemplifies this. 

Then we have 
\item{$(\ast )$} if  $\langle u_{\beta} \colon \beta < \lambda \rangle $  is
a sequence of pairwise disjoint non-empty subsets of  $\lambda $, 
each of cardinality  $\le \sigma $ (not $<\sigma$!) and
$\alpha ^\ast  < \mu$,  then we
can find $B \subseteq  \lambda $  such that:   

\itemitem{(a)} $\otp (B) = \alpha ^\ast $,   

\itemitem{(b)} if  $\beta \in  B$,  $\gamma \in  B$  and $\beta <
\gamma $  then  $\sup  u_{\beta}  < \min  u_{\gamma} $,   

\itemitem{(c)} we can find $s_{\zeta}  \in  J$  for  $\zeta \in
\bigcup_{{i} \in  {B}}u_{i}$ such that:
if \ $\zeta \in \bigcup_{{\beta} \in  {B}}u_{\beta} $, 
$\xi  \in \bigcup_{{\beta} \in  {B}}u_{\beta} $,  $\zeta < \xi$ 
and $i \in  \delta \backslash s_{\zeta} \backslash s_{\xi}$, then 
$f_{\zeta} (i) < f_{\xi} (i).$
\endproclaim

\demo{Proof} 
For each regular  $\theta, \theta^+ < \mu $,  there is a stationary  $S_{\theta}
\subseteq  \{\delta < \lambda \colon \cf(\delta)= \theta < \delta \}$ 
which is in  $I[\lambda ]$  (see [Sh420, 1.5]) which is equivalent (see
[Sh420, 1.2(1)]) to: 
\item{$(\ast )$} there is  $\bar C^\theta = \langle C^\theta_{\alpha}
\colon i < \lambda \rangle ,$   

\itemitem{$(\alpha )$} $C^\theta_{\alpha} $ a subset of  $\alpha $, 
with no accumulation points (in  $C^\theta_{\alpha} )$,   

\itemitem{$(\beta )$} $[\alpha \in  \nacc(C^\theta_{\beta}) 
\Rightarrow  C^\theta_{\alpha} = C^\theta_{\beta}  \cap  \alpha ]$,

\itemitem{$(\gamma )$} for some club $E^0_{\theta} $ of  $\lambda $,
$$
[\delta \in  S_{\theta}  \cap  E^0_{\theta}  \Rightarrow  \cf(\delta)=
\theta < \delta \und \delta = \sup  C^\theta_{\delta}  \und \otp 
(C^\theta_{\delta}) = \theta ].
$$

\noindent
Without loss of generality  $S_{\theta}  \subseteq  E^0_{\theta} $, 
and $ \bigwedge_{{\alpha} < {\delta}} \otp(C^\theta_{\delta})
\leq  \theta $.  By [Sh365, 2.3, Def.\ 1.3] for some club 
$E_{\theta} $ of  $\lambda $,  
$\langle \gl (C^\theta_{\alpha} ,E_{\theta} )\colon \alpha \in 
S_{\theta} \rangle $ guess clubs  (i.e.\ for every club
$E \subseteq  E_{\theta} $ of  $\lambda $,  
for stationarily many  $\zeta \in  S_{\theta} $,  
$\gl (C^\theta_{\zeta} ,E_{\theta} ) \subseteq  E)$  (remember  
$\gl (C^\theta_{\delta} ,E_{\theta} ) = \{\sup (\gamma \cap 
E_{\theta} )\colon \gamma \in  C^\theta_{\delta} ;\gamma >
\Min (E_{\theta})\})$.  Let  
$C^{\theta ,\ast }_{\alpha}  = \{\gamma \in  C^\theta_{\alpha} \colon
\gamma =
\Min (C^\theta_{\alpha} \backslash \sup (\gamma \cap  E_{\theta} )\}$, 
they have all the properties of the  $C^\theta_{\alpha}$'s  and guess
clubs in a weak sense: for every club $E$  of  $\lambda $  for some 
$\alpha \in  S_{\theta}  \cap  E$,  if  $\gamma_{1} < \gamma_{2}$ are
successive members of  $E$  then  $|(\gamma_{1},\gamma_{2}] \cap 
C^{\theta ,\ast }_{\alpha} | \leq  1$;  moreover, 
the function  $\gamma \mapsto  \sup (E \cap  \gamma )$  is one to one on
$C^{\theta ,\ast }_{\zeta} .$ 

Now we define by induction on  $\zeta < \lambda $,  an ordinal  
$\alpha_{\zeta} $ and functions  $g^\zeta_{\theta}  \in \prod_{{i}
< {\delta}} \lambda_{i}$ (for each  $\theta \in  {\Theta } =: 
\{\theta \colon \theta < \mu,$ $\theta$ regular uncountable$\}).$ 

For given  $\zeta $,  let  $\alpha_{\zeta}  < \lambda $  be minimal
such that:
$$
\eqalign{
\xi  < \zeta & \Rightarrow  \alpha_{\xi}  < \alpha_{\zeta}, \cr
\xi  < \zeta \und \theta \in  {\Theta } & \Rightarrow  g^\zeta_{\theta}  <
f_{\alpha_{\zeta} }\hbox{ mod }\ J.\cr}
$$
Now  $\alpha_{\zeta} $ exists as  $\langle f_{\alpha} \colon \alpha < \lambda
\rangle$  is $<_{J}$-increasing cofinal in  $ \prod_{{i} <
{\lambda_i}} /J$.  Now for each  $\theta \in  {\Theta }$  we define  
$g^\zeta_{\theta} $ as follows: 
\item{} for  $i < \delta ^\ast $,  $g^\zeta_{\theta} (i)$  is  
$\sup \left[ \{g^\xi_{\theta} (i) + 1\colon \xi  \in  C^\theta_{\zeta} \} \cup
\{f_{\alpha_{\zeta} }(i) + 1\}\right] $ if this number is 
$< \lambda_{i}$, and $f_{\alpha_{\zeta} }(i)$  otherwise.

Having made the definition we prove the assertion. 
We are given\break
$\langle u_{\beta} \colon \beta < \lambda \rangle $, 
a sequence of pairwise
disjoint non-empty subsets of  $\lambda $,  each of cardinality 
$< \sigma $  and $\alpha ^\ast  < \mu $. 
We should find $B$ as promised; let  $\theta =: (|\alpha^* | + |\delta
|)^+$ so  $\theta < \mu $  is regular  $> |\delta |$.  
Let  $E = \{\delta \in  E_{\theta}:$ for every $\zeta \colon
[\zeta <
\delta \Leftrightarrow  \sup  u_{\zeta}  < \delta \Leftrightarrow 
u_{\zeta}  \subseteq  \delta\Leftrightarrow\alpha_\zeta <\delta ]\}$. 
Choose  $\alpha \in  S_{\theta}  \cap  \acc(E)$ such that
$\gl(C_\zeta^\theta, E_\theta)\subseteq E$;  hence letting 
$C_{\alpha}^{\theta,*}  = \{\gamma_{i}\colon i < \theta \}$ 
(increasing) we know  $ \bigwedge_{i} (\gamma_{i},\gamma_{i+1}) \cap 
E \neq \emptyset $. 
Now  $B=\left\lbrace \gamma_{5i+3}\colon i < \alpha^*
\right\rbrace $  are as required. 
For $\alpha\in \bigcup_{\zeta<\alpha^*} u_{5\zeta + 3}$ let
$s_\alpha=s_\alpha^o\cup s^1_\alpha$.
For $\alpha\in u_{5\zeta + 3}$, $\zeta<\alpha^*$, let
$s_\alpha^o = \{ i<\delta\colon g_\theta^{5\zeta + 1} (i) <
f_\alpha(i) < g^{5\zeta + 4} (i)\}$, for each $\zeta <
\alpha^*$;
let $\langle\alpha_\epsilon\colon\epsilon<|u_{5\zeta +
3}|\rangle$ enumerate $u_{5\zeta + 3}$ and 
$$
s^1_{\alpha_\epsilon} = \{ i\colon \hbox{ for every }
\xi<\epsilon, f_{\alpha_\xi}(i)< f_{\alpha_\epsilon}(i)
\Leftrightarrow \alpha_\xi <\alpha_\epsilon \Leftrightarrow
f_{\alpha_\xi}(i) \le f_{\alpha_\epsilon}(i) \}.  \sqed{6.2}
$$

\rem{6.2A Remark} 
In 6.2:
(1) We can avoid guessing clubs.

(2)  Assume  $\sigma < \theta_{1} < \theta_{2} < \mu $  are regular
and there is 
$S\subseteq\{\delta <\lambda\colon\cf(\delta)= \theta_{1}\}$  from  
$I[\lambda ]$  such that for every  $\zeta < \lambda $ (or at least
a club) of cofinality  $\theta_{2}$,  $S \cap  \zeta $  is stationary
and $\langle f_\alpha\colon\alpha<\lambda\rangle$ obey suitable
$\bar C^\theta$ (see [Sh345a, \S2]). 
Then for some  $A \subseteq  \lambda $  unbounded, for every 
$\langle u_{\beta} \colon \beta < \theta_{2}\rangle $ 
sequence of pairwise disjoint non-empty subsets of  A,  each of
cardinality  $< \sigma $  with  $[\min  u_{\beta} ,\sup  u_{\beta} ]$ 
pairwise disjoint we have:  for every  $B_{0} \subseteq  A$  of order
type  $\theta_{2}$,  for some  $B \subseteq  B_{0}$, $|B| = \theta_{1}$,
(c) of $(\ast )$ of 6.2 holds.

(3)  In $(\ast )$ of 6.2,  $``\alpha ^\ast  < \mu "$  can be replaced by
``$\alpha ^\ast  < \mu ^+$'' (prove by induction on  $\alpha^*$).

\proclaim{6.3 Observation} 
Assume  $\lambda < \lambda ^{<\lambda }$,  $\mu  = 
\Min \{\mu \colon 2^\mu  > \lambda \}$. 
Then there are  $\delta $, $\chi$ and ${\cal T}$,
satisfying the condition $(\ast )$ below for  $\chi  = 2^\mu $ or at
least arbitrarily large regular  $\chi  \leq  2^\mu .$ 
\item{$(\ast )$} ${\cal T}$  a tree with  $\delta $  levels,
(where $\delta \leq \mu )$  with a set  $X$ of $\geq  \chi$
$\delta$-branches, and for  $\alpha < \delta $, 
$ \bigcup_{{\beta} < {\alpha}} |{\cal T}_{\beta} | < \lambda .$
\endproclaim

\demo{Proof of Observation} 
So let  $\chi  \leq  2^\mu $ be regular,  $\chi  > 
\lambda .$

\subdemoinfo{Case 1}%
{$ \bigwedge_{{\alpha} < {\mu}} 2^{|\alpha |} < \lambda$}
Then  ${\cal T} =\ ^{\mu >} 2$,  $\calT_{\alpha}  =\ ^\alpha 2$
are O.K.
(the set of branches  $^\mu 2$  has cardinality  $2^\mu )$.

\subdemoinfo{Case 2}{Not Case 1} 
So for some  $\theta < \mu $,  $2^\theta \geq  
\lambda $,  but by the choice of  $\mu $,  $2^\theta \leq  \lambda $, 
so  
$2^\theta = \lambda $,  $\theta < \mu $  and so  $\theta \leq  \alpha < 
\mu  \Rightarrow  2^{|\alpha |} = 2^\theta $. 
Note $|^{\mu >}2| = \lambda$ as $\mu\le\lambda$.

\subdemoinfo{Subcase 2A}{$\cf(\lambda)\neq  \cf(\mu)$} 
Let  $^{\mu >}2 = \bigcup_{{j} < {\lambda}} B_{j}$, 
$B_{j}$ increasing with  $j$,  $|B_{j}| < \lambda $.  For each  
$\eta \in\  ^\mu 2$,  (as  $\cf(\lambda)\neq  \cf(\mu)$)  for some 
$j_{\eta}  < \lambda ,$
$$
\mu  = \sup \left\lbrace \zeta < \mu \colon \eta \upharpoonright \zeta \in  
B_{j_{\eta} }\right\rbrace .
$$
So as  $\cf (\chi ) > \mu $,  for some ordinal  $j^\ast  < \lambda $ 
we have
$$
\{\eta \in\  ^\mu 2\colon j_{\eta}  \leq  j^\ast \}
\hbox{  has cardinality } \ge\chi. 
$$
As  $\cf(\lambda)\neq  \cf(\mu)$  and $\mu  \leq  \lambda $  (by its
definition) clearly  $\mu  < \lambda $,  hence  $|B_{j^\ast }| \times 
\mu  < \lambda $.  
Let
$$
{\cal T} = \left\lbrace \eta \upharpoonright \epsilon \colon \epsilon 
< \lg(\eta)\hbox{  and } \eta \in  B_{j^\ast }\right\rbrace .
$$
It is as required.

\subdemoinfo{Subcase 2B}{Not 2A so  $\cf(\lambda)= \cf(\mu)$}
As  $(\forall \sigma )[\theta \leq  \sigma < \mu  \Rightarrow  \lambda =
2^\sigma \Rightarrow  \cf(\lambda)= \cf (2^\sigma ) > \sigma ]$, 
clearly  $\cf(\lambda)\geq  \mu $  so  $\mu $  is regular. 
If  $\lambda = \mu $  we get
$\lambda = \lambda ^{<\lambda }$ contradicting an assumption. 

So  $\lambda > \mu $,  so  $\lambda $  singular. 
So if  $\alpha < \mu $,  $\mu  < \sigma_{i} = \cf(\sigma_{i})< \lambda$
for  $i < \alpha $  then (see [Sh-g, 345a, 1.3(10)])
$\max\pcf \{\sigma_{i}\colon i <
\alpha \} \leq \prod_{{i} < {\alpha}} \sigma_{i} \leq 
\lambda^{|\alpha |} \leq  
\left( 2^\theta \right) ^{|\alpha |} \leq  2^{<\mu } = \lambda $, 
but as  $\lambda $  is singular and $\max  \pcf \{\sigma_{i}\colon
i <
\alpha \}$  is regular 
(see [Sh345a, 1.9]),  clearly the inequality is strict, i.e.\  $\max  
\pcf \{\sigma_{i}\colon i < \alpha \} < \lambda $.  So let  $\langle
\sigma_{i}\colon i < 
\mu \rangle $  be a strictly increasing sequence of regulars in  
$(\mu ,\lambda )$  with limit  $\lambda $,  and by [Sh355, 3.4] there
is 
$T \subseteq \prod_{{i} < {\mu}} \sigma_{i}$,  $|\{\nu 
\upharpoonright i\colon
\nu  \in  T\}| \leq  \max  \pcf \{\lambda_{j}\colon j
< i\} < \lambda $,  and number of $\mu $-branches  $> \lambda $. 
In fact we can get any regular cardinal in  $(\lambda,\pp^+(\lambda ))$
in the same way.
Let $\lambda^\ast = \min\{\lambda' \colon \mu <\lambda'\le \lambda$,
$\cf(\lambda')=\mu$ and $pp(\lambda')>\lambda\}$, so (by [Sh355,
2.3]), also $\lambda^\ast$ has those properties and
$\pp(\lambda^\ast)\ge \pp(\lambda)$.
So if  $\pp^+(\lambda^\ast) = (2^\mu )^+$ or $\pp(\lambda^\ast)
 = 2^\mu$ is singular, we are done. 
So assume this fails.

If  $\mu  > \aleph_{0}$, then (as in 3.4)  $\alpha <
2^\mu  \Rightarrow  \cov(\alpha ,\mu ^+,\mu ^+,\mu ) < 2^\mu $ and we
can finish as in subcase 2A (as in 3.4; actually $\cov
(2^{<\mu}, \mu^+, \mu^+, \mu)<2^\mu$ suffices which holds by the
previous sentence and [Sh355, 5.4]). 
If  $\mu  = \aleph_{0}$ all is easy. $_{\qed_{6.3}}$

\proclaim{6.4 Claim} 
Assume  $\frb_{k} \subseteq  \frb_{k+1} \subseteq \cdots$  for  $k < 
\omega $,  $\fra = \bigcup_{{k} \leq  {\omega}} \frb_{k}$  (and
$|\fra| < \Min\fra $)  and $\lambda \in  \pcf \fra\backslash
\bigcup_{{k} < {\omega}} \pcf (\frb_k).$
\item{(1)}  Then we can find finite  $\frd_{k} \subseteq  \pcf (\frb_{k}
\backslash \frb_{k-1})$  (stipulating  $\frb_{-1} = \emptyset )$ 
such that $\lambda \in  \pcf \bigcup_{{k} < {\omega}} \frd_{k}.$

\item{(2)} Moreover, we can demand $\frd_{k} \subseteq  (\pcf
\frb_{k})\backslash (\pcf  (\frb_{k-1})).$
\endproclaim

\demo{Proof} 
We start to repeat the proof of [Sh371, 1.5] for  $\kappa = \omega $.  
But there we apply [Sh371, 1.4] to  $\langle \frb_{\zeta} \colon
\zeta < \kappa \rangle $ 
and get  $\langle \langle\frc_{\zeta ,\ell }\colon \ell  \leq  n_\zeta 
\rangle \colon \zeta < 
\kappa \rangle $  and let  $\lambda_{\zeta ,\ell } = \max  
\pcf(\frc_{\zeta ,\ell})$. 
Here we apply the same claim ([Sh371, 1.4]) to  $\langle
\frb_{k}\backslash\frb_{k-1}\colon k<\omega\rangle$  to get part (1). 
As for part (2), in the proof of [Sh371, 1.5] we let $\delta =|\fra|^+ +
\aleph_{2}$ choose $\langle N_{i}\colon i<\delta\rangle$,  but now we
have to adapt the proof of [Sh371, 1.4] (applied to $\fra$,
$\langle \frb_{k}\colon k < \omega \rangle $, $\langle N_{i}\colon i
< \delta \rangle )$;  we have gotten there, toward the end, 
$\alpha < \delta $  such that  $E_{\alpha}  \subseteq  E$. 
Let  $E_{\alpha}  = \{i_{k}\colon k < \omega \}$,  $i_{k} < i_{k+1}$. 
But now instead of applying 
[Sh371, 1.3] to each  $\frb_{\ell} $ separately,  we try to choose  
$\langle c_{\zeta ,\ell }\colon \ell  \leq  n(\zeta )\rangle $  by induction
on  $\zeta < \omega $. 
For  $\zeta = 0$  we apply [Sh371, 1.3].  For  $\zeta > 
0$,  we apply [Sh371, 1.3] to  $\frb_{\zeta} $ but there defining by
induction on  
$\ell \ \ \frc_{\ell}=\frc_{\zeta,\ell}\subseteq\fra$  such that  
$\max \left( \pcf (\fra\backslash \frc_{\zeta ,0}\backslash \cdots
\backslash
\frc_{\zeta , \ell -1}) \cap  \pcf \frb_{\zeta} \right) $  is
strictly decreasing with  $\ell $.  
We use:

\rem{6.4A Observation} 
If $|\fra_i|<\Min(\fra_i)$ for $i<i^*$, then
$\frc=\bigcap_{i< i^*}\pcf(\fra_i)$ has a last element or is
empty.
\smallskip

\demo{Proof}
Wlog $\langle |\fra_i|\colon i< i^*\langle $ is
nondecreasing.
By [Sh345b, 1.12]
$$
\frd\subseteq\frc\ \&\  |\frd|<\Min\frd\Rightarrow
\pcf(\frd)\subseteq\frc .
\leqno{(*)_1}
$$
By [Sh371, 2.6]
$$
\displaylines{
\hbox{ if } \lambda\in\pcf(\frd),\ \frd\subseteq \pcf(\frc),
\ |\frd|<\Min(\frd)\hbox{ then }\cr
\hbox{ for some } \fre\subseteq \frd \hbox{ we have }
|\fre|\le\Min|\fra_0|, \ \lambda\in\pcf(\fre).\cr}
$$
Now choose by induction on $\zeta<|\fra_0|^+$, $\theta_\zeta\in\frc$,
satisfying $\theta_\zeta>\max\pcf\{\theta_\epsilon\colon \epsilon<\zeta\}$.
If we are stuck in $\zeta$,
$\max \pcf\{\theta_\epsilon\colon\epsilon<\zeta\}$ is the desired maximum by
$(\ast)_1$.
If we succeed $\theta=\max\pcf\{\theta_\epsilon\colon
\epsilon<|\fra_0|^+\}$ is in
$\pcf\{\theta_\epsilon\colon\epsilon<\zeta\}$ for some
$\zeta<|\fra_0|^+$ by $(\ast)_2$; easy contradiction.
\sqed{6.4A}

\ \sqed{6.4}

\rem{6.5 Conclusion} 
Assume  $\aleph_{0} = \cf(\mu) \leq  \kappa \leq  \mu_{0} < 
\mu $,  $[\mu ' \in  (\mu_{0},\mu ) \und \cf(\mu')\leq  \kappa
\Rightarrow  \pp_{\kappa} (\mu ') <\lambda]$ and $\pp^+_{\kappa}(\mu)
> \lambda = \cf(\lambda)> \mu $. 
Then we can find $\lambda_{n}$ for  $n < \omega $,  
$\mu_{0} < \lambda_{n} < \lambda_{n+1} < \mu $,  $\mu  =
\bigcup_{{n} < {\omega}} \lambda_{n}$ and $\lambda = \tcf
 \prod_{n<\omega} \lambda_{n}/J$  for some ideal  
$J$  on  $\omega $  (extending  $J^{bd}_{\omega} ).$

\demo{Proof} 
Let $\fra\subseteq (\mu,\mu)\cap\Reg$, $|\fra|\leq\kappa $,  
$\lambda \in  \pcf(\fra)$. Without loss of generality  $\lambda = \max  
\pcf\fra$,  let  $\mu  = \bigcup_{{n} < {\omega}} \mu ^0_{n}$, 
$\mu_{0} \leq  \mu ^0_{n} < \mu ^0_{n+1} < \mu $, 
let  $\mu ^1_{n} = \mu^0_n +\sup \{\pp_{\kappa} (\mu ')\colon \mu_{0} <
\mu ' \leq \mu ^0_{n}$ and $\cf (\mu ') \leq  \kappa \}$, 
by [Sh355, 2.3] 
$\mu ^1_{n} < \mu $,  $\mu ^1_{n} = \mu^0_n + \sup \{\pp_{\kappa}
(\mu ')\colon \mu_{0}
< \mu ' < \mu ^1_{n}$ and $\cf (\mu ') \leq  
\kappa \}$  and obviously  $\mu ^1_{n} \leq  \mu ^1_{n+1}$; 
by replacing by a 
subsequence without loss of generality  $\mu ^1_{n} < \mu ^1_{n+1}$. 
Now let  $\frb_{n} = \fra \cap  \mu ^1_{n}$ and apply the previous
claim:
to  $\frb_{k} =: \fra \cap  (\mu ^1_{n})^+$,  note:
$$
\max\pcf(\frb_{k}) \le\mu^1_{k} <\Min(\frb_{k+1}\backslash \frb_{k}).
\sqed{6.5}
$$

\proclaim{6.6 Claim} 
\item{(1)} Assume  $\aleph_{0} < \cf(\mu) = \kappa < \mu_{0} < \mu $,  
$2^\kappa < \mu $  and $[\mu_{0} \leq  \mu ' < \mu  \und \cf(\mu')\leq 
\kappa \Rightarrow  \pp_\kappa\mu ' < \mu ]$. 
If  $\mu  < \lambda =\cf(\lambda) < \pp^+(\mu )$  then there is a
tree  ${\cal T}$ with $\kappa$ levels, each level of cardinality 
$<\mu$,  ${\cal T}$  has exactly  $\lambda \ \ \kappa $-branches.

\item{(2)}  Suppose  $\langle \lambda_{i}\colon i < \kappa \rangle $  is
a strictly increasing sequence of regular cardinals,  $2^\kappa <
\lambda_{0}$,  $\fra =: \{\lambda_{i}\colon i
< \kappa \}$,  $\lambda = \max  \pcf \fra$,  $\lambda_{j} > \max  
\pcf \{\lambda_{i}\colon i < j\}$  for each  $j < \kappa $  (or at least
$\sum_{i<\kappa }\lambda_{i} > \max  \pcf \{\lambda_{i}\colon i < j\})$ 
and $\fra \notin J$  where  $J = \{\frb \subseteq  \fra\colon \frb$ 
is the union of
countably many members of $J_{<\lambda }[\fra]\}$  (so  $J  
 \supseteq J^{bd}_{\fra}$, $\cf \kappa > \aleph_{0})$. 
Then the conclusion of (1) holds with 
$\mu = \sum_{i<\kappa }\lambda_{i}.$
\endproclaim

\demo{Proof} 
(1) By (2) and [Sh371, \S 1] (or can use the conclusion of [Sh-g,
AG 5.7]).

(2)  For each  $\frb \subseteq  \fra$  define the function 
$g_{\frb}\colon \kappa \rightarrow  \Reg$  by 
$$
g_{\frb}(i) =\max
\pcf [\frb \cap  \{\lambda_{j}\colon j < i\}].
$$ 
Clearly $[\frb_{1}\subseteq \frb_{2}\Rightarrow g_{\frb_{1}}\leq
g_{\frb_{2}}]$. 
As  $\cf (\kappa ) > \aleph_{0}$,\ $J\  \aleph_{1}$-complete,  there is
 $\frb \subseteq  \fra$,  $\frb \notin  J$  such that:
$$
\frc\subseteq\frb\und\frc\notin J\Rightarrow \neg g_{\frc}<_{J} g_{\frb}.
$$
Let  $\lambda ^\ast_{i} = \max  \pcf (\frb \cap  \{\lambda_{j}\colon
j < i\})$. 
For each  $i$  let  $\frb_{i} = \frb \cap  \{\lambda_{j}\colon j < i\}$
and $\langle \langle f^{\frb}_{\lambda ,\alpha }\colon \alpha <
\lambda \rangle \colon \lambda \in  
\pcf \frb\rangle $  be as in [Sh371, {\S}1]. 
Let 
$$
{\cal T}^0_{i} = \left\lbrace  \Max_{\ell =1,n}
f^\frb_{\lambda_{\ell} , \alpha_{\ell} }\upharpoonright \frb_i\colon
\lambda_{\ell}\in \pcf (\frb_{i}), \
\alpha_{\ell} < \lambda_{\ell} ,\ n < \omega \right\rbrace .
$$ 
Let  ${\cal T}_{i} = \{f \in  T^0_{i}\colon$ for every
$j < i,f\upharpoonright
\frb_{j} \in  {\cal T}^0_{j}$ moreover for some $f^\prime\in
\prod_{j<\kappa}\lambda_j$, for every $j$, $f^\prime\upharpoonright
j\in \calT_i^0$ and $f\subseteq f^\prime\}$,  and ${\cal T} =
 \bigcup_{{i} < {\kappa}} {\cal T}_{i}$,  clearly it is a tree, 
${\cal T}_{i}$ its $i$th level (or 
empty),  $|{\cal T}_{i}| \leq  \lambda ^\ast_{i}$. 
By [Sh371, 1.3, 1.4] for every  
$g \in  \prod \frb$  for some  $f \in  \prod \frb$, 
$ \bigwedge_{{i} < {\kappa}} f\upharpoonright \frb_{i} \in 
{\cal T}^0_{i}$ hence  $ \bigwedge_{{i} < {\kappa}}
f\upharpoonright \frb_{i} \in  {\cal T}_{i}$. 
So  $|{\cal T}_{i}| = \lambda ^\ast_{i}$,  
and ${\cal T}$ has $\geq  \lambda \ \ \kappa $-branches. 
By the observation
below we can finish (apply it essentially to  $F = \{\eta\colon$
for some $f\in\prod\frb$ for $i<\kappa$ we have $\eta(i)= f
\upharpoonright \frb_i$ and
for every $i < \kappa $, $f\upharpoonright\frb_{i} \in 
{\cal T}^0_{i}\})$,  then
find $A \subseteq  \kappa $, $\kappa\hefresh A\in J$ and $g^*\in
\prod_{i< \kappa}(\lambda_i +1)$ such that
 $Y' =: \{f \in  F\colon f \upharpoonright A < 
g^\ast \upharpoonright A\}$  has cardinality $\lambda$ and then the
tree will be  ${\cal T}'$ where ${\cal T}'_{i} =:
\{f\upharpoonright \frb_{i}\colon f \in  Y'\}$
and ${\cal T}' = \bigcup_{{i} < {\kappa}} {\cal T}'_{i}.$
(So actually this proves that if we have such a tree with
$\ge\theta$ ($\cf(\theta)>2^\kappa$) $\kappa$-branches then there is one with
exactly $\theta$ $\kappa$-branches.)

\proclaim{6.6A Observation} 
{\rm (1)} If  $F \subseteq \prod_{{i} < {\kappa}} \lambda_{i}$,  $J$ 
an $\aleph_{1}$-complete ideal on  $\kappa $,  and 
$[f \neq  g \in  F \Rightarrow  f \neq_{J} g]$  and $|F| \geq  \theta $,
$\cf \theta > 2^\kappa $, then for some  $g^\ast  \in \prod_{{i}
< {\kappa}} (\lambda_{i} + 1)$  we have:\  

\item{(a)} $Y=\{f\in F\colon f<_{J} g^\ast\}$ has cardinality $\theta$,

\item{(b)} for $f'<_{J}g^\ast$, we have $|\{f\in F\colon f\leq_{J}f'\}|<
\theta$,

\item{(c)}
there\footnote*{Or strightening clause (i) see the proof of 6.6B}
are  $f_{\alpha}  \in  Y$  for  $\alpha < \theta $
such that:  
$f_{\alpha}  <_{J} g^\ast $,  $[\alpha < \beta < \theta \Rightarrow  
\neg f_{\beta}  <_{J} f_{\alpha} ].$
\endproclaim

\demo{Proof} 
Let  $Z =: \left\lbrace g\colon g \in \prod_{{i} < {\kappa}}
(\lambda_{i} + 1)\right. $  and $Y_{g} =: \{f \in  F\colon  f \leq_{J}
g\}$ has cardinality  $\geq  \theta \left. \right\rbrace .$
Clearly $\langle\lambda_{i}\colon i <\kappa\rangle\in Z$ so there is 
\ $g^\ast  \in  Z$  such that:  $[g' \in  Z \Rightarrow  \neg g' <_{J}
g^\ast ]$; so (b) holds. 
Let  $Y = \{f \in  F\colon f <_{J} g^\ast \}$,  easily $Y\subseteq
 Y_{g^\ast}$ and $|Y_{g^\ast}\hefresh Y|\le 2^\kappa$ hence $|Y| \geq  
\theta $,  also clearly  $[f_{1} \neq  f_{2} \in  F \und f_{1}
\leq_J f_{2} \Rightarrow  f_{1} <_{J} f_{2}]$;
if (a) fails, necessarily (by (b)) $|Y| > \theta $. 
For each  $f \in  Y$  let  $Y_{f} = \{h \in  Y\colon h \leq_{D} f\}$,  so  
$|Y_{f}| < \theta $  hence by the Hajnal free subset theorem for some 
$Z'
\subseteq Z$,  $|Z'| = \lambda ^+$,  and $f_{1} \neq  f_{2} \in  Z'
\Rightarrow  f_{1} \notin Y_{f_{2}}$ so  $[f_{1} \neq  f_{2} \in  Z'
\Rightarrow  \neg f_{1} <_{J} f_{2}]$. 
But there is no such  $Z'$  of cardinality  $> 2^\kappa$ ([Sh111,
2.2, p.\ 264])
so (a) holds. 
As for (c): choose  $f_{\alpha}  \in  F$  by induction on  $\alpha $, 
such that  $f_{\alpha}  \in  Y\setminus \bigcup_{{\beta} <
{\alpha}} Y_{f_{\beta} }$;  it exists by cardinality considerations and 
$\langle f_{\alpha} \colon \alpha < \theta \rangle $  is as required (in
(c)).    
\sqed{6.6{\rm A}} \sqed{6.6}
 
\proclaim{6.6B Observation} 
Let  $\kappa < \lambda $  be regular uncountable,  
$2^\kappa < \mu_{i} < \lambda $ (for $i<\kappa$), 
$\mu_{i}$ increasing in  $i$. 
The following are equivalent:

\item{(A)} there is $F\subseteq\ ^\kappa\!\lambda$  such that:    
\itemitem{(i)}  $|F| = \lambda$,
\itemitem{(ii)} $|\{f\upharpoonright i\colon f \in  F\}| \leq  \mu_{i}$,
\itemitem{(iii)} $\left[ f \neq  g \in  F \Rightarrow  f
\neq_{J^{bd}_{\kappa}} g\right] $;
\item{(B)} there be a sequence  $\langle \lambda_{i}\colon i < \kappa
\rangle $  such that:
\itemitem{(i)} $2^\kappa < \lambda_{i} = \cf(\lambda_{i})\leq  \mu_{i}$,
\itemitem{(ii)} $\max\pcf\{\lambda_{i}\colon i<\kappa\}=\lambda$,
\itemitem{(iii)} for  $j < \kappa $,  $\mu_{j} \geq  \max  \pcf
\{\lambda_{i}\colon i < j\}$; 
\item{(C)} there is an increasing sequence  $\langle \fra_{i}\colon
i < \kappa
\rangle $  such 
that  $\lambda \in  \pcf \bigcup_{{i} < {\kappa}} \fra_{i}$, 
$\pcf \fra_{i} \subseteq  \mu_{i}$  $($so  $\Min(\bigcup_{{i} <
{\kappa}} \fra_{i})> | \bigcup_{{i} < {\kappa}} \fra_{i}|).$
\endproclaim

\demo{Proof} 

\subdemo{(B)$\Rightarrow$(A)}
By [Sh355, 3.4].

\subdemo{(A)$\Rightarrow$(B)}
If  $(\forall \theta )[\theta \geq  2^\kappa \Rightarrow  \theta ^\kappa
\leq  \theta ^+]$  we can directly prove (B) if for a club of
$i<\kappa$, $\mu_i> \bigcup_{j<i}\mu_j$, and contradict (A) if
this fails. 
Otherwise every normal filter 
$D$  on  $\kappa $  is\ nice (see [Sh386, {\S}1]). 
Let  $F$  exemplify (A).

Let  $K = \left\lbrace (D,g)\colon D\right. $  a\ normal filter on 
$\kappa $,
$g \in\ ^\kappa\!(\lambda +1)$, $\lambda = |\{f \in  F\colon f <_{D} 
g\}|\left. \right\rbrace .$ 
Clearly  $K$  is not empty (let  $g$  be constantly  $\lambda )$  so
by [Sh386] we can find $(D,g) \in  K$  such that: 
\item{$(\ast )_{1}$}  if  $A \subseteq  \kappa $,  $A \neq  \emptyset
\mod D$,  $g_{1}
<_{D+A} g$  then  $\lambda > |\{f \in  F\colon f <_{D+A} g_{1}\}|.$

\noindent
Let $F^\ast =\{f\in F\colon f<_{D} g\}$, so (as in the proof of 6.6) 
$|F^\ast | = \lambda .$

We claim: 
\item{$(\ast )_{2}$}  if  $h \in  F^\ast $ then  $\{f \in  F^\ast
\colon \neg h \leq_{D} f\}$  has cardinality  $< \lambda .$

\noindent
[Why?  Otherwise for some  $h \in  F^\ast $,  $F' =: \{f \in  F^\ast
\colon \neg h 
\leq_{D} f\}$  has cardinality  $\lambda $,  for  $A \subseteq  \kappa$
let $F'_{A} = \{f \in  F^\ast \colon f\upharpoonright A \leq  h
\upharpoonright A\}$  so 
$F' = \bigcup \{F'_{A}\colon A 
\subseteq  \kappa ,A \neq  \emptyset $\ mod\ $D\}$,  hence for some  $A 
\subseteq  \kappa $,  $A \neq  \emptyset $\ mod\ $D$  and $|F'_{A}| =
\lambda $; now  $(D + A,h)$  contradicts $(\ast )_{1}].$ 

By $(\ast )_{2}$ we can choose by induction on  $\alpha < \lambda $, 
a function $f_{\alpha}  \in  F^\ast $ such that  $
\bigwedge_{{\beta} < {\alpha}} f_{\beta}  <_{D} f_{\alpha} $.
By [Sh355, 1.2A(3)]  
$\langle f_{\alpha} \colon \alpha < \lambda \rangle $  has an e.u.b.
$f^\ast $.  
Let  $\lambda_{i} = \cf (f^\ast (i))$, \ clearly  $\{i < \kappa
\colon \lambda_{i} \leq  2^\kappa \} = \emptyset $\ mod\ $D$,  so
without
loss of generality  $ \bigwedge_{{i} < {\kappa}} \cf (f^\ast (i))
> 2^\kappa $ so  $\lambda_{i}$ is regular  $\in  
(2^\kappa ,\lambda ]$,  and $\lambda = \tcf\left( \prod_{{i} <
{\kappa}} \lambda_{i}/D\right) $. 
Let  $J_{i} = \left\lbrace A \subseteq  i\colon \max  
\pcf \{\lambda_{j}\colon j < i\} \leq  \mu_{i}\right\rbrace $; 
so (remembering (ii) of (A)) we can find $h_{i} \in \prod_{{j} <
{i}}f^\ast (i)$  such that: 

\item{$(\ast )_{3}$}  if  $\{j\colon j < i\} \notin  J_{i}$,  then for every 
$f \in  F$,  $f\upharpoonright i <_{J_{i}}h_{i}.$ 

Let  $h \in \prod_{{i} < {\kappa}} f^\ast (i)$  be defined\ by: 
$h(i) = \sup \left\lbrace h_{j}(i)\colon j \in  
(i,\kappa )\right. $\ \ and\ \ $\{j\colon j < i\} \notin  J_{i}\left.
\right\rbrace $. 
As  $ \bigwedge_{i} \cf [f^\ast (i)] > 2^\kappa $,  clearly 
$h < f^\ast $ hence by the choice of  
$f^\ast $ for some  $\alpha (\ast ) < \lambda $  we have:  $h <_{D} 
f_{\alpha (\ast )}$ and let  $A =: \{i < \kappa \colon h(i) < f_{\alpha
(\ast )}\}$, so  $A \in  D$. 
Define  $\lambda'_i$ as follows:  $\lambda '_{i}$  is  
$\lambda_{i}$ if  $i \in  A$,  and is  $\left( 2^\kappa \right) ^+$ if
 $i \in  \kappa \backslash A$. 
Now  $\langle \lambda '_{i}\colon i < \kappa \rangle $  is as 
required in (B). 

\subdemo{(B)$\Rightarrow$(C)}
Straightforward.

\subdemo{(C)$\Rightarrow$(B)}
By [Sh371, {\S}1].      
\sqed{6.6{\rm B}}

\proclaim{6.6C Claim} 
If  $F \subseteq\  ^\kappa\!\!\Ord $,  $2^\kappa < \theta = \cf(\theta)
\leq  |F|$  then we can find $g^\ast \in\ ^\kappa\!\!\Ord $ 
and a proper ideal $I$ on  $\kappa $  and $A \subseteq  \kappa $, 
$A \in  I$  such that: 
\item{(a)} $ \prod_{{i} < {\kappa}} g^\ast (i)/I$  has true
cofinality  $\theta$,  and for each  $i \in
\kappa\hefresh A$  we have \break
$\cf [g^\ast (i)] > 2^\kappa$,

\item{(b)} for every  $g \in\  ^\kappa\!\!\Ord $ satisfying 
$g\upharpoonright A = g^\ast \upharpoonright A$,  
$g\upharpoonright (\kappa \backslash A) < g^\ast \upharpoonright
(\kappa \backslash A)$ we can find $f \in  F$  such that:   
$f\upharpoonright A = g^\ast \upharpoonright A$,  $g
\upharpoonright (\kappa \backslash A) < 
f\upharpoonright (\kappa \backslash A) < g^\ast \upharpoonright
(\kappa \backslash A).$
\endproclaim

\demo{Proof} 
As in [Sh410, 3.7 proof of (A)$\Rightarrow$(B)].
(In short let  
$f_{\alpha}  \in  F$  for  $\alpha < \theta $  be distinct,  $\chi $ 
large enough,  $\langle N_{i}\colon i < \left( 2^\kappa \right)%
^+\rangle$
as there,  $\delta_{i} =: \sup (\theta \cap  N_{i})$,  $g_{i}
\in$\break ${}^\kappa\!\!\Ord $,  
$g_{i}(\zeta ) =: \Min \left[ N \cap  \Ord \backslash f_{\delta_{i}}
(\zeta )\right] $,
$A \subseteq  \kappa $  and $S \subseteq  \lbrace i < 
\left( 2^\kappa \right) ^+\colon \cf (i) = \kappa ^+\rbrace $ 
stationary,  $[i \in  S \Rightarrow  g_{i} = g^\ast ]$, 
$\left[ \zeta < \alpha \und i \in  S 
\Rightarrow  [f_{\delta_{i}}(\zeta ) = g^\ast (\zeta ) \equiv  \zeta
\in  A]\right] $  and for some  $i(\ast ) < \left( 2^\kappa
\right)^+$,  $g^\ast  
\in N_{i(\ast)}$, so $[\zeta \in\kappa\hefresh  A \Rightarrow \cf
g^\ast (\zeta ) > 2^\kappa ].)$ 
\sqed{6.6{\rm C}}

\proclaim{6.6D Claim} 
Suppose  $D$  is a filter on  $\theta = \cf(\theta)$,  
$\sigma $-complete,  $\theta > |\alpha|^\kappa$ for
$\alpha<\sigma$,  and for each $\alpha <
\theta $,  $\bar \beta = \langle \beta ^\alpha_{\epsilon} \colon
\epsilon  < 
\kappa \rangle $  is a sequence of ordinals. 
Then for every  $X \subseteq  
\theta $,  $X \neq  \emptyset $\ mod\ $D$  there is  
$\langle \beta ^\ast_{\epsilon} \colon \epsilon  < \kappa \rangle$ (a
sequence of ordinals) and $w \subseteq  \kappa $  such that:
\item{(a)} $\epsilon  \in  \kappa \backslash w \Rightarrow  \sigma \leq 
\cf(\beta^\ast_{\epsilon}) \leq  \theta$, 

\item{(b)} if  $\beta '_{\epsilon}  \leq  \beta ^\ast_{\epsilon} $ and
$[\epsilon \in  w \equiv  \beta '_{\epsilon}  = \beta ^\ast_{\epsilon} ]$, then
$\left\lbrace \alpha \in  X\colon\right.$ for every  $\epsilon
 < \kappa $ we have
$\beta '_{\epsilon}  \leq  \beta ^\alpha_{\epsilon}  \leq 
\beta ^\ast_{\epsilon} $ and 
$[\epsilon \in  w \equiv  \beta ^\alpha_{\epsilon}  = 
\beta ^\ast_{\epsilon} ]\left. \right\rbrace  \neq  \emptyset \mod D$.
\endproclaim

\demo{Proof} 
Essentially by the same proof as 6.6C (replacing  $\delta_{i}$ by  
$\Min \{\alpha \in  X\colon$ for every $Y \in  N_{i} \cap  D$ we have 
$\alpha \in  Y\})$. See more [Sh513, \S6].
\sqed{6.6{\rm D}}

\rem{6.6E Remark}
We can rephrase the conclusion as:
\item{(a)} $B =: \{\alpha\in X\colon$ if $\epsilon\in w$ then
$\beta^\alpha_\epsilon = \beta^\ast_\epsilon$, and: if
$\epsilon\in\kappa\hefresh w$ then $\beta^\alpha_\epsilon$ is
$<\beta^\ast_\epsilon$ but $>\sup\{\beta^\ast_\zeta\colon\zeta
<\epsilon, \beta^\alpha_\zeta<\beta^\ast_\epsilon\}\}$ is $\not=
\emptyset \mod D$.

\item{(b)} If $\beta'_\epsilon < \beta_\epsilon$ for $\epsilon
\in\kappa\hefresh w$ then $\{\alpha\in B\colon$ if
$\epsilon\in\kappa\hefresh w$ then $\beta_\epsilon^\alpha
>\beta'_\epsilon\} \not= \emptyset\mod D$.

\item{(c)} $\epsilon\in\kappa\hefresh w\Rightarrow
\cf(\beta'_\epsilon)$ is $\le\theta$ but $\ge \sigma$.

\rem{6.6F Remark}
(1) If $|\fra| <\min(\fra)$, $F\subseteq\Pi\fra$,
$|F|=\theta=\cf\theta\not\in \pcf(\fra)$ and even $\theta>\sigma
=\sup(\theta^+\cap \pcf(\fra))$ then for some $g\in\Pi\fra$, the
set $\{f\in F\colon f<g\}$ is unbounded in $\theta$ (or use a
$\sigma$-complete $D$ as in 6.6E).
(This is as $\Pi\fra / J_{<\theta}[\fra]$ is $\min(\pcf(\fra)
\hefresh\theta)$-directed as the ideal $J_{<\theta}[\fra]$ is
generated by $\le\sigma$ sets; this is discussed in [Sh513, \S6].)

\rem{6.6G Remark}
It is useful to note that 6.6D is useful to use [Sh462, \S4,
5.14]: e.g.\
for if $n<\omega$, $\theta_0 < \theta_1 < \cdots < \theta_n$, satisfying
$(*)$ below,
for any $\beta'_\epsilon\le \beta^\ast_\epsilon$ satisfying
$[\epsilon\in w \equiv \beta'_\epsilon < \beta^\ast_\epsilon]$
we can find $\alpha<\gamma$ in $X$ such that:
$$
i\in w \equiv \beta^\alpha_\epsilon = \beta^\ast_\epsilon,
$$
$$
\{\epsilon,\zeta\}\subseteq \kappa\hefresh w\ \& \
\{\cf(\beta^\ast_\epsilon), \cf(\beta^\ast_\zeta)\} \subseteq
[\theta_l, \theta_{l+1}))\ \&\ l\hbox{ even }\Rightarrow
\beta^\alpha_\epsilon < \beta^\gamma_\zeta,
$$
$$
\{\epsilon,\zeta\}\subseteq\kappa\hefresh w\ \&\
\{\cf(\beta^\ast_\zeta), \cf(\beta^\ast_\zeta)\}\subseteq
[\theta_l, \theta_{l+1})\  \&\ l \hbox{ odd }\Rightarrow
\beta_\epsilon^\gamma <\beta_\zeta^\alpha
$$
where
\item{$(*)$} (a) $\epsilon \in\kappa\hefresh w\Rightarrow \cf
(\beta^\ast_\epsilon)\in [\theta_0, \theta_n)$, and

\item{\quad} (b) $\max\pcf[\{\cf(\beta^\ast_\epsilon)\colon \epsilon\in\kappa
\hefresh w\}\cap \theta_l] \le \theta_l$ (which holds if $\theta_l
= \sigma^+_l$, $\sigma^\kappa_l = \sigma_l$ for $l\in\{1,\ldots,
n\}$).

\proclaim{6.7 Claim} 
For any  $\fra$, $|\fra|<\Min(\fra)$, we can find
$\bar\frb=\langle\frb_{\lambda}\colon\lambda \in\fra\rangle $ 
such that:  
\item{$(\alpha )$} $\bar \frb$  is a generating sequence, i.e.\      
$$
\lambda \in  \fra \Rightarrow  J_{\leq \lambda }[\fra] = 
J_{<\lambda }[\fra] + \frb_{\lambda},  
$$

\item{$(\beta)$} $\bar\frb$ is smooth, i.e.\ for $\theta<\lambda$ 
in $\fra$,
$$
\theta\in\frb_{\lambda}\Rightarrow\frb_{\theta}\subseteq \frb_{\lambda},
$$ 

\item{$(\gamma)$} $\bar\frb$ is closed, i.e.\ for $\lambda\in 
\pcf(\fra)$ we have $\frb_{\lambda} =\fra\cap\pcf(\frb_{\lambda})$.
\endproclaim

\demo{Proof} 
Let $\langle\frb_{\theta}[\fra]\colon\theta\in\pcf\fra\rangle$
be as in 
[Sh371, 2.6].  For  $\lambda \in  \fra$,  let  $\bar f^{\fra,\lambda } =
\langle f^{\fra,\lambda }_{\alpha} \colon \alpha < \fra\rangle $  be a 
$<_{J_{\lambda} [\fra]}$-increasing cofinal sequence of members of 
$\prod \fra$,  satisfying: 
\item{$(\ast )_{1}$}  if  $\delta < \lambda $,  $|\fra| < \cf(\delta)<
\Min\fra$ and $\theta \in  \fra$  then:
$$
f^{\fra,\lambda }_{\delta} (\theta ) = \Min \left\lbrace
\bigcup_{{\alpha} \in  {C}}f^{\fra,\lambda }_{\alpha} (\theta )
\colon C \hbox{ a club of } \delta\right\rbrace
$$

\noindent
[exists by [Sh345a, Def. 3.3(2)$^{\rm b}$ + Fact 3.4(1)]]. 

Let  $\chi  = \beth
_{\omega} (\sup  \fra)^+$,  $|\fra| < \kappa = \cf \kappa < \Min\fra$
(without loss of 
generality there is such  $\kappa )$  and $\bar N = \langle N_{i}\colon
i < 
\kappa \rangle $  be an increasing continuous sequence of elementary
submodels of  $(H(\chi ),\in ,<^\ast_{\chi} )$,  $N_{i} \cap  \kappa $
an ordinal,  
$\bar N\upharpoonright (i + 1) \in  N_{i+1}$,  $\|N_{i}\| < \kappa $, 
and $\fra$,  $\langle \bar f^{\fra,\lambda }\colon \lambda \in 
\fra\rangle $  belong to $N_{0}$.  Let  
$N_{\kappa}  = \bigcup_{{i} < {\kappa}} N_{i}$. 
For every  $\lambda \in  \fra$,  for some club $E_{\lambda} $ of  
$\kappa ,$ 
\item{$(\ast )$}  $\theta \in  \fra \Rightarrow 
f^{\fra,\lambda}_{\sup (N_{\kappa} \cap \lambda )}(\theta ) =
\bigcup_{{\alpha}\in{E_{\lambda}}}f^{\fra,\lambda}_{\sup
(N_\alpha\cap\lambda)} (\theta).$

\noindent
Let $E = \bigcap_{{\lambda} \in  {\fra}}E_{\lambda} $,  so  $E$ 
is a club of  $\kappa $.  For any  $i < j < \kappa $  
let 
$$
\frb^{i,j}_{\lambda} =\left\lbrace\theta\in \fra\colon \sup (N_{i} \cap 
\theta ) < f^{\fra,\lambda }_{\sup (N_{j}\cap \lambda )}(\theta )
\right\rbrace .
$$ 
As in the proof of [Sh371, 1.3], possibly shrinking $E$, we have:
\item{$(\ast )_{2}$} for  $i < j$  from\footnote*{Actually for
 any $i<j<\kappa$ clauses $(\beta)$, $(\gamma)$,
$(\delta)$ hold.} $E$ 
and $\lambda\in\fra$,  we have:
\itemitem{$(\alpha )$} $J_{\leq \lambda }[\fra] = J_{<\lambda }[\fra] +
\frb^{i,j}_{\lambda} $  
(hence  $\frb^{i,j}_{\lambda}  = \frb_{\lambda} [\fra]\mod
J_{<\lambda }[\fra])$,

\itemitem{$(\beta )$} $\frb^{i,j}_{\lambda} \subseteq \lambda ^+\cap \fra$,\  

\itemitem{$(\gamma )$} $\langle \frb^{i,j}_{\lambda} \colon \lambda \in \fra
\rangle  \in  N_{j+1}$,    

\itemitem{$(\delta )$} 
$f^{\fra,\lambda }_{\sup (N_{\kappa} \cap \lambda )}\upharpoonright
\frb^{i,j}_{\lambda}  = 
\langle (\theta ,\sup (N_{\kappa}  \cap  \theta ))\colon \theta \in  
\frb^{i,j}_{\lambda} \rangle $,

\itemitem{$(\epsilon )$} $f^{\fra,\lambda }_{\sup (N_{\kappa} \cap
\lambda )} \leq  \langle (\theta ,\sup (N_{\kappa}  \cap  \theta ))\colon
\theta \in  \fra\rangle .$

\noindent
We now define by\ induction on  $\epsilon  < |\fra|^+$,  for  $\lambda
\in  \fra$ (and $i<j<\kappa$),  the set 
$\frb^{i,j,\epsilon }_{\lambda} :$  
$$
\eqalign{
\frb^{i,j,0}_{\lambda} & = \frb^{i,j}_{\lambda} \cr
\frb^{i,j,\epsilon +1}_{j} & = \frb^{i,j,\epsilon }_{\lambda}  \cup  
\bigcup \left\lbrace \frb^{i,j,\epsilon }_{\theta} \colon \theta \in  
\frb^{i,j,\epsilon }_{\lambda} \right\rbrace  \cup  \left\lbrace \theta
\in  \fra\colon \theta \in  \pcf \frb^{i,j,\epsilon }\right\rbrace,
\cr
\frb^{i,j,\epsilon }_{\lambda}  & = \bigcup_{{\zeta} < {\epsilon}}
\frb^{i,j,\zeta }_{\lambda} \ \hbox{ for } \epsilon  <
|\fra|^+\hbox{ limit.} \cr}
$$
Clearly for $\lambda\in \fra$,  $\langle \frb^{i,j,\epsilon }_{\lambda}
\colon \epsilon  < |\fra|^+\rangle $  belongs to  $N_{j+1}$ and is a
non-decreasing sequence of subsets of $\fra$,  hence for some 
$\epsilon (i,j,\lambda ) < |\fra|^+$,  
$$
\left[ \epsilon  \in  (\epsilon(i,j,\lambda) ,|\fra|^+) \Rightarrow  
\frb^{i,j,\epsilon }_{\lambda}  =
\frb^{i,j,\epsilon (i,j,\lambda )}_{\lambda} \right] .
$$
So letting  $\epsilon (i,j) =  \sup_{\lambda \in \fra}
\epsilon (i,j,\lambda ) < |\fra|^+$ we have: 
\item{$(\ast )_{3}$}  $\epsilon (i,j) \leq  \epsilon  < |\fra|^+
\Rightarrow \bigwedge_{{\lambda} \in  {\fra}}\frb^{i,j,\epsilon
(i,j)}_{\lambda}  = \frb^{i,j,\epsilon }_{\lambda} .$ 

Which of the properties required from \ $\langle \frb_{\lambda} \colon
\lambda \in  \fra\rangle $  are satisfied by  $\langle \frb^{i,j,
\epsilon (i,j)}_{\lambda} \colon\shvor \lambda 
\in  \fra\rangle ?$  Note $(\beta )$, $(\gamma )$ hold by the inductive 
definition of  $\frb^{i,j,\epsilon }_{\lambda} $ (and the choice of  
$\epsilon (i,j))$,  as for property $(\alpha )$, one half,  
$J_{\leq \lambda }[\fra] \subseteq  J_{<\lambda }[\fra] + 
\frb^{i,j,\epsilon (i,j)}_{\lambda} $ hold by $(\ast )_{2}(\alpha )$ (and 
$\frb^{i,j}_{\lambda}  = \frb^{i,j,0}_{\lambda}  \subseteq  
\frb^{i,j,\epsilon (i,j)}_{\lambda} )$,  so it is enough to prove (for 
$\lambda \in  \fra):$ 
\item{$(\ast )_{4}$}  $\frb^{i,j,\epsilon (i,j)}_{\lambda}  \in  J_{\leq
\lambda }[\fra].$ 

For this end we define by induction on $\epsilon<|\fra|^+$ functions  
$f^{\fra,\lambda ,\epsilon }_{\alpha} $ with domain 
$\frb^{i,j,\epsilon }_{\lambda} $ 
for every  $\alpha < \lambda \in  \fra$,  such that  $\zeta < \epsilon  
\Rightarrow  f^{\fra,\lambda ,\zeta }_{\alpha}  \subseteq  
f^{\fra,\lambda ,\epsilon }_{\alpha}$, so the domain increases
with $\epsilon$.

We let  $f^{\fra,\lambda ,0}_{\alpha}  = 
f^{\fra,\lambda }_{\alpha} \upharpoonright \frb^{i,j}_{\lambda} ,$
$f^{\fra,\lambda ,\zeta }_{\alpha}  = \bigcup_{{\zeta} <
{\epsilon}} f^{\fra,\lambda ,\zeta }_{\alpha} $ for  $\epsilon <
|\fra|^+$ limit, and
$f^{\fra,\lambda ,\epsilon +1}_{\alpha} $ is defined by defining each  
$f^{\fra,\lambda ,\epsilon +1}_{\alpha} (\theta )$  as follows:

\subdemo{Case 1} 
If  $\theta \in  \frb^{i,j,\epsilon }_{\lambda} $ then  
$f^{\fra,\lambda ,\epsilon }_{\alpha} (\theta ).$

\subdemo{Case 2} 
If\  $\mu  \in  \frb^{i,j,\epsilon }_{\lambda} $,  $\theta \in  
\frb^{i,j,\epsilon }_{\mu} $ and not Case 1 and $\mu $  minimal under
those conditions, then  $f^{\fra,\mu ,\epsilon }_{\beta} (\theta )$ 
where we choose  $\beta = f^{\fra,\lambda ,\epsilon }_{\alpha} (\mu ).$

\subdemo{Case 3} 
If  $\theta \in  \fra \cap  \pcf (\frb^{i,j,\epsilon }_{\lambda} )$ 
and not Case 1 or 2, then 
$$
\Min \left\lbrace \gamma < 
\theta \colon f^{\fra,\lambda ,\epsilon }_{\alpha} \upharpoonright
\frb_{\theta} [\fra] 
\leq_{J_{<\theta }[\fra]}f^{\fra,\theta ,\epsilon }_{\gamma}
\right\rbrace .
$$

Now $\langle\langle\frb^{i,j,\epsilon}_{\lambda}\colon\lambda\in
\fra\rangle \colon \epsilon  < 
|\fra|^+\rangle $  can be computed from $\fra$ and 
$\langle \frb^{i,j}_{\lambda} \colon \lambda \in  \fra\rangle $. 
But the latter belong\footnote*{As $\langle \frb^{i, j,
\epsilon}_\lambda: \lambda\in\fra\rangle: \epsilon |\fra|^+\rangle$ is
eventually constant, also each member of the sequence belongs to
$N_{j+1}$.}
to
$N_{j+1}$, so the former belongs to $N_{j+1}$, so as also\break
$\langle \langle f^{\fra,\lambda }_{\alpha} \colon \alpha < 
\lambda \rangle \colon \lambda \in  \pcf \fra\rangle $  belongs to 
$N_{j+1}$ we clearly get 
that
$$
\left\langle\langle \langle f^{\fra,\lambda ,\epsilon }_{\alpha}
\colon \epsilon  < |\fra|^+\rangle\colon \alpha < \lambda \rangle \colon
\lambda \in \fra\right\rangle
$$
belongs to  $N_{j+1}$.
Next we prove by induction on  $\epsilon $  that, for  $\lambda\in\fra$,
we have: 
$$
\theta \in  \frb^{i,j,\epsilon }_{\lambda}  \und \lambda \in  \fra 
\Rightarrow 
f^{\fra,\lambda ,\epsilon }_{\scriptstyle{\sup(N_\kappa\cap\theta)}}
(\theta ) = \sup (N_{\kappa}  \cap  \theta).
\leqno{\otimes_{1}} 
$$ 

For  $\epsilon  = 0$\  this is by $(\ast )_{2}(\delta )$. 
For  $\epsilon $  
limit, by the induction hypothesis and the definition of
$f^{\fra,\lambda ,\epsilon }_{\alpha} $.  For  $\epsilon  + 1$,  we check  
$f^{\fra,\lambda ,\epsilon +1}_{\sup (N_{\kappa} \cap \lambda )}(\theta )$ 
according to the case in its definition; for Case 1 use the induction 
hypothesis applied to  
$f^{\fra,\lambda ,\epsilon }_{\sup (N_{\kappa} \cap \lambda )}$. 
For Case 2 (with  $\mu )$,  by the induction hypothesis applied to  
$f^{\fra,\mu ,\epsilon }_{\sup (N_{\kappa} \cap \mu )}$. 
Lastly, for Case 3 (with  
$\theta )$  we should note:

\item{(i)} $\frb^{i,j,\epsilon }_{\lambda}  \cap  \frb_{\theta} [\fra]
\notin J_{<\theta }[\fra]$  
(by the case's assumption and $(\ast )_{2}(\alpha )$ above),

\item{(ii)} $f^{\fra,\lambda ,\epsilon }_{\sup (N_{\kappa} \cap
\lambda )}
\upharpoonright (\frb^{i,j,
\epsilon }_{\lambda}  \cap  \frb^{i,j,\epsilon }_{\theta} ) \subseteq  
f^{\fra,\theta ,\epsilon }_{\sup (N_{\kappa} \cap \theta )}$
(by the induction hypothesis for $\epsilon$, used concerning
$\lambda$ and $\theta$)
hence (by the definition in case 3 and (i) + (ii)),

\item{(iii)} $f^{\fra,\lambda ,\epsilon +1}_{\sup (N_{\kappa} \cap
\lambda )}(\theta ) \leq  \sup (N_{\kappa}  \cap  \theta )$.

\noindent
Now if  $\gamma < \sup (N_{\kappa}  \cap  \theta )$  then for some  
$\gamma (1)$,  $\gamma < \gamma (1) \in  N_{\kappa}  \cap  \theta $,
so letting $\frb =: \frb_\lambda^{i,j,\epsilon}\cap \frb_\theta
[\fra] \cap\frb_\theta^{i,j,\epsilon}$, it belongs to
$J_{\le\theta}[\fra] \hefresh J_{<\theta} [\fra]$, we have
$$
f^{\fra,\theta }_{\gamma} \upharpoonright \frb
<_{J_{<\theta }[\fra]} 
f^{\fra,\theta }_{\gamma (1)}\upharpoonright \frb  \leq  
f^{\fra,\theta ,\epsilon}_{\sup (N_{\kappa} \cap \theta )}$$ hence 
$f_{\sup(N_\kappa\cap\lambda)}^{\fra,\lambda,\epsilon+1}(\theta)
>\gamma$; as this holds for
every $\gamma<\sup(N_\kappa\cap\theta)$ we have obtained

\item{(iv)} $f^{\fra,\lambda ,\epsilon +1}_{\sup (N_{\kappa} \cap \lambda)}
(\theta ) \geq \sup (N_{\kappa} \cap \theta )$;

\noindent
together we have finished proving
the inductive step for $\epsilon+1$, hence we have proved
$\otimes_{1}$.

This is enough for proving $\frb_\lambda^{i,j,\epsilon}\in
J_{\le\lambda}[\fra]$: Why?
If it fails, as $\frb_\lambda^{i,j,\epsilon}\in N_{j+1}$ and
$\langle
f_\alpha^{\fra,\lambda,\epsilon}\colon\alpha<\lambda\rangle$
belongs to $N_{j+1}$, there is
$g\in\prod\frb_\lambda^{i,j,\epsilon}$ s.t.
$$
\alpha<\lambda \Rightarrow f_\alpha^{\fra,\lambda,\epsilon}
\upharpoonright \frb^{i,j,\epsilon} <
g\mod J_{\le\lambda}[\fra].
\leqno(*)
$$
Wlog $g\in N_{j+1}$; by $(*)$,
$f^{\fra,\lambda,\epsilon}_{\sup(N_\kappa\cap\lambda)} <g \mod
J_{\le\lambda}[\fra]$.
But
$g<\langle\sup(N_\kappa\cap\theta)\colon\theta\in
\frb_\lambda^{i,j,\epsilon}\rangle$.
Together this contradicts $\oplus_1$!

This ends the proof of 6.7.
\sqed{6.7}

\proclaim{6.7A Claim} 
Assume  $|\fra| < \kappa = \cf(\kappa)< \Min(\fra)$,  $\sigma $  an 
infinite ordinal,  $|\sigma |^+ < \kappa $.  Let  $\bar f$, $\bar N = 
\langle N_{i}\colon i < \kappa \rangle $,  $N_{\kappa} $ be as in the
proof of {\rm 6.7}. 
Then we can find $\bar i =\langle i_{\alpha}\colon\alpha\leq\sigma
\rangle$, $\bar\fra =\langle\fra_{\alpha}\colon\alpha< \sigma\rangle$ 
and $\langle \langle \frb^\beta_{\lambda} [\bar \fra]\colon \lambda \in 
\fra_{\beta} \rangle \colon \beta <
\sigma \rangle $  such that:
\item{(a)} $\bar i$ is  a strictly increasing continuous sequence of
ordinals  $< \kappa $,

\item{(b)} for  $\beta < \sigma $  we have  $\langle i_{\alpha}\colon
\alpha \leq  \beta \rangle  \in  N_{i_{\beta +1}}$
(hence\footnote*{We can get $\bar i\upharpoonright (\beta+1)\in
N_{i_{\beta}+1}$ if $\kappa$ succesor of regular and $\bar C$ a
square later.}
$\langle N_{i_{\alpha} }\colon \alpha \leq  \beta \rangle  \in 
N_{i_{\beta +1}})$  and $\langle \frb^\gamma_{\lambda} [\bar \fra]
\colon \lambda \in  \fra_{\gamma} $\ and
$\gamma \leq  \beta \rangle  \in  N_{i_{\beta +1}}$, 

\item{(c)} $\fra_{\beta}  = N_{i_{\beta} } \cap  \pcf(\fra)$,  so 
$\fra_{\beta} $ is increasing 
continuous in $\beta$, $\fra\subseteq\fra_{\beta}\subseteq\pcf\fra$,
$|\fra_{\beta} | < \kappa$, 

\item{(d)} $\frb^\beta_{\lambda} [\bar \fra] \subseteq  \fra_{\beta} $ (for 
$\lambda \in  \fra_{\beta} )$,

\item{(e)} $J_{\leq\lambda}[\fra_{\beta}] = J_{<\lambda }[\fra_{\beta}]
+\frb^\beta_{\lambda} [\fra]$  (so  $\lambda \in  \frb_{\lambda} [\fra]$
and $\frb_{\lambda} [\fra] \subseteq  \lambda ^+)$,  

\item{(f)} if $\mu<\lambda$ are in $\fra_{\beta}$ and $\mu\in  
\frb^\beta_{\lambda} [\bar\fra]$ then $\frb^\beta_{\mu}[\bar\fra]
\subseteq  \frb^\beta_{\lambda} [\bar \fra]$ (i.e.\
smoothness),

\item{(g)} $\frb^\beta_\lambda [\bar \fra] = \fra_{\beta}  \cap  \pcf
\frb^\beta_{\lambda} [\bar \fra]$ (i.e.\ closedness),

\item{(h)} if $\frc\subseteq\fra_{\beta}$, $\beta<\sigma$,  $\frc\in  
N_{i_{\beta +1}}$ then for some finite $\frd\subseteq\fra_{\beta +1}
\cap \pcf(\frc)$,  we have  $\frc \subseteq \bigcup_{{\mu} \in 
{\frd}}\frb^{\beta +1}_{\mu} [\bar \fra]$;  more generally,
\footnote{**}{If in (h)$^+$,  $\theta = \aleph_{0}$, we get (h).}

{\itemindent=32pt
\item{(h)$^+$} if  $\frc \subseteq  \fra_{\beta} $,  $\beta < \sigma $,
$\frc\in N_{i_{\beta +1}}$,  $\theta = \cf(\theta)\in
N_{i_{\beta+1}}$, then for some 
$\frd\in N_{i_{\beta+1}}, \frd\subseteq \fra_{\beta +1}\cap 
\pcf_{\theta {\scriptstyle{\rm -complete}}} (\frc)$  we have
$\frc \subseteq \bigcup_{{\mu} \in  {\frd}}\frb^{\beta+1}_{\mu}
[\bar \fra]$  and $|\frd| < \theta $,

}
\item{(i)} $\frb^\beta_{\lambda} [\bar \fra]$  increases with  $\beta$.

\endproclaim

This will be proved below.

\proclaim{6.7B Claim} 
In {\rm 6.7A} we can also have:
\item{(1)}
if we let  $\frb_{\lambda} [\bar \fra] =
\frb^\sigma_{\lambda} [\fra] = \bigcup_{{\beta} < {\sigma}}
\frb^\beta_{\lambda} [\bar \fra]$,  $\fra_{\sigma}  = \bigcup_{{\beta}
< {\sigma}} \fra_{\beta} $ then also for  $\beta = \sigma $  we have
{\rm (b)} (use  $N_{i_{\beta} +1})$, {\rm (c), (d), (f), (i)}.

\item{(2)}  If  $\sigma = \cf (\sigma ) > |\fra|$  then for  $\beta =
\sigma $  also {\rm (e), (g)}.

\item{(3)}  If  $\cf (\sigma ) > |\fra|$,  $\frc \in  N_{i_{\sigma} }$, 
$\frc \subseteq  \fra_\sigma$  (hence  $|\frc| < \Min(\frc)$ and
$\frc\subseteq\fra_\sigma)$,
then for some finite  $\frd \subseteq  (\pcf \frc) \cap  
\fra_{\sigma} $ we have  $\frc \subseteq \bigcup_{{\mu} \in  {\frd}}
\frb_{\mu} [\bar\fra]$. 
Similarly for $\theta $-complete,  $\theta < \cf (\sigma )$
(i.e.\ we have clauses {\rm (h), (h)$^+$} for $\beta=\sigma$).

\item{(4)}  We can have continuity in  $\delta\le\sigma$ when 
$\cf (\delta ) > |\fra|$, i.e.\ $\frb_\lambda^\delta
=\bigcup_{\beta<\delta}\frb_\lambda^\beta$.
\endproclaim

\rem{6.7C Remark} 
\item{(1)} If we want to use length  $\kappa $,  use  $\bar N$  as 
produced in [Sh420, 2.6] so  $\sigma = \kappa .$

\item{(2)}  Concerning 6.7B, in 6.7C(1) for a club $E$  of 
$\sigma = \kappa $,  
we have $\alpha \in  E \Rightarrow  \frb^\alpha_{\lambda} [\bar \fra] = 
\frb_{\lambda} [\bar \fra] \cap  \fra_{\alpha} .$

\item{(3)}  We can also use 6.7 (6.7A, 6.7B) to give an alternative
proof of part of the localization theorems similar to the one given in
the Spring '89 lectures. For example:

\item{(3A)} If  $|\fra| < \theta = \cf \theta < \Min(\fra)$,
for no  $\lambda_{i} \in  \pcf \fra$ $(i < \theta )$ $\alpha
< \theta $,  do we have  $ \bigwedge_{{\alpha} < {\theta}}
[\lambda_{\alpha}  > \max  \pcf \{\lambda_{i}\colon i < \alpha \}].$

\item{(3B)} if $|\fra|<\Min(\fra)$, $|\frb|<\Min\frb$,
$\frb\subseteq \pcf(\fra)$, $\lambda\in\pcf(\fra)$, then for some
$\frc\subseteq\frb$ we have $|\frc|\le |\fra|$ and $\lambda\in
\pcf(\frc)$.
\endproclaim

\demo{Proof of (3A) from 6.7C(3)} 
Without loss of generality $\Min \fra > \theta ^{+3}$,  let  $\kappa =
\theta ^{+2}$,  let  $\bar N$, $N_{\kappa} $, 
$\bar \fra$, $\frb$ (as a function), 
$\langle i_{\alpha} \colon \alpha \leq \sigma =: |\fra|^+\rangle $
be as in 6.7A but also  $\langle \lambda_{i}\colon i<\theta\rangle \in 
N_{0}$. 
So for $j<\theta$, $\frc_{j} =: \{\lambda_{i}\colon i < j\} \in  N_{0}$
(and $\frc_{j} \subseteq  \fra_{0})$ hence (by clause (h) of 6.7A),
for some finite  $\frd_{j} \subseteq
 \fra_{1} \cap  \pcf \frc_{j} = N_{i_{1}} \cap  
\pcf \fra \cap  \pcf \frc_{j}$ we have  $\frc_{j} \subseteq
 \bigcup_{{\lambda} \in  {\frd_{j}}}\frb^1_{\lambda} [\bar \fra]$. 
Assume  $j(1) < j(2) < \theta $. 
Now if  $\mu  \in  \fra \cap \bigcup_{\lambda \in  \frd_{j(1)}}
\frb^1_{\lambda} [\bar \fra]$  then for some  $\mu_0\in\frd_{j(1)}$
we have $\mu  \in  \frb^1_{\mu_{0}}[\bar \fra]$; 
now  $\mu_{0} \in  \frd_{j(1)} 
\subseteq  \pcf (\frc_{j(1)}) \subseteq  \pcf (\frc_{j(2)}) \subseteq  
\pcf \left( \bigcup_{\lambda \in \frd_{j(2)}}\frb^1_{\lambda}
[\bar \fra]\right)  = \bigcup_{{\lambda} \in  {\frd_{j(2)}}}\pcf
(\frb^1_{\lambda} [\bar \fra])$  hence (by clause (g) of 6.7A as
$\mu_0\in \frd_{j(0)}\subseteq N_1$)
for some  $\mu_{1} \in  
\frd_{j(2)}$,  $\mu_{0} \in  \frb^1_{\mu_{1}}[\bar \fra]$. 
So by clause (f) of 6.7A we have  $\frb^1_{\mu_{0}}[\bar \fra] \subseteq
\frb^1_{\mu_{1}}[\bar \fra]$  
so remembering $\mu\in\frb_{\mu_0}^1[\bar\fra]$, we have 
$\mu  \in  \frb^1_{\mu_{1}}[\bar \fra]$.
Remembering $\mu$ was any member of
$\fra\cap\bigcup_{\lambda\in\frd_{j(1)}}\frb^1_\lambda[\bar\fra]$,
we have
$\fra \cap \bigcup_{{\lambda \in } \frd {j(1)}}\frb^1_{\lambda}
[\bar \fra] \subseteq  \fra \cap \bigcup_{{\lambda \in } \frd {j(2)}}
\frb^1_{\lambda} [\bar \fra]$ (holds without ``$\fra\cap$'' but
not used). 
So $\langle\fra\cap\bigcup_{{\lambda}\in{\frd_{j}}}\frb^1_{\lambda}
[\bar \fra]\colon j < \theta \rangle $ \ is a
non-decreasing sequence of subsets of $\fra$, but $\cf(\theta)>|\fra|$, 
so the sequence is eventually constant, say for  $j \geq  j(\ast )$. 
But 
$$
\eqalign{
\max  \pcf \left( \fra \cap  
 \bigcup_{{\lambda} \in  {\frd_{j}}}\frb^1_{\lambda} [\bar \fra]\right) 
& \leq  \max  \pcf \left( \bigcup_{{\lambda} \in  {\frd_{j}}}
\frb^1_{\lambda} [\bar \fra]\right) \cr
& =  \max_{\lambda \in \frd_j}
\left( \max  \pcf (\frb^1_{\lambda} [\bar \fra])\right) \cr
& =  \max_{\lambda \in \frd_j}\lambda \le \max  \pcf \{\lambda_{i}\colon
i < j\} < \lambda_{j} \cr
& = \max  \pcf \left( \fra \cap \bigcup_{{\lambda \in } \frd_{j+1}}
\frb^1_{\lambda} [\bar \fra]\right) \cr}
$$
(last equality as  $\frb_{\lambda_{j}}[\fra] \subseteq  
\frb^1_{\lambda} [\bar \fra]\mod  J_{<\lambda }[\fra_{1}])$. 
Contradiction. 
\sqed{6.7C}

\demoinfo{Proof of 6.7C(3B)}{like [Sh371, \S 3]}
Included for completeness.
If this fails choose a counterexample $(\fra, \frb,\lambda)$ with
$|\frb|$ minimal, and among those with\break
$\max\pcf(\frb)$ minimal
and among those with 
$\bigcup\{\mu^+\colon\mu\in\lambda\cap\pcf(\frb)\}$
minimal.
So $\max\pcf(\frb)=\lambda$, and
$\mu=\sup[\lambda\cap\pcf(\fra)]$ is not in $\pcf(\frb)$ or
$\mu=\lambda$.
Try to choose by induction on $i<|\fra|^+$, $\lambda_i\in\lambda
\cap\pcf(\frb)$, $\lambda_i >\max\pcf\{\lambda_j\colon j<i\}$,
by 6.7C(3A), we will be stuck at some $i$, and by the previous
sentence (and choice of $(\fra, \frb, \lambda)$, $i$ is limit,
so $\pcf(\{\lambda_j\colon j<i\})\not
\subseteq\lambda$ but it is $\subseteq
\pcf(\frb)\subseteq\lambda^+$, so $\lambda=\max\pcf\{
\lambda_j\colon j<i\}$.
For each $j$, by the minimality condition for some
$\frb_j\subseteq\frb$, we have $|\frb_j|\le |\fra|$,
$\lambda_j\in\pcf(\frb_j)$.
So $\lambda\in\pcf\{\lambda_j\colon j<i\}\subseteq\pcf
(\bigcup_{j<i}\frb_j)$ but $\bigcup_{j<i} \frb_j$ is a subset of
$\frb$ of cardinality $\le |i|\times |\fra| =|\fra|$.

\demo{6.7D Proof of 6.7A} 
Let $\langle\langle
f_\alpha^{\fra,\lambda}\colon\alpha<\lambda\rangle\colon\lambda\in
\pcf\fra\rangle$ be chosen as in the proof of 6.7.
For  $\zeta < \kappa $  we define  $\fra^\zeta =: N_{\zeta}  
\cap  \pcf\fra$;  we also define  $^\zeta \bar f$  as\break  
$\langle \langle f^{\fra^\zeta ,\lambda }_{\alpha} \colon \alpha <
\lambda \rangle \colon \lambda \in
\pcf \fra\rangle $  where  $f^{\fra^\zeta ,\lambda }_{\alpha} 
\in  \prod \fra^\zeta $ is defined as follows:

\item{(a)} if  $\theta \in  \fra$,  $f^{\fra^\zeta ,\lambda }_{\alpha}
(\theta ) = f^{\fra,\lambda }_{\alpha} (\theta )$,

\item{(b)} if  $\theta \in  \fra^\zeta \backslash \fra$ and 
$\cf(\alpha) \notin  (|\fra^\zeta|,\Min\fra$),  then 
$$
f^{\fra^\zeta ,\lambda }_{\alpha} (\theta ) = 
\Min \left\lbrace \gamma < \theta \colon f^{\fra,\lambda }_{\alpha}
\upharpoonright \frb_{\theta} [\fra] \leq_{J_{<\theta }[\frb_{\theta}
[\fra]]}f^{\fra,\theta }_{\gamma}
\upharpoonright \frb_{\theta} [\fra] \right\rbrace,
$$

\item{(c)} if $\theta\in\fra^\zeta\backslash\fra$ and  $\cf(\alpha)\in 
(|\fra^\zeta\vert,\Min\fra$),  define  $f^{\fra^\zeta ,\lambda }_{\alpha}
(\theta )$  so as to satisfy $(\ast )_{1}$ in the proof of 6.7.

Now  $^\zeta \bar f$  is legitimate except that we have only
$$
\beta < \gamma <
\lambda \in  \pcf \fra \Rightarrow  f^{\fra^\zeta ,\lambda }_{\beta}
\leq  f^{\fra^\zeta ,\lambda }_{\gamma} \mod J_{<\lambda}
[\fra^\zeta ]
$$
(instead of strict inequality) and $ \bigwedge_{\beta <
\lambda } \bigvee_{\gamma<\lambda} \left[ f^{\fra^\zeta,
\lambda }_{\beta}  <  f^{\fra^\zeta ,\lambda }_{\gamma} \mod
J_{<\lambda} [\fra^\zeta ] \right] $,  but this suffices.
(The first statement is actually proved in [Sh371, 3.2A], the
second in [Sh371, 3.2B]; by it also $^\zeta \bar f$ is cofinal in the
required sense.)

For every  $\zeta < \kappa $  we can apply 6.7 with  $(N_{\zeta}  \cap  
\pcf\fra)$, $^\zeta\bar f$  and $\langle N_{\zeta +1+i}\colon i <
\kappa \rangle $
here standing for  $\fra$, $\bar f$, $\bar N$  there. 
In the proof of 6.7 get a club $E_{\zeta} $ of  $\kappa $  (so any 
$i < j$  from  $E_{\zeta} $ are O.K.). 
Now we can define for $\zeta<\kappa$ and $i < j$ in $E_{\zeta}$,  
$^\zeta \frb^{i,j}_{\lambda} $ and 
$\langle ^\zeta \frb^{i,j,\epsilon }_{\lambda} \colon \epsilon 
< |\fra^\zeta|^+ \rangle $,  
$\langle\epsilon^\zeta
(i,j,\lambda)\colon\lambda\in\fra^\zeta\rangle$,
$\epsilon^\zeta(i,j)$, as well as in the proof of 6.7.  Let:   
$$
E = \left\lbrace i < \kappa \colon i\right.  \quad\hbox{ is a limit
ordinal }
(\forall j < i)(j + j < i \und j \times  j < i) \ \hbox{  and }
\bigwedge_{{j} < {i}}i \in  E_{j}\left. \right\rbrace .
$$
So by [Sh420, {\S}1] we can find $\bar C = \langle C_{\delta} \colon
\delta
\in  S\rangle $,  $S \subseteq  \{\delta < \kappa \colon \cf \delta =
\cf \sigma \}$  
stationary,  $C_{\delta} $ a club of  $\delta $, $\otp C_{\delta}  = 
\omega ^2\sigma $  such that:
\item{(1)}
for each  $\alpha < \lambda $,  $\{C_{\delta}  \cap  \alpha \colon
\alpha \in 
\nacc(C_{\delta})\}$  has cardinality  $< \kappa $,\footnote{*}{If
$\kappa $  is successor of regular, then we
can get $[\gamma \in  C_{\alpha}  \cap  C_{\beta}  \Rightarrow
C_{\alpha} \cap  \gamma = C_{\beta}  \cap  \gamma]$.}
and

\item{(2)} for every
club $E'$ of  $\theta $  for stationarily many  $\delta \in  S$, 
$C_{\delta}  \subseteq E'$. 

\noindent
Without loss of generality  $\bar C \in  N_{0}$. 
For some  $\delta ^\ast $,  $C_{\delta ^\ast } \subseteq  E$,  and let  
$\{j_{\zeta} \colon \zeta \leq  \omega ^2\sigma \}$  enumerate 
$C_{\delta ^\ast } \cup  \{\delta ^\ast \}$. 
So  $\langle j_{\zeta} \colon \zeta \leq  
\omega ^2\sigma \rangle $  is a strictly increasing continuous
sequence of ordinals from  $E \subseteq  \kappa $  such that  
$\langle j_{\epsilon}\colon\epsilon\leq\zeta\rangle\in N_{j_{\zeta+1}}$.
Let  $j(\zeta ) = j_{\zeta} $,  $i(\zeta ) = i_{\zeta}  =: 
j_{\omega ^2(1+\zeta )}$,  $\fra_{\zeta}  = N_{i_{\zeta} }\cap\pcf\fra$,
and $\bar \fra =: \langle \fra_{\zeta} \colon \zeta < \sigma \rangle $, 
$\frb^\zeta_{\lambda} [\bar \fra] =:$\break 
${}^{i(\zeta )}\frb_\lambda^{j(\omega^2\zeta +1),j(\omega^2\zeta +2),
\epsilon^\zeta
(j(\omega^2\zeta +1), j(\omega^2\zeta +2))}$.  
Most of the requirements follow\break
immediately, as 

\item{$(\ast )$}  for  each $\zeta<\sigma$, we have $\fra_{\zeta} $, 
$\langle \frb^\zeta_{\lambda}
[\bar \fra]\colon \lambda \in  \fra_{\zeta} \rangle $  are as in 6.7
and belong to  $N_{i_{\beta} +3} \subseteq  N_{i_{\beta +1}}.$

\noindent
We are left (for proving 6.7A) with proving (h)$^+$ and (i) (remember
(h) is a special case of (h)$^+$ choosing $\theta=\aleph_0)$. 

For proving clause (i) note that for  $\zeta < \xi  < \kappa $,  
$f^{\fra^\zeta ,\lambda }_{\alpha}  \subseteq  f^{\fra^\xi ,
\lambda}_{\alpha} $ hence  
$^\zeta \frb^{i,j}_{\lambda}  \subseteq\  ^\xi \frb^{i,j}_{\lambda} $. 
Now we can prove by induction on  $\epsilon $  that 
$^\zeta \frb^{i,j, \epsilon }_{\lambda}  \subseteq\ 
^\xi\! \frb^{i,j,\epsilon}_{\lambda}$ for every $\lambda\in\fra_\zeta$
(check the definition after $(*)_2$ in the proof of 6.7) and the
conclusion follows.

Instead of proving (h)$^+$ we prove an apparently weaker version
(h)$'$ below, and then note that $\bar i'=\langle
i_{\omega^2\zeta} \colon\zeta<\sigma\rangle$, $\bar\fra' =\langle
\fra_{\omega^2\zeta}\colon \zeta<\sigma\rangle$, $\langle
N_{i(\omega^2\zeta)}\colon \zeta<\sigma\rangle$,
$\langle \frb_\lambda^{\omega^2\zeta} [\bar\fra']\colon
\zeta<\sigma, \lambda\in \fra'_\zeta
=\fra_{\omega^2\zeta}\rangle$ will exemplify the
conclusion\footnote{**}{Assuming $\sigma>\aleph_0$ hence,
$\omega^2\sigma= \sigma$ for notational simplicity.}
where

\item{(h)$'$} if $\frc\subseteq \fra_\beta$, $\beta<\sigma$,
$\frc\in N_{i_{\beta+1}}$, $\theta=\cf(\theta)\in N_{i_{\beta+1}}$
{\sl then} for some $\frd\in N_{i_{\beta+\omega+1}+1}$,
$\frd\subseteq \fra_{\beta+\omega}
\cap\pcf_{\theta{\scriptstyle{\rm -complete}}}(\frc)$ we
have $\frc\subseteq \bigcup_{\mu\in\frd}\frb_\mu^{\beta+\omega}
[\bar\fra]$ and $|\frd|<\theta$.

\demo{Proof of (h)$^\prime$}
So let  $\theta $, $\beta $, $\frc$ 
be given; let $\langle\frb_{\mu}[\fra]\colon\mu\in\pcf\frc\rangle
(\in N_{i_{\beta+1}})$ be a generating sequence. 
We define by induction on  $n < \omega $,  $A_{n}$,  
$\langle \frc_{\eta} ,\lambda_{\eta} \colon \eta \in  A_{n}\rangle $ 
such that: 

\item{(a)} $A_0=\{\langle\rangle\}$, $\frc_{\langle \rangle} = \frc$, 
$\lambda_{\langle \rangle} = \max \pcf \frc$,

\item{(b)} $A_{n} \subseteq\  ^n\!\theta $, $|A_{n}| < \theta $,

\item{(c)} if  $\eta\in A_{n+1}$ then $\eta\upharpoonright n\in A_{n}$, 
$\frc_{\eta}\subseteq\frc_{\eta\upharpoonright n}$, $\lambda_{\eta} <
\lambda_{\eta \upharpoonright n}$ and $\lambda_{\eta}  = \max  \pcf
(\frc_{\eta})$,

\item{(d)} $A_{n}$, $\langle \frc_{\eta} ,\lambda_{\eta} \colon \eta \in
A_{n}\rangle $  belongs to  $N_{i_{\beta +1+n}}$ hence \ $\lambda_{\eta}
\in  N_{i_{\beta +1+n}}$, 

\item{(e)} if  $\eta \in  A_{n}$ and \ $\lambda_{\eta}  \in  
\pcf_{\theta{\scriptstyle{\rm{-complete}}}}(\frc_{\eta} )$ and  $\frc_{\eta} {\miss}
\frb^{\beta +1+n}_{\lambda_{\eta} }[\bar\fra]$  then\break 
$(\forall \nu )[\nu  \in 
A_{n+1}\shvor \und \eta \subseteq \nu  \Leftrightarrow  \nu  = \eta
\char94 \langle 0\rangle ]$  and 
$\frc_{\eta \char94 \langle 0\rangle } =
\frc_{\eta} \backslash \frb^{\beta +1+n}_{\lambda_\eta}[\bar \fra]$  (so
$\lambda_{\eta \char94 \langle 0\rangle } =$\break
$\max  \pcf \frc_{\eta \char94
\langle 0\rangle } < \lambda_{\eta}  = \max  \pcf \frc_{\eta} )$,

\item{(f)} if  $\eta \in  A_{n}$ and $\lambda_{\eta}  \notin
\pcf_{\theta{\scriptstyle{\rm{-complete}}}}(\frc_\eta)$
 then 
$$
\frc_{\eta} =
\bigcup \left\lbrace \frb_{\lambda_{\gamma \char94 \langle i\rangle
}}[\frc]\colon i < i_n < 
\theta ,\eta \char94 \langle i\rangle  \in A_{n+1}\right\rbrace,
$$ 
and if  $\nu  = 
\eta \char94 \langle i\rangle  \in  A_{n+1}$ then  $\frc_{\nu}  = 
\frb_{\lambda_{\nu} }[\frc]$,

\item{(g)} if $\eta\in A_n$, and
$\lambda_\eta\in\pcf_{\theta{\scriptstyle{\rm{-complete}}}}
(\frc_\eta)$ but
$\frc_\eta\subseteq \frb_{\lambda_n}^{\beta+1-n}[\bar\fra]$, then
$\lnot (\exists\nu)[\eta\triangleleft\nu\in A_{n+1}]$.

\noindent
There is no problem to carry the definition (we use 6.7F(1)
below\footnote*{No vicious circle; 6.7F(1) does not depend on
6.7B.}, the
point is that  $\frc \in  N_{i_{\beta +1+n}}$ implies  
$\langle \frb_{\lambda} [\frc]\colon \lambda \in  \pcf_{\theta}
[\frc]\rangle  \in  
N_{i_{\beta +1+n}}$ and as there is $\frd$ as in 6.7F(1),
there is one in  $N_{i_{\beta +1+n+1}}$ so 
$\frd \subseteq  \fra_{\beta +1+n+1}$). 
Now let
$$
\frd_n =: \left\lbrace \lambda_{\eta}^{} \colon \eta \in
A_{n}\hbox{ and }\lambda_{\eta}  \in
\pcf_{\theta{\scriptstyle{\rm{-complete}}}}(\frc_\eta)\hbox{ and }
\frc_\eta\subseteq \frb^{\beta+1+n}_{\lambda_\eta}[\fra]\right\rbrace
$$
and
$\frd =: \bigcup_{n<\omega} \frd_n$; 
we shall show that it is as required.          

The main point is  $\frc \subseteq \bigcup_{{\lambda} \in  {\frd}}
\frb^{\beta +\omega }_{\lambda} [\bar \fra]$; note that 
$$
\left[ \lambda_{\eta}  \in  \frd,\eta \in  A_n  \Rightarrow 
\frb^{\beta +1+n}_{\lambda_{\eta}} [\bar \fra]
\subseteq\frb^{\beta +\omega }_{\lambda_{\eta} }
[\bar \fra]\right]
$$
hence it suffices to show  $\frc \subseteq \bigcup_{n<\omega}
\bigcup_{{\lambda} \in  {\frd_n}}\frb^{\beta + 1 + n}_{\lambda} [\bar
\fra]$, so assume  $\theta \in$\break  
$\frc\backslash\bigcup_{n<\omega} \bigcup_{{\lambda} \in 
{\frd_n}}\frb^{\beta
+ 1 + n}_{\lambda} [\bar \fra]$,  and we choose by induction on  $n$,  
$\eta_{n} \in  A_{n}$ such that  $\eta_{0} = <>$, $\eta_{n+1}
\upharpoonright n
= \eta_{n}$ and $\theta \in  \frc_{\eta} $;  by clauses (e) $+$ (f)
above
this is possible and $\langle \max  \pcf \frc_{\eta_{n}}\colon
n <
\omega \rangle $  is strictly decreasing, contradiction. 

The minor point is  $|\frd| < \theta $;  if  $\theta > \aleph_{0}$ note
that  $ \bigwedge_{n} |A_{n}| < \theta $  and $\theta = \cf
(\theta)$  so  $|\frd| \leq  | \bigcup_{n} A_{n}| < \theta +
\aleph_{1} = \theta .$ 

If  $\theta = \aleph_{0}$ (i.e.\ clause (h)) we should have 
$\bigcup_{n} A_{n}$ finite; the proof is as above noting the
clause (f) is vacuous now.
So $\bigwedge_n |A_n| =1$ and $\bigvee_n A_n=\emptyset$, so $\bigcup_n
A_n$ is finite.
Another minor point is $\frd\in N_{i_{\beta+\omega+1}}$; this
holds as the construction is unique from
$\langle N_j\colon j< i_{\beta+\omega}\rangle$,
$\langle i_j\colon j\le \beta+\omega\rangle$, $\langle
(\fra_{i(\zeta)},\langle \frb_\lambda^\zeta\colon\lambda\in
\fra_{i(\zeta)}\rangle)\colon \zeta\le\beta+\omega\rangle$; no
``outside'' information is used so $\langle (A_n, \langle
(c_\eta, \lambda_\eta)\colon\eta\in A_n\rangle)\colon
n<\omega\rangle\in N_{i_{\beta+\omega+1}}$,
so (using a choice function) really $\frd \in
N_{i_{\beta+\omega+1}}$.
\sqed{6.7{\rm A}}

\demo{6.7E Proof of 6.7B} 
Let  $\frb_\lambda [\bar \fra] =
\frb_\lambda^\sigma=\bigcup_{\beta<\sigma}
\frb^\beta_\lambda [\fra_\beta]$ and
$\fra_\sigma=\bigcup_{\zeta<\sigma} \fra_\zeta$. 
Part (1) is straightforward.
For part (2), for clause (g), for $\beta=\sigma$, the inclusion
``$\subseteq$'' is straightforward;
so assume $\mu\in
\fra_\beta\cap\pcf\frb^\beta_\lambda[\bar\fra]$.
Then by 6.7A(c) for some $\beta_0<\beta$, we have
$\mu\in\fra_{\beta_0}$, and by 6.7C(3B) (which
depends on 6.7A only) for some $\beta_1<\beta$,
$\mu\in\pcf\frb^{\beta_1}_\lambda[\bar\fra]$;
by monotonicity wlog $\beta_0=\beta_1$, by clause (g) of 6.7A
applied to $\beta_0$, $\mu\in\frb^{\beta_0}_\lambda[\bar\fra]$.
Hence by clause (i) of 6.7A,
$\mu\in\frb_\lambda^\beta[\bar\fra]$, thus proving the other
inclusion.

The proof of clause (e) (for 6.7B(2)) is similar, and also
6.7B(3).
For 6.7(B)(4) for $\delta<\sigma$, $\cf(\delta)>|\fra|$ redefine
$\frb_\lambda^\delta [\bar a]$ as 
$\bigcup_{\beta<\delta} \frb_\lambda^{\beta+1}[\fra]$.
\sqed{6.7{\rm B}}

\proclaim{6.7F Claim} 
Let  $\theta $  be\ regular.
\item{(0)}  If  $\alpha < \theta $, 
$\pcf_{\theta{\scriptstyle\rm{-complete}}}
\left( \bigcup_{{i} < {\alpha}} \fra_{i}\right)  = \bigcup_{{i} <
{\alpha}}  \pcf_{\theta{\scriptstyle{\rm{-complete}}}}(\fra_{i}).$

\noindent
\item{(1)} If  $\langle \frb_{\theta} [\fra]\colon \theta \in  \pcf \fra
\rangle $  is a generating 
sequence for $\fra$, $\frc\subseteq\fra$, then for some $\frd\subseteq  
\pcf_{\theta{\scriptstyle{\rm{-complete}}}}(\frc)$ we have: $|\frd|<\theta$ and $\frc 
\subseteq \bigcup_{{\theta} \in  {\fra}}\frb_{\theta} [\fra].$

\item{(2)}  If  $|\fra \cup  \frc| < \Min\fra$,  $\frc \subseteq  
\pcf_{\theta{\scriptstyle{\rm{-complete}}}}(\fra)$,  $\lambda \in  
\pcf_{\theta{\scriptstyle{\rm{-complete}}}}(\frc)$  then 
$\lambda \in  
\pcf_{\theta{\scriptstyle{\rm{-complete}}}}(\fra).$

\item{(3)}  In (2) we can weaken  $|\fra \cup  \frc| < \Min\fra$ to 
$|\fra| < \Min\fra$,  $|\frc| < \Min\frc.$

\item{(4)}  We cannot find $\lambda_{\alpha}  \in  \pcf_{\theta
{\scriptstyle{\rm{-complete}}}}(\fra)$  
for  $\alpha < |\fra|^+$ such that  $\lambda_{i} >$\break
$\sup \pcf_{\theta{\scriptstyle{\rm{-complete}}}}(\{\lambda_{j}\colon j < i\}).$

\item{(5)}  Assume $\theta\leq |\fra|$, $\frc \subseteq
\pcf_{\theta{\scriptstyle{\rm{-complete}}}}\fra$ 
(and $|\frc| < \Min \ \frc$;  of course  $|\fra| < \Min\fra)$. 
If  $\lambda \in 
\pcf_{\theta{\scriptstyle{\rm{-complete}}}}(\frc)$
then for some  $\frd \subseteq  \frc$  we
have  $|\frd| \leq  |\fra|$  and $\lambda \in  \pcf_{\theta
{\scriptstyle{\rm{-complete}}}}
(\frd).$
\endproclaim

\demo{Proof} 
(0) and (1): Check.

\noindent
(2)  See [Sh345b, 1.10--1.12].

\noindent
(3)  Similarly.

\noindent
(4)  If  $\theta = \aleph_{0}$ we already know it (e.g.\
6.7C(3A)), so assume  $\theta > 
\aleph_{0}$ and, without loss of generality,  $\theta $  is regular
$\le|\fra|$. 
We use 6.7A with  $\{\theta$, $\langle \lambda_{i}\colon i < |\fra|^+
\rangle  \}\in  N_{0}$, $\sigma = |\fra|^+$,  $\kappa = |\fra|^{+3}$
where, without loss of generality,
$\kappa < \Min(\fra)$.
For each  $\alpha < |\fra|^+$ by (h)$^+$ of 6.7A there is  
$\frd_{\alpha}  \in  N_{i_{1}}$,  $\frd_{\alpha}  \subseteq  
\pcf_{\theta{\scriptstyle{\rm{-complete}}}}(\{\lambda_{i}\colon i < \alpha \})$, 
$|\frd_{\alpha} | < 
\theta $  such that  $\{\lambda_{i}\colon i < \alpha \} \subseteq
\bigcup_{{\theta} \in  {\frd_{\alpha}} }\frb^1_{\theta} [\bar \fra]$; 
hence by clause (g) of 6.7A and 6.7F(0) we have 
$\fra_{1} \cap  \pcf_{\theta{\scriptstyle{\rm{-complete}}}}(\{\lambda_{i}\colon 
i < \alpha \}) \subseteq \bigcup_{{\theta} \in  {\frd_{\alpha}} }
\frb^1_{\theta} [\bar \fra]$. 
So for $\alpha < \beta < |\fra|^+$,  $\frd_{\alpha}  \subseteq 
\fra_1\cap\pcf_{\theta{\scriptstyle{\rm{-complete}}}}\{\lambda_{i}\colon i < \alpha \} 
\subseteq \fra_1\cap \pcf_{\theta{\scriptstyle{\rm{-complete}}}}\{\lambda_{i}\colon i < \beta \}
\subseteq  
 \bigcup_{{\theta} \in  {\frd_{\beta}} }\frb^1_{\theta} [\bar \fra]$. 
As the sequence is smooth (i.e.\ clause (f) of 6.7A) clearly 
$\alpha < \beta \Rightarrow \bigcup_{{\mu} \in  {\frd_{\alpha}} }
\frb^1_{\mu} [\bar \fra] \subseteq \bigcup_{{\mu} \in  {\frd_{\beta}} }
\frb^1_{\mu} [\bar \fra].$ 

So  $\langle \bigcup_{{\mu} \in  {\frd_{\alpha}} }\frb^1_{\mu}
[\bar \fra] \cap  \fra\colon \alpha < |\fra|^+\rangle $  is a 
non-decreasing sequence of subsets of $\fra$ of length  $|\fra|^+$,
 hence for some  $\alpha (\ast ) < |\fra|^+$ we have: 

\item{$(\ast )_{1}$}  $\alpha (\ast ) \leq  \alpha < |\fra|^+
\Rightarrow \bigcup_{{\mu} \in  {\frd_{\alpha}} }\frb^1_{\mu}
[\bar \fra] \cap  \fra = \bigcup_{{\mu} \in  {\frd_{\alpha (\ast )}}}
\frb^1_{\mu} [\bar \fra] \cap  \fra.$ 

If $\tau\in\fra_{1} \cap  \pcf_{\theta{\scriptstyle{\rm{-complete}}}}(\{\lambda_{i}
\colon i<\alpha\})$ then
$\tau\in\pcf_{\theta{\scriptstyle{\rm{-complete}}}}(\fra)$ (by 6.7F(2),(3)),
and
$\tau\in\frb^1_{\mu_{\tau} }[\bar \fra]$ 
for some $\mu_{\tau}\in\frd_{\alpha} $ so  $\frb^1_{\tau} [\bar \fra]
\subseteq\frb^1_{\mu_{\tau}}[\bar\fra]$, also $\tau\in\pcf\nolimits_{\theta
{\scriptstyle{\rm{-complete}}}}(\frb^1_{\tau} [\bar \fra] \cap  \fra)$  
(by clause (e) of 6.7A),  hence  
$$
\eqalign{
\tau\in \pcf\nolimits_{\theta{\scriptstyle{\rm{-complete}}}}(\frb^1_{\tau} [\bar \fra] \cap 
\fra)
& \subseteq \pcf\nolimits_{\theta{\scriptstyle{\rm{-complete}}}}(\frb^1_{\mu_{\tau} }
[\bar \fra] \cap \fra)\cr
& \subseteq  \pcf\nolimits_{\theta{\scriptstyle{\rm{-complete}}}}\left( \bigcup_{{\mu}
\in  {\frd_{\alpha}} }\frb^1_{\mu} [\bar \fra] \cap  \fra\right) .
\cr}
$$

\noindent
So  $\fra_{1} \cap  \pcf_{\theta{\scriptstyle{\rm{-complete}}}}(\{\lambda_{i}\colon
i <
\alpha \}) 
\subseteq  \pcf_{\theta{\scriptstyle{\rm{-complete}}}}\left( \bigcup_{{\mu} \in 
{\frd_{\alpha}} }\frb^1_{\mu} [\bar \fra] \cap  \fra\right) $. 
But for each  $\alpha < |\fra|^+$ we have 
$\lambda_{\alpha}  > \sup  \pcf_{\theta{\scriptstyle{\rm{-complete}}}}(\{\lambda_{i}
\colon i < \alpha \})$, whereas $\frd_\alpha\subseteq$\break
$\pcf_{\sigma{\scriptstyle{\rm{-complete}}}}\{\lambda_i\colon i<\alpha\}$,
hence  $\lambda_{\alpha} >\sup \frd_{\alpha}$ hence 

\item{$(\ast )_{2}$}  $\lambda_{\alpha}  >  \sup_{\mu
\in \frd_\alpha} \max  \pcf \frb^1_{\mu} [\bar \fra] \geq  \sup  
\pcf_{\theta{\scriptstyle{\rm{-complete}}}}\left( \bigcup_{{\mu} \in 
{\frd_{\alpha}} }\frb^1_{\mu} [\bar \fra] \cap  \fra\right) .$ 

On the other hand, 

\item{$(\ast )_{3}$} $\lambda_{\alpha}  \in  \pcf_{\theta
{\scriptstyle{\rm{-complete}}}}\{\lambda_{i}\colon i < 
\alpha + 1\} \subseteq  \pcf_{\theta{\scriptstyle{\rm{-complete}}}}\left(
\bigcup_{{\mu} \in  {\frd_{\alpha +1}}}\frb^1_{\mu} [\bar \fra] \cap  \fra\right) .$

\noindent
For  $\alpha = \alpha (\ast )$  we get contradiction by $(\ast )_{1}
+ (\ast )_{2} + (\ast )_{3}.$

\noindent
(5)  Assume $\fra$, $\frc$, $\lambda $  form a counterexample with
$\lambda$ minimal.
Without loss of generality  $|\fra|^{+3} < \Min(\fra)$ and
$\lambda = \max
\pcf\fra$ and $\lambda = \max  \pcf \frc$  (just let  $\fra' =:
\frb_{\lambda} [\fra]$,  $\frc' =: \frc \cap  \pcf_{\theta} [\fra']$; 
if  $\lambda \notin  \pcf_{\theta{\scriptstyle{\rm{-complete}}}}(\frc')$  
then necessarily  $\lambda \in  \pcf (\frc\backslash \frc')$ 
(by 6.7F(0)) and similarly  $\frc\backslash \frc' \subseteq  
\pcf_{\theta{\scriptstyle{\rm{-complete}}}}(\fra\backslash \fra')$  hence by 6.7F(2),(3) 
$\lambda\in\pcf_{\theta{\scriptstyle{\rm{-complete}}}}
(\fra\backslash \fra')$,\break 
contradiction).

Also without loss of generality  $\lambda \notin  \frc$. 
Let  $\kappa $, $\sigma $, $\bar N$,  $\langle i_{\alpha}=i(\alpha)
\colon\alpha\leq\sigma\rangle$, $\bar\fra=\langle\fra_i \colon i
\leq  \sigma\rangle$ be as in 6.7A with $\fra\in N_{0}$, $\frc\in
N_{0}$,  $\lambda \in 
N_{0}$,  $\sigma = |\fra|^+$,  $\kappa = |\fra|^{+3} < \Min\fra$. 
We choose by induction on  $\epsilon  
< |\fra|^+$,  $\lambda_{\epsilon} $,  $\frd_{\epsilon} $ such that: 

\item{(a)}
$\lambda_{\epsilon}  \in  \fra_{\omega ^2\epsilon +\omega +3}$,
$\frd_{\epsilon}  \in  N_{i(\omega ^2\epsilon +\omega +1)}$, 

\item{(b)} $\lambda_{\epsilon}  \in  \frc$,

\item{(c)} $\frd_{\epsilon}\subseteq\fra_{\omega^2\epsilon +\omega+1}
\cap  \pcf_{\theta{\scriptstyle{\rm{-complete}}}}(\{\lambda_{\zeta}
\colon \zeta < \epsilon \})$, 

\item{(d)} $|\frd_{\epsilon} | < \theta $, 

\item{(e)} $\{\lambda_{\zeta} \colon \zeta < \epsilon \} \subseteq
\bigcup_{{\theta} \in  {\frd_{\epsilon}} }
\frb^{\omega ^2\epsilon +\omega +1}_{\theta} [\bar \fra]$,

\item{(f)} $\lambda_{\epsilon}  \notin 
\pcf_{\theta{\scriptstyle{\rm{-complete}}}}
\left( \bigcup_{{\theta} \in  {\frd_{\epsilon}} }\frb^{\omega ^2
\epsilon +\omega +1}_{\theta} [\bar \fra]\right) .$

\noindent
For every  $\epsilon  < |\fra|^+$ we first choose  $\frd_{\epsilon} $
as the $<^\ast_{\chi} $-first element satisfying (c) $+$ (d) $+$ (e) and
then if
possible $\lambda_{\epsilon} $ as the $<^\ast_{\chi} $-first element
satisfying (b) $+$ (f). 
It is easy to check the requirements and in fact  
$\langle \lambda_{\zeta} \colon \zeta < \epsilon \rangle  \in  
N_{\omega ^2\epsilon +1}$,\break  $\langle \frd_{\zeta} \colon \zeta
< \epsilon \rangle  \in  N_{\omega ^2\epsilon +1}$ 
(so clause (a) will hold). 
But why can we choose at all?  Now  $\lambda \notin  
\pcf_{\theta{\scriptstyle{\rm{-complete}}}}\{\lambda_{\zeta} \colon \zeta
< \epsilon \}$ as $\fra$, $\frc$, $\lambda $  form a counterexample with
$\lambda$ minimal
and $\epsilon  < |\fra|^+$ (by 6.7F(3)).
As  $\lambda = \max  \pcf\fra$  necessarily  
$\pcf_{\theta{\scriptstyle{\rm{-complete}}}}(\{\lambda_{\zeta}
\colon \zeta < \epsilon \}) \subseteq  \lambda $  hence 
$\frd_{\epsilon}  \subseteq  \lambda $ (by clause (c)). 
By part (0) of the claim (and clause (a)) we know:  
$$
\eqalign{
\pcf\nolimits_{\theta{\scriptstyle{\rm{-complete}}}}\left[ \bigcup_{{\mu} \in 
{\frd_{\epsilon}} }\frb^{\omega ^2\epsilon +\omega +1}_{\mu} [\bar \fra]\right] &
=  \bigcup_{{\mu} \in  {\frd_{\epsilon}} }\pcf\nolimits_{\theta
{\scriptstyle{\rm{-complete}}}}\left[ \frb^{\omega ^2+\omega +1}_{\mu}
[\bar \fra]\right] \cr
& \subseteq \bigcup_{{\mu} \in  {\frd_{\epsilon}} }(\mu  + 1)
\subseteq  \lambda \cr}
$$
(note  $\mu  = \max  \pcf \frb^\beta_{\mu} [\bar \fra])$. 
So  $\lambda \notin  \pcf_{\theta{\scriptstyle{\rm{-complete}}}}\left(
\bigcup_{{\mu} \in  {\frd_{\epsilon}} }\frb^{\omega ^2\epsilon +\omega
+1}_{\mu} [\bar \fra]\right) $  hence by 
part (0) of the claim  $\frc {\miss} \bigcup_{{\mu} \in
{\frd_{\epsilon}}}\frb^{\omega ^2\epsilon +\omega +1}_{\mu}
[\bar \fra]$  so  $\lambda_{\epsilon} $ exists. 
Now  $\frd_{\epsilon} $ exists by 6.7A clause (h)$^+$.

Now clearly  $\left\langle\fra \cap \bigcup_{{\mu} \in 
{\frd_{\epsilon}} }\frb^{\omega ^2\epsilon +\omega +1}_{\mu}
[\bar \fra]\colon \epsilon < |\fra|^+\right\rangle$ is non-decreasing
(as in the earlier proof) hence eventually
constant, say for  $\epsilon \geq  \epsilon (\ast )$  (where 
$\epsilon (\ast ) < |\fra|^+)$. 
But
\item{$(\alpha )$} $\lambda_{\epsilon}  \in \bigcup_{{\mu} \in 
{\frd_{\epsilon +1}}}\frb^{\omega ^2\epsilon +\omega +1}_{\mu}
[\bar \fra]$ [clause (e) in the choice of $\lambda_\epsilon,
\frd_\epsilon$], 

\item{$(\beta )$} $\frb^{\omega ^2\epsilon +\omega
+1}_{\lambda_{\epsilon} }[\bar \fra] 
\subseteq \bigcup_{{\mu} \in  {\frd_{\epsilon +1}}}\frb^{\omega ^2
\epsilon +\omega +1}_{\mu} [\bar \fra]$  [by clause (f) of 6.7A
and $(\alpha)$ alone],

\item{$(\gamma )$} $\lambda_{\epsilon}  \in  \pcf_{\theta
{\scriptstyle{\rm{-complete}}}}(\fra)$  [as  
$\lambda_{\epsilon}  \in  \frc$  and a hypothesis], 

\item{$(\delta )$} $\lambda_{\epsilon}  \in  \pcf_{\theta
{\scriptstyle{\rm{-complete}}}}(\frb^{\omega ^2\epsilon +
\omega +1}_{\lambda_{\epsilon} }
[\bar \fra])$  [by $(\gamma )$ above and clause (e) of 6.7A], 

\item{$(\epsilon)$} $\lambda_\epsilon\not\in \pcf(\fra\hefresh
\frb^{\omega^2\epsilon+\omega+1}_{\lambda_\epsilon})$,

\item{$(\zeta )$} $\lambda_{\epsilon}  \in  
\pcf_{\theta{\scriptstyle{\rm{-complete}}}}\left(\fra\cap\bigcup_{{\mu} \in 
{\frd_{\epsilon +1}}}\frb^{\omega ^2\epsilon +\omega
+1}_{\mu}[\bar \fra]\right)$ [by $(\delta ) + (\epsilon) +(\beta )].$

\noindent
But for  $\epsilon  = \epsilon (\ast )$,  the statement $(\zeta )$ 
contradicts the choice of $\epsilon(*)$ and clause (f) above.
\sqed{6.7{\rm F}}


%
%
\references{letterkey}{ARSh153}

\refpaper
\key ARSh153
\by U. Abraham, M. Rubin and S. Shelah
\reftitle On the consistency of some 
partition theorems for continuous colorings and the structure of 
$\aleph_{1}$-dense real order type
\journal Annals of Pure and Applied Logic
\volume 29
\yr 1985
\pages 123--206
\endref

\refpaper
\key BoSh210
\by R. Bonnet and S. Shelah
\reftitle Narrow Boolean algebras
\journal Annals of Pure and Applied Logic
\volume  28
\yr 1985
\pages 1--12
\endref

\refpaper
\key CEG
\by  F. S. Cater, Paul Erd\"os and Fred Galvin
\reftitle On the density of $\lambda$-box products
\journal General Topology and its Applications
\volume 9
\yr 1978
\pages 307--312
\endref

\refpaper
\key CN1
\by  W. W. Comfort and S. Negrepontis
\reftitle On families of large oscillation
\journal Fundamenta Mathematicae
\volume 75
\yr 1972
\pages 275--290
\endref

\refbook
\key CN2
\by  W. W. Comfort and S. Negrepontis
\reftitle The 
Theory of Ultrafilters, Grundleh\-ren der Mathematischen
Wissenschaften {\bf 211}
\publisher Springer-Verlag, Berlin--\break Heidelberg--New York 
\yr 1974
\endref

\refpaper
\key CR
\by W. W. Comfort and Lewis C. Robertson
\reftitle Cardinality constraints for pseudocompact and for totally
dense subgroups of compact topological groups
\journal \PJM
\volume 119
\yr 1985
\pages 265--285
\endref

\refpaper
\key DoJe
\by  T. Dodd and R. Jensen
\reftitle The Covering Lemma for  $K$
\journal \AML
\volume 22
\yr 1982
\pages 1--30
\endref

\refpaper
\key EK
\by  R. Engelking and M. Karlowicz
\reftitle Some theorems of set theory and their topological consequences
\journal Fundamenta Math
\volume 57
\yr 1965
\pages 275--285
\endref

\refpaper
\key H
\by Edwin Hewitt
\reftitle A remark on density characters
\journal \BAMS
\volume 52
\yr 1946
\pages 641--643
\endref

\refpaper
\key Ma
\by Edward Marczewski
\reftitle S\'eparabilit\'e et multiplication cart\'esienne 
des espaces topologiques
\journal Fundamenta Mathematicae
\volume  34
\yr 1947
\pages 127--143
\endref

\refall
\key Mo
\by D. Monk
\reftitle Cardinal functions of Boolean algebras
\rest circulated notes.
\endref

\refpaper
\key P
\by E. S. Pondiczery
\reftitle Power problems in abstract spaces
\journal \DMJ
\volume 11
\yr 1944
\pages 835--837
\endref

\refall
\key Sh-b
\by  S. Shelah
\reftitle Proper Forcing
\rest Lectures Notes in Mathematics {\bf 940}, Springer-Verlag,
Berlin--New York, 1982.
\endref

\refall
\key Sh-f
\by S. Shelah
\reftitle Proper and Improper Forcing 
\rest Perspectives in Mathematical Logic, to appear.
\endref

\refall 
\key Sh-g
\by S. Shelah
\reftitle Cardinal Arithmetic
\rest Oxford University Press, 1994.
\endref

\refall
\key Sh88
\by  S. Shelah
\reftitle Classification theory for non-elementary classes II
\rest in {\sl Abstract Elementary Classes},
Proc. U.S.A--Israel Conference on 
Classification Theory, Chicago $12/85$ (J. T. Baldwin, ed.),
Springer-Verlag Lecture Notes {\bf 1292} (1987), 419--497.
\endref

\refpaper
\key Sh92
\by  S. Shelah
\reftitle Remarks on Boolean algebras
\journal Algebra Universalis
\volume 11
\yr 1980
\pages 77--89
\endref

\refpaper
\key Sh93
\by S. Shelah
\reftitle Simple unstable theories
\journal Annals of Mathematical Logic
\volume 19
\yr 1980
\pages 177--203
\endref

\refpaper
\key Sh107
\by  S. Shelah
\reftitle Models with second order properties IV, a general
method and eliminating diamonds
\journal \AML
\volume 25
\yr 1983
\pages 183--212
\endref

\refpaper
\key Sh111
\by S. Shelah
\reftitle On powers of singular cardinals
\journal \NDJFL
\volume 27
\yr 1986
\pages 263--299
\endref

\refpaper
\key RuSh117
\by M. Rubin and S. Shelah
\reftitle Combinatorial problems on trees: partitions,
$\Delta$-systems and large free subtrees
\journal Annals of Pure and Applied Logic
\volume 33
\yr 1987
\pages 43--81
\endref

\refpaper
\key ShLH162
\by S. Shelah, C. Laflamme and B. Hart
\reftitle Models with second order 
properties $V:$ a general principle
\journal Annals of Pure and Applied Logic
\volume64
\yr 1993
\pages 169--194
\endref

\refall
\key Sh233
\by S. Shelah
\reftitle Remarks on the number of ideals of Boolean algebra and 
open sets of a topology
\rest Springer-Verlag Lecture Notes {\bf 1182} (1986), 151--187.
\endref

\refpaper
\key Sh282
\by S. Shelah
\reftitle Successors of singulars,
cofinalities of reduced products of cardinals and 
productivity of chain conditions
\journal \IJM
\volume 62
\yr 1988
\pages 213--256
\endref

\refall
\key Sh288
\by S. Shelah
\reftitle Strong partition relations below the power set:
Consistency, was Sierpinski right, II?
\rest Proceedings of the Conference on Set Theory and its
Applications in honor of A. Hajnal and V. T. Sos, Budapest,
1/91, Colloquia Mathematica Societatis Janos Bolyai, Sets,
Graphs and Numbers, Vol. 60, 1991.
\endref

\refall
\key Sh300
\by S. Shelah
\reftitle Universal classes
\rest Chapter I--IV,
Proceedings of the  USA--Israel Conference on Classification Theory,
Chicago, December 1985 (J. Baldwin, ed.), Springer-Verlag Lecture
Notes {\bf 1292} (1987), 264--418.
\endref

\refall
\key Sh326
\by S. Shelah
\reftitle Viva la difference I: Nonisomorphism of ultrapowers of
countable models 
\rest in {\sl Set Theory of the Continuum}, Mathematical
Sciences Research Institute Publications, Vol. 26,
Springer-Verlag, Berlin, 1992, pp. 357--405.
\endref

\refpaper
\key GiSh344
\by M. Gitik and S. Shelah
\reftitle On certain indestructibility of strong cardinals and a
question of Hajnal
\journal Archive for Mathematical Logic
\volume 28
\yr 1989
\pages 35--42
\endref

\refpaper
\key Sh345
\by S. Shelah
\reftitle Products of regular cardinals and cardinal invariants of 
products of Boolean algebras
\journal \IJM
\volume 70
\yr 1990
\pages 129--187
\endref

\refall
\key Sh345a
\by S. Shelah
\reftitle Basic: Cofinalities of small reduced products
\rest in{\sl Cardinal Arithmetic}, Chapter I, Oxford University
Press, 1994.
\endref

\refall
\key Sh345b
\by S. Shelah
\reftitle Entangled orders and narrow Boolean algebras 
\rest in {\sl Cardinal Arithmetic}, Appendix 2, Oxford
University Press, 1994.
\endref

\refall
\key Sh355
\by S. Shelah
\reftitle $\aleph_{\omega +1}$ has a Jonsson algebra
\rest in {\sl Cardinal Arithmetic}, Chapter II,
Oxford University Press, 1994.
\endref

\refall
\key Sh365
\by S. Shelah
\reftitle There are Jonsson algebras in many inaccessible
cardinals
\rest in {\sl Cardinal Arithmetic}, Chapter III,
Oxford University Press, 1994.
\endref

\refall
\key Sh371
\by S. Shelah
\reftitle Advanced: cofinalities of small reduced products
\rest in {\sl Cardinal Arithmetic},
Chapter VIII, Oxford University Press, 1994.
\endref

\refall
\key Sh386
\by S. Shelah
\reftitle Bounding  $pp(\mu )$  when  $\cf (\mu ) > \mu  > 
\aleph_{0}$ using ranks and normal ideals
\rest in {\sl Cardinal Arithmetic},
Chapter V, Oxford University Press, 1994.
\endref

\refall
\key Sh400
\by S. Shelah
\reftitle Cardinal arithmetic
\rest in {\sl Cardinal Arithmetic},
Chapter IX, Oxford University Press, 1994.
\endref

\refpaper
\key Sh410
\by S. Shelah
\reftitle More on cardinal arithmetic
\journal \ARML
\volume 32
\yr 1993
\pages 399--428
\endref

\refpaper
\key GiSh412
\by  M. Gitik and S. Shelah
\reftitle More on ideals with simple forcing notions
\journal\break \APAL
\volume 59
\yr 1993
\pages 219--238
\endref

\refall
\key Sh420
\by S. Shelah
\reftitle Advances in cardinal arithmetic
\rest Proceedings of the Banff Conference in Alberta;
4/91, {\sl Finite and Infinite Combinatorics in Sets and Logic}
(N. W. Saure et al., eds.), Kluwer Academic Publishers, Dordrecht 1993,
pp. 355--383.
\endref

\refall
\key Sh460
\by S. Shelah
\reftitle The Generalized Continuum Hypothesis revisited
\rest submitted to Israel Journal of Mathematics.
\endref

\refall
\key Sh513
\by S. Shelah
\reftitle PCF and infinite free subsets
\rest  submitted to Archive for Mathematical Logic.
\endref

\refall
\key Sh589
\by S. Shelah
\reftitle PCF theory: applications
\rest in preparation.
\endref

\endreferences


\bye